\documentclass[11pt]{article}
\usepackage{amssymb}
\usepackage{amsmath}
\usepackage{latexsym}
\usepackage{longtable}

\input{prepictex}
\input{pictex}
\input{postpictex}
\begin{document}
\bibliographystyle{plain}
\baselineskip=14pt
\newtheorem{lemma}{Lemma}[section]
\newtheorem{theorem}[lemma]{Theorem}
\newtheorem{prop}[lemma]{Proposition}
\newtheorem{cor}[lemma]{Corollary}
\newtheorem{definition}[lemma]{Definition}
\newtheorem{definitions}[lemma]{Definitions}
\newtheorem{remark}[lemma]{Remark}
\newtheorem{conj}[lemma]{Conjecture}
\renewcommand{\arraystretch}{2}
\newcommand{\TF}{\widetilde{F}}
\newcommand{\C}{{\mathbb{C}}}
\newcommand{\N}{{\mathbb{N}}}
\newcommand{\Q}{{\mathbb{Q}}}
\newcommand{\R}{{\mathbb{R}}}
\newcommand{\Z}{{\mathbb{Z}}}
\newcommand{\B}{{\mathbf{B}}}
\newcommand{\CA}{{\mathcal{A}}}
\newcommand{\CL}{{\mathcal{L}}}
\newcommand{\0}{-\,}
\newcommand{\cdi}{\mbox{CD$(\mathbf{i})$}}
\newcommand{\cdii}{\mbox{CD$(\mathbf{i}')$}}
\begin{center}
\baselineskip=11pt
{\LARGE Regions of linearity, {L}usztig cones and canonical basis
elements for the quantized enveloping algebra of type ${A}_4$}

{\Large Roger Carter}

{\em Mathematics Institute, University of Warwick, Coventry CV4 7AL, England}
\\
{E-mail: rwc@maths.warwick.ac.uk}

{\Large Robert Marsh}

{\em Department of Mathematics and Computer Science, University of Leicester,
University Road, Leicester LE1 7RH, England} \\
{E-mail: R.Marsh@mcs.le.ac.uk}

\baselineskip=14pt
{\normalsize To Professor Helmut Wielandt on his 90th birthday.}

{\bf Abstract} \\
\parbox[t]{5.15in}{
Let $U_q$ be the quantum group associated to a Lie algebra
$\bf g$ of rank $n$. The negative part $U^-$ of $U$ has
a canonical basis $\mathbf{B}$ with favourable properties
(see Kashiwara~\cite{kash2} and Lusztig~\cite[\S14.4.6]{lusztig6}). The
approaches of Lusztig and Kashiwara lead to a set of
alternative parametrizations of the canonical basis, one for each reduced
expression for the longest word in the Weyl group of $\bf g$.
We show that if $\bf g$ is of type $A_4$ there are close relationships
between the Lusztig cones, canonical basis elements and the regions of
linearity of reparametrization functions arising from the above
parametrizations. A graph can be defined on the set of simplicial regions
of linearity with respect to adjacency, and we further show that this graph is
isomorphic to the graph with vertices given by the reduced expressions
of the longest word of the Weyl group modulo commutation and edges given by
long braid relations.

} \\ \ \\
\end{center}

\section{Introduction}

Let $U=U_q({\mathbf{g}})$ be the quantum group associated to a semisimple
Lie algebra $\mathbf{g}$ of rank $n$. The negative part $U^-$ of $U$ has
a canonical basis $\mathbf{B}$ with favourable properties
(see Kashiwara~\cite{kash2} and Lusztig~\cite[\S14.4.6]{lusztig6}).
For example, via action on highest weight vectors it gives rise to bases
for all the finite-dimensional irreducible highest weight $U$-modules.

Let $W$ be the Weyl group of ${\mathbf{g}}$, with Coxeter generators
$s_1,s_2,\ldots ,s_n$, and let $w_0$ be the element of maximal length in $W$.
Let ${\mathbf{i}}$ be a reduced expression for $w_0$, i.e. $w_0=s_{i_1}s_{i_2}
\cdots s_{i_k}$ is reduced. Lusztig obtains a parametrization of the
canonical basis $\mathbf{B}$ for each such reduced expression $\mathbf{i}$,
via a correspondence between a basis of PBW-type associated to $\mathbf{i}$
and the canonical basis. This gives a bijection
$$\phi_{\mathbf{i}}:\B\rightarrow \N^k,$$
where $\N=\{0,1,2,\ldots \}$.

Kashiwara, in his approach to the canonical basis (which he calls the global
crystal basis), defines certain root operators $\TF_i$ on the canonical
basis (see~\cite[\S3.5]{kash2}) which lead to a parametrization of the
canonical basis for each reduced expression $\mathbf{i}$ by a certain subset
$Y_{\mathbf{i}}$ of $\N^k$. This gives a bijection $\psi_{\mathbf{i}}:\B\rightarrow
Y_{\mathbf{i}}$. The subset $Y_{\mathbf{i}}$ is called the string cone.

In Lusztig's theory, an important role is played by two specific reduced
expressions ${\mathbf{j}}$ and ${\mathbf{j}}'$, which in type $A_4$ with
Dynkin diagram as in Figure $1$ are $\mathbf{j}=(1,3,2,4,1,3,2,4,1,3)$ and
$\mathbf{j}'=(2,4,1,3,2,4,1,3,2,4)$.

\beginpicture

\setcoordinatesystem units <1cm,1cm>             
\setplotarea x from -3 to 6, y from 0.8 to 3       


\put{Figure $1$: Dynkin diagram of type $A_4$} [c] at 4.4 0.75
\scriptsize{

\multiput {$\circ$} at 3   2 *3 1 0 /      %

\linethickness=1pt           

\putrule from 3.05 2 to 3.95 2  %
\putrule from 4.05 2 to 4.95 2  
\putrule from 5.05 2 to 5.95 2  %

\put{$1$}   [l] at 2.9 1.75
\put{$2$}   [l] at 3.9 1.75
\put{$3$}   [l] at 4.9 1.75
\put{$4$}   [l] at 5.9 1.75

}
\endpicture

The function $R_{\mathbf{j}}^{{\mathbf{j}}'}=\phi_{{\mathbf{j}}'}\phi_{\mathbf{j}}^{-1}:\N^k\rightarrow \N^k$
was shown by Lusztig to be piecewise-linear, and the regions of linearity of
this function were shown to be relevant to understanding the behaviour of
the canonical basis.

The function $S_{\mathbf{i}}^{\mathbf{j}}=\phi_{\mathbf{j}}\psi_{\mathbf{i}}^{-1}:Y_{\mathbf{i}}\rightarrow \N^k$
is useful in relating Kashiwara's and Lusztig's parametrizations of $\B$.
These re-parametrization functions have recently
been studied using an approach involving totally positive 
varieties, in the preprint~\cite{bz3} of Berenstein and Zelevinsky 
(although not in terms of regions of linearity).

There is no simple way to express the elements of $\B$ in terms of the
natural generators $F_1,F_2,\ldots ,F_n$ of $U^-$. This has been done only
in types $A_1,A_2,A_3$ and $B_2$
(see~\cite[\S3.4]{lusztig2},~\cite[\S13]{lusztig7},~\cite{xi2}
and~\cite{xi3}) and appears to become arbitrarily complicated in general.
However it seems that there is an interesting subset of $\B$ whose elements
are expressible as monomials in $F_1,F_2,\ldots ,F_n$.

A monomial $F_{i_1}^{(a_1)}F_{i_2}^{(a_2)}\cdots F_{i_k}^{(a_k)}$, where
$F_i^{(a)}=F_i^a/[a]!$, is said to be tight if it belongs to ${\B}$.
Lusztig~\cite{lusztig7} described a method which in low rank cases leads
to the construction of tight monomials. He defined, for each reduced
expression $\mathbf{i}$ of $w_0$, a certain cone $C_{\mathbf{i}}$ in $\N^k$ which we shall
call the Lusztig cone, and showed that, for types $A_1,A_2,A_3$, the above
monomial is tight for all $\mathbf{a}\in C_{\mathbf{i}}$. In~\cite{me7} the second author
showed that this is also true in type $A_4$. Let
$$M_{\mathbf{i}}=\{F_{i_1}^{(a_1)}F_{i_2}^{(a_2)}\cdots F_{i_k}^{(a_k)}\,:\,\mathbf{a}\in
C_{\mathbf{i}}\}$$
be the set of monomials obtained from elements of $C_{\mathbf{i}}$.

There appears to be an intriguing relationship between the Lusztig cones,
the regions of linearity of Lusztig's function $R_{\mathbf{j}}^{{\mathbf{j}}'}$, and the
tight monomials in the canonical basis. This relationship is not fully
understood in general, but the aim of this paper is to describe this
relationship in the type $A_4$ case. Types $A_1,A_2,A_3$ were considered by
Lusztig, but the $A_4$ case is considerably more complicated, while at the
same time being amenable to explicit calculation.

In type $A_4$ we can show the following. For each $\mathbf{a}\in C_{\mathbf{i}}$ we have
$$\TF_{i_1}^{a_1}\TF_{i_2}^{a_2}\cdots \TF_{i_k}^{a_k}\cdot 1\equiv
F_{i_1}^{(a_1)}F_{i_2}^{(a_2)}\cdots F_{i_k}^{(a_k)}\mod v\mathcal{L}',$$
where $\mathcal{L}'$ is the $\CA$-lattice spanned by $\B$ and $\CA$ is the subring
of $\Q(v)$ of functions regular at $v=0$. Since it has been shown
(independently by the second author~\cite{me9} and Premat~\cite{premat1})
that $C_{\mathbf{i}}\subseteq Y_{\mathbf{i}}$, it follows that $M_{\mathbf{i}}\subseteq \B$ and
that, under the Kashiwara parametrization of $\B$, we have
$$\psi_{\mathbf{i}}(M_{\mathbf{i}})=C_{\mathbf{i}}.$$
We can also show that, under Lusztig's parametrization of $\B$, we have
$$\phi_{\mathbf{j}}(M_{\mathbf{i}})=X_{\mathbf{i}},$$
where $X_{\mathbf{i}}$ is one of the regions of linearity of the
piecewise-linear
function $R_{\mathbf{j}}^{\mathbf{j}'}$ and that the transition function
$$S_{\mathbf{i}}^{\mathbf{j}}:C_{\mathbf{i}}\rightarrow X_{\mathbf{i}}$$
is linear for all $\mathbf{i}$.

The set of all regions of linearity of $R_{\mathbf{j}}^{{\mathbf{j}}'}$ in
type $A_4$ was determined by the first author using the ideas outlined
in~\cite{carter3}. These
regions are described in the present paper. Each such region is defined by
a certain set of inequalities. It turns out that the regions $X_{\mathbf{i}}$ are
all defined by the minimal possible number of inequalities, and that they
give the set of all regions defined by this minimal number of inequalities.
We call these the simplicial regions of $R_{\mathbf{j}}^{{\mathbf{j}}'}$. We thus have a
parametrization of the simplicial regions in terms of the reduced expressions $\mathbf{i}$
of $w_0$ modulo commutation. A graph can be defined on the set of simplicial
regions with respect to adjacency, and we further show that this graph is
isomorphic to the graph with vertices given by the reduced expressions
$\mathbf{i}$ of $w_0$ modulo commutation and edges given by long braid relations.

Thus in type $A_4$ there are close relationships between the Lusztig cones
$C_{\mathbf{i}}$, the tight monomials $M_{\mathbf{i}}$, and the regions of linearity
$X_{\mathbf{i}}$. Examples of N. H. Xi~\cite{xi4} and
M. Reineke~\cite{reineke2} show, however, that
these relationships cannot be expected to hold in the same way in type $A_n$
for arbitrary $n$.

\section{Parametrizations of the canonical basis}

Let $\mathbf{g}$ be the simple Lie algebra over $\C$ of type $A_n$ and $U$ be the
quantized enveloping algebra of $\mathbf{g}$. Then $U$ is a $\Q(v)$-algebra generated by the
elements $E_i$, $F_i$, $K_{\mu}$, $i\in \{1,2,\ldots ,n\}$, $\mu\in Q$,
the root
lattice of $\mathbf{g}$. Let $U^+$ be the subalgebra generated by the $E_i$ and $U^-$
the subalgebra generated by the $F_i$.

Let $W$ be the Weyl group of $\mathbf{g}$. It has a unique element $w_0$
of maximal length. For each reduced expression $\mathbf{i}$ for $w_0$ there are two
parametrizations of the canonical basis $\B$ for $U^-$. The first arises from
Lusztig's approach to the canonical basis~\cite[\S14.4.6]{lusztig6}, and
the second arises from Kashiwara's approach~\cite{kash2}.

\noindent {\bf Lusztig's Approach} \\
There is an $\Q$-algebra automorphism of $U$ which takes each $E_i$ to
$F_i$, $F_i$ to $E_i$, $K_{\mu}$ to $K_{-\mu}$ and $v$ to $v^{-1}$.
We use this automorphism to transfer Lusztig's
definition of the canonical basis in~\cite[\S3]{lusztig2} to $U^-$.

Let $T_i$, $i=1,2,\ldots ,n$, be the automorphism of $U$ as
in~\cite[\S1.3]{lusztig1} given by:

$$T_i(E_j)=\left\{ \begin{array}{cc}
-F_jK_j, & {\rm if\ }i=j, \\
E_j,     & {\rm if\ }|i-j|>1 \\
-E_iE_j+v^{-1}E_jE_i & {\rm if\ }|i-j|=1
\end{array} \right.$$

$$T_i(F_j)=\left\{ \begin{array}{cc}
-K_j^{-1}E_j, & {\rm if\ }i=j, \\
F_j,     & {\rm if\ }|i-j|>1 \\
-F_jF_i+vF_iF_j & {\rm if\ }|i-j|=1
\end{array} \right.$$

$$T_i(K_{\mu})=K_{\mu-\langle \mu,\alpha_i\rangle h_i},{\rm \ for\ }\mu\in
Q,$$ where $\alpha_i$ are the simple roots and $h_i$ the simple coroots of
$\mathbf{g}$.

For each $i$, let $r_i$ be the automorphism of $U$ which fixes $E_j$ and
$F_j$ for $j=i$ or $|i-j|>1$ and fixes $K_{\mu}$ for all $\mu$,
and which takes $E_j$ to $-E_j$ and $F_j$ to $-F_j$ if $|i-j|=1$.
Let $T''_{i,-1}=T_ir_i$ be the automorphism of $U$ as
in~\cite[\S37.1.3]{lusztig6}.
Let ${\mathbf{c}}\in {\N}^k$, where $k=\ell (w_0)$, and
$\mathbf{i}$ be a reduced expression for $w_0$. Let
$$F_{\mathbf{i}}^{\mathbf{c}}:=F_{i_1}^{(c_1)}T''_{i_1,-1}(F_{i_2}^{(c_2)})\cdots
T''_{i_1,-1}T''_{i_2,-1}\cdots T''_{i_{k-1},-1}(F_{i_k}^{(c_k)}).$$
Define $B_{\mathbf{i}}=\{F_{\mathbf{i}}^{\mathbf{c}}\,:\, {\mathbf{c}}\in {\N}^k\}$.
Then $B_{\mathbf{i}}$ is
the basis of PBW-type corresponding to the reduced expression $\mathbf{i}$.
Let\,\, $\mathbf{\bar{\ }}\,$ be the ${\Q}$-algebra automorphism from $U$ to $U$ taking $E_i$ to $E_i$, $F_i$ to $F_i$, and $K_{\mu}$ to $K_{-\mu}$,
for each $i\in [1,n]$ and $\mu\in Q$, and $v$ to $v^{-1}$.
Lusztig proves the following result in~\cite[\S\S2.3, 3.2]{lusztig2}.

\begin{theorem} (Lusztig) \\
The $Z[v]$-span $\CL$ of $B_{\mathbf{i}}$ is independent of $\mathbf{i}$. 
Let $\pi:{{\CL}}\rightarrow {{\CL}}/{v{\CL}}$ be the natural
projection. The image $\pi(B_{\mathbf{i}})$ is also independent
of $\mathbf{i}$; we denote it by $B$. The restriction of $\pi$ to
${{\CL}}\cap \overline{{\CL}}$ is an isomorphism of $\mathbb{Z}$-modules
$\pi_1:{{\CL}}\cap \overline{{\CL}}\rightarrow {{\CL}}/{v{\CL}}$.
Also $\B=\pi_1^{-1}(B)$ is a
${\Q}(v)$-basis of $U^-$, which is the canonical basis of $U^-$.
\end{theorem}

Lusztig's theorem provides us with a parametrization of ${\B}$, dependent on
$\mathbf{i}$. If $b\in {\B}$, we write $\phi_{\mathbf{i}}(b)={\mathbf{c}}$,
where ${\mathbf{c}}\in {\N}^k$ satisfies $b\equiv F_{\mathbf{i}}^{\mathbf{c}}
\mod v{{\CL}}$.
Note that $\phi_{\mathbf{i}}$ is a bijection.

Lusztig considers in~\cite{lusztig2} two particular reduced expressions
for $w_0$. Let the nodes in the Dynkin diagram of $A_n$ be labelled as in
Figure $2$.

\begin{figure}[htbp]

\beginpicture

\setcoordinatesystem units <1cm,1cm>             
\setplotarea x from -2.7 to 12, y from 0.8 to 3       

\put{Figure $2$: Dynkin diagram of type $A_n$}[c] at 5.5 0.75   

\scriptsize{

\multiput {$\circ$} at 3   2 *1 1 0 /      %
\multiput {$\circ$} at 7   2 *1 1 0 /      

\linethickness=1pt           

\putrule from 3.05 2 to 3.95 2  %
\putrule from 7.05 2 to 7.95 2  

\setdashes <2mm,1mm>          %
\putrule from 4.05 2 to 6.95 2  

\put{$1$}   [l] at 2.9 1.75
\put{$2$}   [l] at 3.9 1.75
\put{$n-1$} [l] at 6.65 1.75
\put{$n$} [l] at 7.9 1.75

}
\endpicture
\end{figure}

Let $\mathbf{j} = 135\cdots 246\cdots 135\cdots 246\cdots$, and let
$\mathbf{j}' = 246\cdots 135\cdots 246\cdots 135\cdots$,
where both expressions have length $k=\ell (w_0)$ --- they are both
reduced expressions for $w_0$. We have bijections
$\phi_{\mathbf{j}}:\B \rightarrow \N^k$ and
$\phi_{\mathbf{j}'}:\B \rightarrow \N^k$.

Lusztig defines in~\cite[\S2.6]{lusztig2}
a function $R_{\mathbf{j}}^{{\mathbf{j}}'}=\phi_{{\mathbf{j}}'}
\phi_{\mathbf{j}}^{-1}\,:\,{\N}^k \rightarrow {\N}^k$. This function
was shown by Lusztig to be piecewise linear and its regions of linearity were
shown to have significance for the canonical basis, in the sense that
elements $b$ of the canonical basis with $\phi_{\mathbf{j}}(b)$ in the same
region of linearity of $R_{\mathbf{j}}^{{\mathbf{j}}'}$ often have
similar form.

\noindent {\bf Kashiwara's approach}

Let $\widetilde{E}_i$ and $\widetilde{F}_i$ be the Kashiwara operators on
$U^-$ as defined in~\cite[\S3.5]{kash2}. Let ${\CA}\subseteq {\Q}(v)$ be
the subring of elements regular at $v=0$, and let ${\CL}'$ be the
${\CA}$-lattice spanned by arbitrary products
$\widetilde{F}_{j_1}\widetilde{F}_{j_2}\cdots \widetilde{F}_{j_m}\cdot 1$ in $U^-$.
We denote the set of all such elements by $S$. The following results
were proved by Kashiwara in~\cite{kash2}.

\begin{theorem} \label{kashiwara} (Kashiwara) \\
(i) Let $\pi':{{\CL}'}\rightarrow {{\CL}'}/{v{{\CL}'}}$ be the natural
projection, and let $B'=\pi'(S)$. Then $B'$ is a ${\Q}$-basis of
${{\CL}'}/{v{{\CL}'}}$ (the crystal basis). \\
(ii) Furthermore,
$\widetilde{E}_i$ and $\widetilde{F}_i$ each preserve ${\CL}'$ and thus act on
${{\CL}'}/{v{{\CL}'}}$. They satisfy $\widetilde{E}_i(B')\subseteq B'\cup\{0\}$ and
$\widetilde{F}_i(B')\subseteq B'$. Also for $b,b'\in B'$ we have $\widetilde{F}_ib=b'$,
if and only if $\widetilde{E}_ib'=b$. \\
(iii) For each $b\in B'$, there is a unique element $\widetilde{b}\in {{\CL}'}
\cap \overline{{\CL}'}$ such that $\pi'(\widetilde{b})=b$. The set of elements
$\{\widetilde{b}\,:\, b\in B'\}$ forms a basis of $U^-$, the {\em global
crystal basis} of $U^-$.
\end{theorem}

It was shown by Lusztig~\cite[2.3]{lusztig3} that the global crystal basis
of Kashiwara coincides with the canonical basis of $U^-$.

There is a parametrization of ${\B}$ arising from Kashiwara's approach,
again dependent on a reduced expression $\mathbf{i}$ for $w_0$.
Let $\mathbf{i}=(i_1,i_2,\ldots ,i_k)$ and $b\in B$.
Let $a_1$ be maximal such that $\widetilde{E}_{i_1}^{a_1}b\not\equiv 0 \mod
v{{\CL}'}$;
let $a_2$ be maximal such that
$\widetilde{E}_{i_2}^{a_2}\widetilde{E}_{i_1}^{a_1}b\not\equiv 0 \mod v{{\CL}'}$,
and so on, so that
$a_k$ is maximal such that $\widetilde{E}_{i_k}^{a_k}\widetilde{E}_{i_{k-1}}^{a_{k-1}}
\cdots \widetilde{E}_{i_2}^{a_2}\widetilde{E}_{i_1}^{a_1}b\not\equiv 0 \mod
v{{\CL}'}$.
Let $\mathbf{a}=(a_1,a_2,\ldots ,a_k)$. We write
$\psi_{\mathbf{i}}(b)=\mathbf{a}$.
This is the crystal string of $b$ --- see~\cite[\S2]{bz1}
and the end of Section $2$ in~\cite{nz1}; see also~\cite{kash4}.
It is known that $\psi_{\mathbf{i}}(b)$ uniquely determines $b\in {\B}$
(see~\cite[\S2.5]{nz1}).
We have $b\equiv \TF_{i_1}^{a_1}\TF_{i_2}^{a_2}\cdots \TF_{i_k}^{a_k}\cdot 1
\mod v{\CL}'$.
The image of $\psi_{\mathbf{i}}$ is a cone which first appears
in~\cite{bz1}.
We shall call this the {\em string cone} $Y_{\mathbf{i}}=
\psi_{\mathbf{i}}(\B)$.

We next define a function which compares Kashiwara's approach with
Lusztig's approach. Consider the maps
$$Y_{\mathbf{i}}\xrightarrow[\psi_{\mathbf{i}}^{-1}]{} B
\xrightarrow[\phi_{\mathbf{j}}]{} \N^k$$
where $\mathbf{j}=135\cdots 246\cdots 135\cdots 246\cdots $.
We define $S_{\mathbf{i}}^{\mathbf{j}}=\phi_{\mathbf{j}}
\psi_{\mathbf{i}}^{-1}\,:\,Y_{\mathbf{i}}\rightarrow {\N}^k$,
a reparametrization function. We shall show that this function has some
interesting properties in the case in which $\mathbf{g}$ has type $A_4$.

\section{The Lusztig cones and their spanning vectors}

Lusztig~\cite{lusztig7}
introduced certain regions which, in low rank, give rise to
canonical basis elements of a particularly simple form. We consider
reduced expressions $\mathbf{i}=(i_1,i_2,\ldots ,i_k)$ for $w_0$. We shall
identify this $k$-tuple with the reduced expression $s_{i_1}s_{i_2}\cdots
s_{i_k}$. Given two such reduced expressions we say that $\mathbf{i}\sim
\mathbf{i}'$ if there is a sequence of commutations (of the form
$s_is_j=s_js_i$ with $|i-j|>1$) which, when applied to $\mathbf{i}$, give
$\mathbf{i}'$. This is an equivalence relation on the set of reduced
expressions for $w_0$, and the equivalence classes are called
commutation classes.

The {\em Lusztig cone}, $C_{\mathbf{i}}$, corresponding
to a reduced expression $\mathbf{i}$ for $w_0$ is defined to be the set of
points ${{\mathbf{a}}}\in {\N}^k$
satisfying the following inequalities: \\
(*) For every pair $s,s'\in [1,k]$ with $s<s'$, $i_s=i_{s'}=i$ and
$i_p\not=i$ whenever $s<p<s'$, we have
$$(\sum_p a_p)-a_s - a_{s'}\geq 0,$$
where the sum is over all $p$ with $s<p<s'$ such that $i_p$ is joined to $i$
in the Dynkin diagram.

It was shown by Lusztig~\cite{lusztig7} that if $\mathbf{a}\in C_{\mathbf{i}}$
then the monomial $F_{i_1}^{(a_1)}F_{i_2}^{(a_2)}\cdots
F_{i_k}^{(a_k)}$ lies in the canonical basis ${\B}$, provided $n=1,2,3$.
The second author~\cite{me7} showed that this remains true if $n=4$.
The Lusztig cones have been studied in the second author's
preprints~\cite{me8} and~\cite{me9} in type $A$ for every reduced expression
$\bf i$ for the longest word, and have also been studied by Bedard
in~\cite{bedard2} for arbitrary finite
(simply-laced) type for reduced expressions compatible with a quiver whose
underlying graph is the Dynkin diagram. Bedard describes these vectors
using the Auslander-Reiten quiver of the quiver and homological algebra,
showing they are closely connected to the representation theory of the quiver.

The reduced expression $\mathbf{i}$ defines an ordering on the set $\Phi^+$ of
positive roots of the root system associated to $W$.
We write $\alpha^j=s_{i_1}s_{i_2}\cdots s_{i_{j-1}}(\alpha_{i_j})$ for
$j=1,2,\ldots ,k$. Then $\Phi^+=\{\alpha^1,\alpha^2,\ldots \alpha^k\}$.
For ${{\mathbf{a}}}=(a_1,a_2,\ldots
,a_k)\in \mathbb{Z}^k$, write $a_{\alpha^j}=a_j$.
If $\alpha=\alpha_{ij}:=\alpha_i+\alpha_{i+1}+\cdots +\alpha_{j-1}$
with $i<j$, we also write $a_{ij}$ for $a_{\alpha_{ij}}$.

We can regard $C_{\mathbf{i}}$ as a subset of
$\mathbb{Z}^k$ defined by the inequalities in (*) above, together
with the $n$ inequalities $a_{\alpha_j}\geq 0$ for $j=1,2,\ldots ,n$
(see the paragraph before Lemma 4.2 in~\cite{me8}).
The number of inequalities in (*) is $k-n$, thus we have $k$ inequalities
altogether defining $C_{\mathbf{i}}$. There is therefore a
matrix $P_{\mathbf{i}}\in M_k(\mathbb{Z})$ such that
$$C_{\mathbf{i}}=\{{{\mathbf{a}}}\in \mathbb{Z}^k\,:\,
P_{\mathbf{i}}{{\mathbf{a}}}\geq 0\},$$
where, for $\mathbf{z}\in \mathbb{Z}^k$, $\mathbf{z}\geq 0$
means that each entry in $\mathbf{z}$ is
nonnegative.

There is a description of $C_{\mathbf{i}}$ in~\cite{me8} which will be useful
to us. We shall need the chamber diagram (chamber ansatz) for $\mathbf{i}$
defined in~\cite[\S\S1.4, 2.3]{bfz1}.
We take $n+1$ strings, numbered from top to bottom, and write 
$\mathbf{i}$ from left to right along the bottom of the diagram.
Above a letter $i_j$ in $\mathbf{i}$,
the $i_j$th and $(i_j+1)$st strings from the top above $i_j$ cross.
Thus, for example, in the case $n=3$ with $\mathbf{i}=(1,3,2,1,3,2)$, the
chamber diagram is shown in Figure $3$.
\begin{figure}[htbp]
\beginpicture

\setcoordinatesystem units <0.5cm,0.3cm>             
\setplotarea x from -9 to 14, y from -2 to 7       
\linethickness=0.5pt           

\put{$1$}[c] at -0.3 6 %
\put{$2$}[c] at -0.3 4 %
\put{$3$}[c] at -0.3 2 %
\put{$4$}[c] at -0.3 0 %

\put{$1$}[c] at 1.5 -1.5 %
\put{$3$}[c] at 3.5 -1.5 %
\put{$2$}[c] at 5.5 -1.5 %
\put{$1$}[c] at 7.5 -1.5 %
\put{$3$}[c] at 9.5 -1.5 %
\put{$2$}[c] at 11.5 -1.5 %

\setlinear \plot 0 0  3 0 / %
\setlinear \plot 3 0  4 2 / %
\setlinear \plot 4 2  5 2 / %
\setlinear \plot 5 2  6 4 / %
\setlinear \plot 6 4  7 4 / %
\setlinear \plot 7 4  8 6 / %
\setlinear \plot 8 6  13 6 / %

\setlinear \plot 0 2  3 2 / %
\setlinear \plot 3 2  4 0 / %
\setlinear \plot 4 0  9 0 / %
\setlinear \plot 9 0  10 2 / %
\setlinear \plot 10 2  11 2 / %
\setlinear \plot 11 2  12 4 / %
\setlinear \plot 12 4  13 4 / %

\setlinear \plot 0 4  1 4 / %
\setlinear \plot 1 4  2 6 / %
\setlinear \plot 2 6  7 6 / %
\setlinear \plot 7 6  8 4 / %
\setlinear \plot 8 4  11 4 / %
\setlinear \plot 11 4  12 2 / %
\setlinear \plot 12 2  13 2 / %

\setlinear \plot 0 6  1 6 / %
\setlinear \plot 1 6  2 4 / %
\setlinear \plot 2 4  5 4 / %
\setlinear \plot 5 4  6 2 / %
\setlinear \plot 6 2  9 2 / %
\setlinear \plot 9 2  10 0 / %
\setlinear \plot 10 0  13 0 / %

\put{$134$}[c] at 4.5 5 %
\put{$13$}[c] at 8.5 3 %
\put{$3$}[c] at 6.5 1 %

\put{Figure 3: Chamber Diagram}[c] at 6.5 -4.5 %
\endpicture

\end{figure}

We denote the chamber diagram of $\mathbf{i}$ by CD($\mathbf{i}$).
We shall be concerned with the bounded chambers of the diagram.
A {\em chamber} will be defined as a pair $(c,{\mathbf{i}})$, where $c$ is a
bounded component of the complement of CD($\mathbf{i}$).
Each chamber $(c,{\mathbf{i}})$ can be labelled with the numbers of the
strings passing below it, denoted $l(c,{\mathbf{i}})$.
Following~\cite{bfz1}, we call such a label a chamber set.
For example, the chamber sets corresponding to the $3$ bounded chambers in
Figure $3$ are $134,3,13$.

There is another way to think of the set of bounded chambers in a chamber
diagram. We recall that a quiver of type $A_n$ is a directed graph such that
the underlying undirected graph is the Dynkin diagram of type $A_n$.
Following~\cite[\S5.1]{me8}, we define a {\em partial quiver} $P$
of type $A_n$ to be a quiver of type $A_n$ which has some (or none) of its
arrows replaced by undirected edges in such a way that the subgraph obtained
by deleting undirected edges and vertices incident only with undirected edges
is non-empty and connected.
We shall now number the edges of a partial quiver of type $A_n$ from $2$
to $n$, starting from the right, as in Figure $4$.

\begin{figure}[htbp]

\beginpicture

\setcoordinatesystem units <1cm,1cm>             
\setplotarea x from -2.7 to 12, y from 0.8 to 3       

\put{Figure $4$: Edge Numbering of the Dynkin Diagram.}[c] at 5.5 0.75
\scriptsize{

\multiput {$\circ$} at 2   2 *2 1 0 /      %
\multiput {$\circ$} at 7   2 *2 1 0 /      

\linethickness=1pt           

\putrule from 2.05 2 to 2.95 2  %
\putrule from 3.05 2 to 3.95 2  %
\putrule from 8.05 2 to 8.95 2  %
\putrule from 7.05 2 to 7.95 2  

\setdashes <2mm,1mm>          %
\putrule from 4.05 2 to 6.95 2  

\put{$2$}   [c] at 8.5 2.2
\put{$3$}   [c] at 7.5 2.2
\put{$n-1$} [c] at 3.5 2.2
\put{$n$}   [c] at 2.5 2.2

}

\endpicture

\end{figure}

If $P$ is a partial quiver of type $A_n$, we denote by $l(P)$ the subset of
$[1,n+1]$ defined as follows:
Let $l_1(P)=\{j\in [2,n]\,:\,\mbox{edge\ $j$\ of\ $P$\ is\ an\ 
$L$}\}$ (this means that edge $j$ has an arrow pointing to the left).
If the rightmost directed edge of $P$ is an $R$, and this is
in position $i$, then let $l_2(P)=[1,i-1]$. Otherwise $l_2(P)$ is
the empty set. If the leftmost directed edge of $P$ is an $R$, and this is
edge $j$, let $l_3(P)=[j+1,n+1]$. Otherwise $l_3(P)$ is the empty set.
We then define $l(P)=l_1(P)\cup l_2(P)\cup l_3(P)$.

It is shown in~\cite[\S5.4]{me8} that
the map $l$ is a bijection from the set of partial quivers of type
$A_n$ to the set of all chamber sets.  For example, the chamber sets
associated to the $8$ partial quivers of type $A_3$ are shown in Figure $5$.

\begin{figure}[htbp]

\beginpicture

\setcoordinatesystem units <1cm,1cm>             
\setplotarea x from -4 to 8, y from -1 to 8.5       

\put{Partial Quiver}[c] at 1 8
\put{Chamber Set}[c] at 6 8
\scriptsize{

\multiput {$\circ$} at 0 0  *7 0 1 /      %
\multiput {$\circ$} at 1 0  *7 0 1 /      %
\multiput {$\circ$} at 2 0  *7 0 1 /      %

\linethickness=1pt           %

\putrule from 0.05 0 to 0.95 0  %
\putrule from 0.05 1 to 0.95 1  %
\putrule from 0.05 2 to 0.95 2  %
\putrule from 0.05 3 to 0.95 3  %
\putrule from 0.05 4 to 0.95 4  %
\putrule from 0.05 5 to 0.95 5  %
\putrule from 0.05 6 to 0.95 6  %
\putrule from 0.05 7 to 0.95 7  

\putrule from 1.05 0 to 1.95 0  %
\putrule from 1.05 1 to 1.95 1  %
\putrule from 1.05 2 to 1.95 2  %
\putrule from 1.05 3 to 1.95 3  %
\putrule from 1.05 4 to 1.95 4  %
\putrule from 1.05 5 to 1.95 5  %
\putrule from 1.05 6 to 1.95 6  %
\putrule from 1.05 7 to 1.95 7  

\setlinear \plot 1.4 0.1  1.5 0 / %
\setlinear \plot 1.4 -0.1  1.5 0 / %
\put{$R$} at 1.5 0.3

\setlinear \plot 0.4 1.1  0.5 1 / %
\setlinear \plot 0.4 0.9  0.5 1 / %
\put{$R$} at 0.5 1.3

\setlinear \plot 0.4 4.1  0.5 4 / %
\setlinear \plot 0.4 3.9  0.5 4 / %
\put{$R$} at 0.5 4.3

\setlinear \plot 1.4 4.1  1.5 4 / %
\setlinear \plot 1.4 3.9  1.5 4 / %
\put{$R$} at 1.5 4.3

\setlinear \plot 0.4 5.1  0.5 5 / %
\setlinear \plot 0.4 4.9  0.5 5 / %
\put{$R$} at 0.5 5.3

\setlinear \plot 1.4 6.1  1.5 6 / %
\setlinear \plot 1.4 5.9  1.5 6 / 
\put{$R$} at 1.5 6.3

\setlinear \plot 1.6 2.1  1.5 2 / %
\setlinear \plot 1.6 1.9  1.5 2 / %
\put{$L$} at 1.5 2.3

\setlinear \plot 0.6 3.1  0.5 3 / %
\setlinear \plot 0.6 2.9  0.5 3 / %
\put{$L$} at 0.5 3.3

\setlinear \plot 1.6 5.1  1.5 5 / %
\setlinear \plot 1.6 4.9  1.5 5 / %
\put{$L$} at 1.5 5.3

\setlinear \plot 0.6 6.1  0.5 6 / %
\setlinear \plot 0.6 5.9  0.5 6 / %
\put{$L$} at 0.5 6.3

\setlinear \plot 0.6 7.1  0.5 7 / %
\setlinear \plot 0.6 6.9  0.5 7 / %
\put{$L$} at 0.5 7.3

\setlinear \plot 1.6 7.1  1.5 7 / %
\setlinear \plot 1.6 6.9  1.5 7 / 
\put{$L$} at 1.5 7.3

\put{$\{1,3,4\}$}[c] at 6 0
\put{$\{1,2,4\}$}[c] at 6 1
\put{$\{2\}$}[c] at 6 2
\put{$\{3\}$}[c] at 6 3
\put{$\{1,4\}$}[c] at 6 4
\put{$\{2,4\}$}[c] at 6 5
\put{$\{1,3\}$}[c] at 6 6
\put{$\{2,3\}$}[c] at 6 7

}
\normalsize \put{Figure $5$: Partial Quivers of Type $A_3$}[c] at 3 -1

\endpicture

\end{figure}

We now consider spanning vectors for the Lusztig cone $C_{\mathbf{i}}$.
The $k\times k$ matrix $P_{\mathbf{i}}$ has an inverse
$Q_{\mathbf{i}}$ with entries in $\N$ (see~\cite[\S4.2]{me8}).
We have $k=\frac{1}{2}n(n+1)$. Note that $k-n=\frac{1}{2}n(n-1)$
of the rows of $P_{\mathbf{i}}$ correspond to inequalities arising from
consecutive occurrences of a letter in $\mathbf{i}$.
Each such consecutive pair corresponds naturally to a bounded chamber
$(c,{\mathbf{i}})$. Thus for each chamber $(c,{\mathbf{i}})$,
there is a corresponding row of $P_{\mathbf{i}}$ and
therefore a corresponding column of $Q_{\mathbf{i}}$. Now the columns
of $Q_{\mathbf{i}}$ give spanning vectors of $C_{\mathbf{i}}$ in the sense
that the elements of $C_{\mathbf{i}}$ are nonnegative linear
combinations of such spanning vectors. We denote this spanning vector by
$\mathbf{a}(c,{\mathbf{i}})$.
The remaining $n$ rows of $P_{\mathbf{i}}$ correspond to inequalities
$a_{\alpha_j}\geq 0$ for $\alpha_j$ a simple root.
We denote the corresponding spanning vectors by $\mathbf{a}(j,{\mathbf{i}})$.
Thus the $k$ spanning vectors of $C_{\mathbf{i}}$ are:
\begin{eqnarray*}
\mathbf{a}(c,\mathbf{i}) & \ \ & \mbox{($k-n$\ vectors)}, \\
\mathbf{a}(j,\mathbf{i}) & \ \ & \mbox{($n$\ vectors)}.
\end{eqnarray*}
It is shown in~\cite{me8} how to obtain the spanning vectors 
$\mathbf{a}(c,\mathbf{i})$, $\mathbf{a}(j,\mathbf{i})$. We recall that the
coordinates of $\mathbf{a}=(a_1,a_2,\ldots ,a_k)$ correspond naturally
to the positive roots $\alpha^1,\alpha^2,\ldots ,\alpha^k$, given the
reduced expression $\mathbf{i}$ for $w_0$. It is possible to find the
spanning vectors by attaching to each positive root $\alpha^m$ a multiplicity
$a_m$, which will be the appropriate coordinate of the spanning vector
$\mathbf{a}$. It is shown in~\cite{me8} that the multiplicity function
$\alpha^m
\mapsto a_m$ for $\mathbf{a}(c,\mathbf{i})$ depends only upon the partial
quiver $P$ with $l(P)=l(c,\mathbf{i})$, and that the multiplicity function
for $\mathbf{a}(j,\mathbf{i})$ depends only upon $j$. These multiplicity
functions are constructed by the following algorithm.

For $j\in [1,n]$, the multiplicity of a positive root $\alpha_{pq}=
\alpha_p+\alpha_{p+1}+\cdots +\alpha_{q-1}$ is $1$ if $1\leq p\leq j<j+1\leq q
\leq n+1$ and is $0$ otherwise. This determines the spanning vector
$\mathbf{a}(j,\mathbf{i})$, whose components are all either $0$ or $1$.

For a partial quiver $P$ the multiplicity of a positive root $\alpha_{pq}$ is
given as follows. We consider the components of $P$, i.e. the maximal
connected subquivers all of whose arrows point in the same direction.
For each component $Y$ of $P$ let $a(Y)$ be the number of the leftmost edge
to the right of $Y$ and let $b(Y)$ be the number of the rightmost edge to the
left of $Y$. We have $a(Y)<b(Y)$. The component $Y$ determines a set
$\Phi^+(Y)$ of positive roots $\alpha_{pq}$ such that $1\leq p\leq a(Y)<b(Y)\leq
q\leq n+1$. Let $m_{pq}$ be the number of components $Y$ of $P$ such that
$\alpha_{pq}\in \Phi^+(Y)$. Then the multiplicity $a_{pq}$ of $\alpha_{pq}$
is given by $a_{pq}=\lceil \frac{1}{2}m_{pq} \rceil$,
which is the smallest integer $m$ with $\frac{1}{2}m_{pq}\leq m$.
The multiplicity $a_{pq}$
is the coordinate of $\mathbf{a}(c,\mathbf{i})$ corresponding
to the positive root $\alpha_{pq}$.

\section{Transforms of the Lusztig cones and their spanning vectors}

We now bring into play the parametrization of the canonical basis $\B$
arising from Kashiwara's approach. We recall
there is a bijection
$$\psi_{\mathbf{i}}:\B\rightarrow Y_{\mathbf{i}}$$
between the canonical basis and the string cone $Y_{\mathbf{i}}$. Now it
has been shown by the second author~\cite{me9} and independently by
Premat~\cite{premat1} that $C_{\mathbf{i}}\subseteq Y_{\mathbf{i}}$, i.e.
that the Lusztig cone lies in the string cone. Thus there is a corresponding
subset $\psi_{\mathbf{i}}^{-1}(C_{\mathbf{i}})$ of the canonical basis
$\B$.

We also have a bijection
$$\phi_{\mathbf{i}}:\B\rightarrow \N^k$$
and a transition function
$$S_{\mathbf{i}}^{\mathbf{j}}:Y_{\mathbf{i}}\rightarrow \N^k$$
given by $S_{\mathbf{i}}^{\mathbf{j}}=\phi_{\mathbf{j}}
\psi_{\mathbf{i}}^{-1}$. We consider the question:
what is the subset $S_{\mathbf{i}}^{\mathbf{j}}(C_{\mathbf{i}})$ of
$\N^k$ corresponding to the Lusztig cone $C_{\mathbf{i}}$? We state the
following conjecture:

\begin{conj}
For any reduced expression $\mathbf{i}$ of $w_0$,
$S_{\mathbf{i}}^{\mathbf{j}}(C_{\mathbf{i}})$ is a region of linearity
of Lusztig's piecewise-linear
function $R=R_{\mathbf{j}}^{\mathbf{j}'}$.
\end{conj}

We now describe an algorithm proved by the authors, as yet unpublished, which
gives the transforms $S_{\mathbf{i}}^{\mathbf{j}}(\mathbf{a}(c,\mathbf{i}))$,
$S_{\mathbf{i}}^{\mathbf{j}}(\mathbf{a}(j,\mathbf{i}))$ of the spanning vectors of
$C_{\mathbf{i}}$.

We first define a $(p,q,r,s)$-rectangle. Suppose we begin with the
array of numbers as in Figure $6$.

\begin{figure}[htbp]

\beginpicture

\setcoordinatesystem units <0.5cm,0.5cm>             
\setplotarea x from -8.3 to 10, y from 0 to 7       

\put{Figure $6$: Array of numbers for defining a $(p,q,r,s)$ rectangle.}[c] at 6.5 0

\linethickness=0.5pt           

\multiput {$\cdot$} at 1   6 *2 1 0 /      %
\multiput {$\cdot$} at 1   5 *2 1 0 /      %
\multiput {$\cdot$} at 1   4 *2 1 0 /      %
\multiput {$\cdot$} at 10   6 *2 1 0 /      %
\multiput {$\cdot$} at 10   5 *2 1 0 /      %
\multiput {$\cdot$} at 10   4 *2 1 0 /      %

\multiput {$\cdot$} at 4   3 *5 1 0 /      %
\multiput {$\cdot$} at 4   2 *5 1 0 /      %

\multiput {$0$} at 4   6 *5 1 0 /      %
\multiput {$1$} at 4   5 *5 1 0 /      %
\multiput {$2$} at 4   4 *5 1 0 /      %

\endpicture

\end{figure}
Let $(p,q,r,s)\in {\N}^4$. A $(p,q,r,s)$-rectangle is a rectangle in
this array, with top vertex on line $p$, middle vertices on lines $q,r$,
and bottom vertex on line $s$ and with all vertices lying half-way
between two numbers.
All lines of the rectangle are at an angle of $\pm \pi/4$ to the horizontal.
This forces $p<q<s$, $p<r<s$ and $p+s=q+r$.
We take only alternate columns of numbers in the rectangle, starting with
the first column if $q$ is odd, and with the second column if $q$ is even.
A $(0,3,5,8)$-rectangle is shown in Figure $7$.
\begin{figure}[htbp]

\beginpicture

\setcoordinatesystem units <0.5cm,0.5cm>             
\setplotarea x from -10 to 9, y from 0 to 10       

\linethickness=1pt        
\setlinear \plot 0.5 6 	3.5 9 	/ %
\setlinear \plot 3.5 9 	8.5 4	/ %
\setlinear \plot 8.5 4 	5.5 1 	/ %
\setlinear \plot 5.5 1 	0.5 6 	/ %

\multiput {$0$}[c] at 1  9  *7 1 0 /      %
\multiput {$1$}[c] at 1  8  *2 1 0 /      %
\multiput {$1$}[c] at 5  8  *3 1 0 /      %
\multiput {$2$}[c] at 1  7  *3 2 0 /      %
\multiput {$2$}[c] at 8  7  *0 1 0 /      %
\multiput {$3$}[c] at 1  6  *3 2 0 /      %
\multiput {$3$}[c] at 8  6  *0 1 0 /      %
\multiput {$4$}[c] at 1  5  *3 2 0 /      %
\multiput {$5$}[c] at 1  4  *2 1 0 /      %
\multiput {$5$}[c] at 5  4  *1 2 0 /      %
\multiput {$6$}[c] at 1  3  *2 1 0 /      %
\multiput {$6$}[c] at 5  3  *1 2 0 /      %
\multiput {$7$}[c] at 1  2  *4 1 0 /      %
\multiput {$7$}[c] at 7  2  *1 1 0 /      %
\multiput {$8$}[c] at 1  1  *7 1 0 /      %

\put{Figure $7$: A $(0,3,5,8)$ rectangle.}[c] at 4.5 -1

\endpicture

\end{figure}
Let $P$ be a partial quiver. Let $Y$ be
a component of $P$ and let $a=a(Y)$ and let $b=b(Y)$ be the integers
defined above.
Let $\rho(Y)$ be a $(0,a,n+2-b,n+a-b+2)$-rectangle if $Y$ has type
$L$, and a $(b-a-1,b-1,n+1-a,n+1)$-rectangle if $Y$ is of type $R$.

Now let $j\in [1,n]$. Let $a=j$ and $b=j+1$. Let $\rho(j)$ be a
$(0,a,n+2-b,n+a-b+2)=(0,j,n+1-j,n+1)$-rectangle (note that this is also
a $(b-a-1,b-1,n+1-a,n+1)$-rectangle).

We define the diagram $E(P)$ of $P$ in the following way. We go through the
components $Y$ of $P$ one by one, from left to right.
It is possible to fit the rectangles $\rho(Y)$ together as follows.
If a component of type $L$ is followed by a component of type $R$, the
corresponding rectangles share leftmost corners, and if a component of
type $R$ is followed by a component of type $L$, they share rightmost
corners. In each case, it is easy to see that when the rectangles are
superimposed, sharing a common leftmost or rightmost corner,
the overlapping numbers agree. The resulting diagram is defined to be $E(P)$.
It is convenient to label each rectangle by the corresponding component of
$P$; on the left hand corner if $Y$ is of type $R$, and on the right
hand corner if $Y$ is of type $L$.

We also define the diagram $E(j)$ for each $j\in [1,n]$. $E(j)$ is defined
to be the rectangle $\rho(j)$.

\noindent {\bf Example} \\
We consider the case when $n=6$ and $P=LRLL-$. Then $P$ has three components,
$L_1=L\0\0\0-$, $R_1=\0R\0\0-$ and $L_2=\0\0LL-$.
The rectangles $\rho(Y)$ are as in Figure $8$.
\begin{figure}[htbp]

\beginpicture

\setcoordinatesystem units <0.5cm,0.5cm>             
\setplotarea x from -9 to 3, y from 4 to 12       

\linethickness=0.5pt           

\put{$\rho(L_1)=$} at -3 11.5 %


\put{$1$}[c] at 3.5 13.5 %
\put{$2$}[c] at 3.5 12.5 %
\put{$3$}[c] at 1.5 11.5 %
\put{$4$}[c] at 1.5 10.5 %
\put{$5$}[c] at -0.5 9.5 %

\setlinear \plot -1 9.5  4 14.5 / %
\setlinear \plot 4 14.5  5 13.5 / %
\setlinear \plot 5 13.5  0 8.5 / %
\setlinear \plot 0 8.5  -1 9.5 / %

\put{$\rho(R_1)=$} at 11.5 11.5 %

\put{$2$}[c] at 17.5 12.5 %
\put{$3$}[c] at 15.5 11.5 %
\put{$3$}[c] at 17.5 11.5 %
\put{$4$}[c] at 15.5 10.5 %
\put{$4$}[c] at 17.5 10.5 %
\put{$5$}[c] at 15.5 9.5 %
\put{$5$}[c] at 13.5 9.5 %
\put{$6$}[c] at 15.5 8.5 %

\setlinear \plot 13 9.5  17 13.5 / %
\setlinear \plot 17 13.5  19 11.5 / %
\setlinear \plot 19 11.5  15 7.5 / %
\setlinear \plot 15 7.5  13 9.5 / %


\put{$\rho(L_2)=$} at -2 3.5 %

\put{$1$}[c] at 1.5 5.5 %
\put{$2$}[c] at 1.5 4.5 %
\put{$2$}[c] at 3.5 4.5 %
\put{$3$}[c] at 1.5 3.5 %
\put{$3$}[c] at 3.5 3.5 %
\put{$4$}[c] at 3.5 2.5 %

\setlinear \plot 0 4.5  2 6.5 / %
\setlinear \plot 2 6.5  5 3.5 / %
\setlinear \plot 5 3.5  3 1.5 / %
\setlinear \plot 3 1.5  0 4.5 / %


\put{Figure $8$: The rectangles $\rho(Y)$.}[c] at 6.5 -2 %

\endpicture

\end{figure}
These rectangles fit together to give the diagram $E(P)$. The left hand
corners of $\rho(L_1),\rho(R_1)$ match, and so do the right hand corners
of $\rho(R_1),\rho(L_2)$. The diagram $E(P)$ is shown in Figure $9$.
\begin{figure}[htbp]

\beginpicture

\setcoordinatesystem units <0.5cm,0.5cm>             
\setplotarea x from -12 to 3, y from 5 to 10       

\linethickness=0.5pt           


\put{$1$}[c] at 3.5 13.5 %
\put{$2$}[c] at 3.5 12.5 %
\put{$3$}[c] at 1.5 11.5 %
\put{$4$}[c] at 1.5 10.5 %
\put{$5$}[c] at -0.5 9.5 %

\setlinear \plot -1 9.5  4 14.5 / %
\setlinear \plot 4 14.5  5 13.5 / %
\setlinear \plot 5 13.5  0 8.5 / %
\setlinear \plot 0 8.5  -1 9.5 / %

\put{$L_1$} at 5.5 13.5


\put{$2$}[c] at 3.5 12.5 %
\put{$3$}[c] at 1.5 11.5 %
\put{$3$}[c] at 3.5 11.5 %
\put{$4$}[c] at 1.5 10.5 %
\put{$4$}[c] at 3.5 10.5 %
\put{$5$}[c] at 1.5 9.5 %
\put{$5$}[c] at -0.5 9.5 %
\put{$6$}[c] at 1.5 8.5 %

\setlinear \plot -1 9.5  3 13.5 / %
\setlinear \plot 3 13.5  5 11.5 / %
\setlinear \plot 5 11.5  1 7.5 / %
\setlinear \plot 1 7.5  -1 9.5 / %

\put{$R_1$} at -1.5 9.5


\put{$1$}[c] at 1.5 13.5 %
\put{$2$}[c] at 1.5 12.5 %
\put{$2$}[c] at 3.5 12.5 %
\put{$3$}[c] at 1.5 11.5 %
\put{$3$}[c] at 3.5 11.5 %
\put{$4$}[c] at 3.5 10.5 %


\setlinear \plot 0 12.5  2 14.5 / %
\setlinear \plot 2 14.5  5 11.5 / %
\setlinear \plot 5 11.5  3 9.5 / %
\setlinear \plot 3 9.5  0 12.5 / %

\put{$L_2$} at 5.5 11.5

\put{Figure $9$: The diagram $E(P)$.}[c] at 2.5 5.5 %

\endpicture

\end{figure}
If $P$ is a partial quiver,
the sides of the rectangles $\rho(Y)$ for $Y$ a component
of $P$ divide $E(P)$ into diagonal rows of smaller rectangles called
boxes. Consider the boxes in a diagonal running from
top left to bottom right. The number of boxes in such a diagonal
will be odd for the first $p$ diagonals, starting from the top, for some
$p$, and even for the remaining diagonals (or vice versa). Let $G$ be
the line in the diagram dividing the two adjacent diagonals 
containing an even and an odd number of boxes. Similarly, we can consider
the boxes in a diagonal running from top right to bottom left.
Again, it is possible to draw a line $H$ separating the diagonals
with an odd number of boxes from those with an even number. Let $C$ be
the point of intersection of the lines $G,H$. We call $C$ the centre of
the diagram $E(P)$. Let $V$ be the vertical line through $C$. In our example,
the picture is shown in Figure $10$.
\begin{figure}[htbp]

\beginpicture

\setcoordinatesystem units <0.5cm,0.5cm>             
\setplotarea x from -12 to 3, y from 7 to 15       

\linethickness=0.5pt           


\put{$1$}[c] at 3.5 13.5 %
\put{$2$}[c] at 3.5 12.5 %
\put{$3$}[c] at 1.5 11.5 %
\put{$4$}[c] at 1.5 10.5 %
\put{$5$}[c] at -0.5 9.5 %

\setlinear \plot -1 9.5  4 14.5 / %
\setlinear \plot 4 14.5  5 13.5 / %
\setlinear \plot 5 13.5  0 8.5 / %
\setlinear \plot 0 8.5  -1 9.5 / %

\put{$L_1$} at 5.5 13.5


\put{$2$}[c] at 3.5 12.5 %
\put{$3$}[c] at 1.5 11.5 %
\put{$3$}[c] at 3.5 11.5 %
\put{$4$}[c] at 1.5 10.5 %
\put{$4$}[c] at 3.5 10.5 %
\put{$5$}[c] at 1.5 9.5 %
\put{$5$}[c] at -0.5 9.5 %
\put{$6$}[c] at 1.5 8.5 %

\setlinear \plot -1 9.5  3 13.5 / %
\setlinear \plot 3 13.5  5 11.5 / %
\setlinear \plot 5 11.5  1 7.5 / %
\setlinear \plot 1 7.5  -1 9.5 / %

\put{$R_1$} at -1.5 9.5


\put{$1$}[c] at 1.5 13.5 %
\put{$2$}[c] at 1.5 12.5 %
\put{$2$}[c] at 3.5 12.5 %
\put{$3$}[c] at 1.5 11.5 %
\put{$3$}[c] at 3.5 11.5 %
\put{$4$}[c] at 3.5 10.5 %


\setlinear \plot 0 12.5  2 14.5 / %
\setlinear \plot 2 14.5  5 11.5 / %
\setlinear \plot 5 11.5  3 9.5 / %
\setlinear \plot 3 9.5  0 12.5 / %

\put{$L_2$} at 5.5 11.5


\setdashes <2mm,1mm>          %

\setlinear \plot -1.5 14  5.5 7 / 
\put{$G$} at -0.5 14
\setlinear \plot 7.5 16  -2.5 6 / 
\put{$H$} at 6.5 16
\setlinear \plot 2 16  2 6 / 
\put{$V$} at 2.5 16
\put{$C$} at 2.5 10.5

\put{Figure $10$: The complete diagram.}[c] at 2.5 4.5 %

v\endpicture

\end{figure}
Let $X$ be an extremal left or right corner of $E(P)$. Let $R(X)$ be the
rectangle which has the given corner as a vertex and whose edges through
this point extend as far as possible in the figure (the vertex of the
rectangle opposite to the given corner point may not be explicitly shown
in the figure). The vertical line $V$ divides $R(X)$ into two parts;
let $\Phi^+(X)$ be the set of positive roots obtained from
the part of $R(X)$ on the same side of $V$ as $X$, by reading
the numbers downwards in vertical lines. Thus $i,i+1,\ldots ,j$ gives
the positive root $\alpha_i+\alpha_{i+1}+\cdots +\alpha_j$.
Let $\Phi^+(P)$ be the union of $\Phi^+(X)$ for all left and right corners
$X$. Thus in the example in Figure $9$ we have:
\begin{eqnarray*}
\Phi^+(L_1) & = & \{\alpha_1+\alpha_2\} \\
\Phi^+(L_2) & = & \{\alpha_2+\alpha_3+\alpha_4\} \\
\Phi^+(R_1) & = & \{\alpha_5, \alpha_3+\alpha_4+\alpha_5+\alpha_6\} \\
\Phi^+(X_0) & = & \{\alpha_1+\alpha_2+\alpha_3\}
\end{eqnarray*}
where $X_0$ is the only unlabelled left or right corner. Finally,
$$\Phi^+(P)=\{
\alpha_1+\alpha_2,
\alpha_2+\alpha_3+\alpha_4,
\alpha_5, \alpha_3+\alpha_4+\alpha_5+\alpha_6,
\alpha_1+\alpha_2+\alpha_3\}.$$

We also define a set $\Phi^+(j)$ for $j\in [1,n]$. $\Phi^+(j)$ is the set
of positive roots obtained from the rectangle $E(j)=\rho(j)$ by reading the
numbers downwards in vertical lines.

\begin{theorem} (Carter and Marsh). Suppose the coordinates of the vector
$S_{\mathbf{i}}^{\mathbf{j}}(\mathbf{a}(c,\mathbf{i}))$ are labelled by the
positive roots $\Phi^+$ by means of the reduced expression $\mathbf{j}$ of
$w_0$. Let $P$ be the partial quiver $P(c,{\mathbf{i}})$. Then the
coordinate of $S_{\mathbf{i}}^{\mathbf{j}}({\bf a}(c,\mathbf{i}))$
labelled by $\alpha\in\Phi^+$ is $1$ if $\alpha\in\Phi^+(P)$ and is $0$
if $\alpha\not\in\Phi^+(P)$. Also, the coordinate of
$S_{\mathbf{i}}^{\mathbf{j}}(\mathbf{a}(j,\mathbf{i}))$ labelled by
$\alpha\in\Phi^+$ is $1$ if $\alpha\in\Phi^+(j)$ and is $0$
if $\alpha\not\in\Phi^+(j)$.
\end{theorem}

The proof of this result depends on Theorem 5.17 in~\cite{me9},
together with the description of the reparametrization function
associated with Lusztig's parametrization provided in~\cite{bfz1}.
It is hoped that this proof will appear in due course.
In the present paper we shall illustrate it in type $A_4$.

\section{Regions of linearity of the function $R$ in type $A_4$}

We shall now consider in detail the case when $\mathbf{g}$ has type $A_4$.
The Dynkin diagram will be labelled as in Figure $1$.
In this case, $\mathbf{j},\mathbf{j}'$ are the reduced words for $w_0$
given by
\begin{eqnarray*}
{\mathbf{j}} & = & (1,3,2,4,1,3,2,4,1,3) \\
{\mathbf{j}'} & = & (2,4,1,3,2,4,1,3,2,4)
\end{eqnarray*}
The piecewise-linear function $R=R_{\mathbf{j}}^{\mathbf{j}'}$ can be
written as a composition of functions corresponding to a sequence of reduced
words beginning with $\mathbf{j}$ and ending with $\mathbf{j}'$ such that
consecutive words differ by a braid relation. When $s_is_j$ is replaced by
$s_js_i$ the corresponding pair of components $a,b$ is replaced by $b,a$.
When $s_is_js_i$ is replaced by $s_js_is_j$ the corresponding triple of
components $a,b,c$ is replaced by $a',b',c'$ where:
$$(a',b',c')= \left\{ \begin{array}{ll}
(b+c-a,a,b) & \mbox{\ \ \ if\ }a\leq c \\
(b,c,b+a-c) & \mbox{\ \ \ if\ }a\geq c.
\end{array}\right.$$
See~\cite[\S2]{lusztig3}.
We list in Table I a sequence of reduced words of this kind from
$\mathbf{j}$ to $\mathbf{j}'$, underlying the long braid relations which
are used. These long braid relations are denoted by the letters
$A,B,C,D,E,F,G,H,I,J$.

\baselineskip=14pt
\begin{eqnarray*}
\mbox{Table\ I\ \ \ \ \ \ \ \ \ \ \ \ \ } & & \\
\mbox{\em Sequence\ of\ reduced\ words\ for\ $w_0$.\!\!\!\!\!\!\!\!\!\!} & & \\
 & & \\
1\ \ 3\ \ 2\ \ 4\ \ 1\ \ 3\ \ 2\ \ 4\ \ 1\ \ 3\ \  & & \\
3\ \ 1\ \ 2\ \ 4\ \ 1\ \ 3\ \ 2\ \ 4\ \ 1\ \ 3\ \  & & \\
3\ \ \underline{1\ \ 2\ \ 1}\ \ 4\ \ 3\ \ 2\ \ 4\ \ 1\ \ 3\ \  & & A \\
3\ \ 2\ \ 1\ \ 2\ \ 4\ \ 3\ \ 2\ \ 4\ \ 1\ \ 3\ \  & & \\
3\ \ 2\ \ 1\ \ 4\ \ \underline{2\ \ 3\ \ 2}\ \ 4\ \ 1\ \ 3\ \  & & B \\
3\ \ 2\ \ 1\ \ 4\ \ 3\ \ 2\ \ 3\ \ 4\ \ 1\ \ 3\ \  & & \\
3\ \ 2\ \ 1\ \ 4\ \ 3\ \ 2\ \ \underline{3\ \ 4\ \ 3}\ \ 1\ \  & & C \\
3\ \ 2\ \ 1\ \ 4\ \ 3\ \ 2\ \ 4\ \ 3\ \ 4\ \ 1\ \  & & \\
3\ \ 2\ \ 1\ \ \underline{4\ \ 3\ \ 4}\ \ 2\ \ 3\ \ 4\ \ 1\ \  & & D \\
3\ \ 2\ \ 1\ \ 3\ \ 4\ \ 3\ \ 2\ \ 3\ \ 4\ \ 1\ \  & & \\
\underline{3\ \ 2\ \ 3}\ \ 1\ \ 4\ \ 3\ \ 2\ \ 3\ \ 4\ \ 1\ \  & & E \\
2\ \ 3\ \ 2\ \ 1\ \ 4\ \ \underline{3\ \ 2\ \ 3}\ \ 4\ \ 1\ \  & & F \\
2\ \ 3\ \ 2\ \ 1\ \ 4\ \ 2\ \ 3\ \ 2\ \ 4\ \ 1\ \  & & \\
2\ \ 3\ \ \underline{2\ \ 1\ \ 2}\ \ 4\ \ 3\ \ 2\ \ 4\ \ 1\ \  & & G \\
2\ \ 3\ \ 1\ \ 2\ \ 1\ \ 4\ \ 3\ \ 2\ \ 4\ \ 1\ \  & & \\
2\ \ 3\ \ 1\ \ 2\ \ 4\ \ 1\ \ 3\ \ 2\ \ 4\ \ 1\ \  & & \\
2\ \ 3\ \ 1\ \ 2\ \ 4\ \ 3\ \ 1\ \ 2\ \ 4\ \ 1\ \  & & \\
2\ \ 3\ \ 1\ \ 2\ \ 4\ \ 3\ \ \underline{1\ \ 2\ \ 1}\ \ 4\ \  & & H \\
2\ \ 3\ \ 1\ \ 2\ \ 4\ \ 3\ \ 2\ \ 1\ \ 2\ \ 4\ \  & & \\
2\ \ 3\ \ 1\ \ 4\ \ \underline{2\ \ 3\ \ 2}\ \ 1\ \ 2\ \ 4\ \  & & I \\
2\ \ 3\ \ 1\ \ 4\ \ 3\ \ 2\ \ 3\ \ 1\ \ 2\ \ 4\ \  & & \\
2\ \ 1\ \ \underline{3\ \ 4\ \ 3}\ \ 2\ \ 3\ \ 1\ \ 2\ \ 4\ \  & & J \\
2\ \ 1\ \ 4\ \ 3\ \ 4\ \ 2\ \ 3\ \ 1\ \ 2\ \ 4\ \  & & \\
2\ \ 4\ \ 1\ \ 3\ \ 4\ \ 2\ \ 3\ \ 1\ \ 2\ \ 4\ \  & & \\
2\ \ 4\ \ 1\ \ 3\ \ 2\ \ 4\ \ 3\ \ 1\ \ 2\ \ 4\ \  & & \\
2\ \ 4\ \ 1\ \ 3\ \ 2\ \ 4\ \ 1\ \ 3\ \ 2\ \ 4\ \  & &
\end{eqnarray*}
\baselineskip=14pt

Let $v\in \R^{10}$ and let $R(v)$ be the image of $v$ under
$R=R_{\mathbf{j}}^{\mathbf{j}'}$. Since
$$
(a',b',c')= \left\{ \begin{array}{ll}
(b,c,b)+(c-a,a-c,0) & \mbox{\ \ \ if\ }a\leq c \\
(b,c,b)+(0,0,a-c)   & \mbox{\ \ \ if\ }a\geq c,
\end{array} \right.
$$
we consider the function $g(v)$ obtained from $v$ by applying the map
$(a,b,c)\mapsto (b,c,b)$ at each length $3$ braid relation and $(a,b)\mapsto
(b,a)$ at each length $2$ braid relation. We write
$$R(v)=g(v)+\varepsilon(v).$$
In our case we write $v=(a,b,c,d,e,f,g,h,i,j)$. By following the sequence of
steps in Table I we see that $g(v)=(c,h,e,j,g,h,i,j,g,h)$.

We also define vectors $\alpha_X$ for $X\in \{A,B,C,D,E,F,G,H,I,J\}$, by
$\alpha_X(v)=\xi-\eta$ where $\xi,\eta$ are the first and third components of
the triple in the vector to which the braid relation at $X$ is applied. These
vectors $\alpha_X$ can be read off from Table II and are as follows:
\begin{eqnarray*}
\alpha_A(v) & = & a-e \\
\alpha_B(v) & = & c-g \\
\alpha_C(v) & = & f-j \\
\alpha_D(v) & = & d-h \\
\alpha_E(v) & = & b-f \\
\alpha_F(v) & = & f-j \\
\alpha_G(v) & = & c-g \\
\alpha_H(v) & = & e-i \\
\alpha_I(v) & = & 0 \\
\alpha_J(v) & = & f-j
\end{eqnarray*}
We note that the vectors $\alpha_A,\alpha_B,\alpha_C,\alpha_D,\alpha_E,
\alpha_H$ are linearly independent and that the remaining vectors are
expressible in terms of these by
\baselineskip=14pt
$$\alpha_J=\alpha_F=\alpha_C,\ \ \ \alpha_G=\alpha_B,\ \ \ \alpha_I=0.$$

Now for each of the $10$ values of $X$ we have a choice between $2$ linear
functions. It may appear from this that the total number of regions of
Lusztig's function $R$ is $2^{10}$. This is not so, however, for some of the
systems of inequalities determining the sequence of choices may be
inconsistent. Calculation along the lines outlined in~\cite{carter3}
shows that there
are $204$ consistent choices. Again it is possible for two different
consistent sequences to give the same linear function $R$. Thus we obtain an
equivalence relation on the set of consistent sequences. Calculation shows
that there are $144$ equivalence classes. Thus the function $R$ in type
$A_4$ has $144$ regions of linearity.

Each region of linearity may be defined by a system of independent
inequalities of the form
\begin{eqnarray*}
& \sum_Xn_X\alpha_X \geq 0  & \\
\mbox{or\ \ \ \ \ } & \sum_Xn_X\alpha_X \leq 0  &
\end{eqnarray*}
where $n_X\in \N$.

\begin{eqnarray*}
\mbox{Table\ II\ \ \ \ \ \ \ \ \ \ \ \ \ \ } & & \\
\mbox{\em Determination\ of\ vector\ $g(v)$.} & & \\
\  & &  \\
a\ \ b\ \ c\ \ d\ \ e\ \ f\ \ g\ \ h\ \ i\ \ j\ \  & &  \\
b\ \ a\ \ c\ \ d\ \ e\ \ f\ \ g\ \ h\ \ i\ \ j\ \  & &  \\
b\ \ \underline{a\ \ c\ \ e}\ \ d\ \ f\ \ g\ \ h\ \ i\ \ j\ \  & & A \\
b\ \ c\ \ e\ \ c\ \ d\ \ f\ \ g\ \ h\ \ i\ \ j\ \  & &  \\
b\ \ c\ \ e\ \ d\ \ \underline{c\ \ f\ \ g}\ \ h\ \ i\ \ j\ \  & & B \\
b\ \ c\ \ e\ \ d\ \ f\ \ g\ \ f\ \ h\ \ i\ \ j\ \  & &  \\
b\ \ c\ \ e\ \ d\ \ f\ \ g\ \ \underline{f\ \ h\ \ j}\ \ i\ \  & & C \\
b\ \ c\ \ e\ \ d\ \ f\ \ g\ \ h\ \ j\ \ h\ \ i\ \  & &  \\
b\ \ c\ \ e\ \ \underline{d\ \ f\ \ h}\ \ g\ \ j\ \ h\ \ i\ \  & & D \\
b\ \ c\ \ e\ \ f\ \ h\ \ f\ \ g\ \ j\ \ h\ \ i\ \  & &  \\
\underline{b\ \ c\ \ f}\ \ e\ \ h\ \ f\ \ g\ \ j\ \ h\ \ i\ \  & & E \\
c\ \ f\ \ c\ \ e\ \ h\ \ \underline{f\ \ g\ \ j}\ \ h\ \ i\ \  & & F \\
c\ \ f\ \ c\ \ e\ \ h\ \ g\ \ j\ \ g\ \ h\ \ i\ \  & &  \\
c\ \ f\ \ \underline{c\ \ e\ \ g}\ \ h\ \ j\ \ g\ \ h\ \ i\ \  & & G \\
c\ \ f\ \ e\ \ g\ \ e\ \ h\ \ j\ \ g\ \ h\ \ i\ \  & &  \\
c\ \ f\ \ e\ \ g\ \ h\ \ e\ \ j\ \ g\ \ h\ \ i\ \  & &  \\
c\ \ f\ \ e\ \ g\ \ h\ \ j\ \ e\ \ g\ \ h\ \ i\ \  & &  \\
c\ \ f\ \ e\ \ g\ \ h\ \ j\ \ \underline{e\ \ g\ \ i}\ \ h\ \  & & H \\
c\ \ f\ \ e\ \ g\ \ h\ \ j\ \ g\ \ i\ \ g\ \ h\ \  & &  \\
c\ \ f\ \ e\ \ h\ \ \underline{g\ \ j\ \ g}\ \ i\ \ g\ \ h\ \  & & I \\
c\ \ f\ \ e\ \ h\ \ j\ \ g\ \ j\ \ i\ \ g\ \ h\ \  & &  \\
c\ \ e\ \ \underline{f\ \ h\ \ j}\ \ g\ \ j\ \ i\ \ g\ \ h\ \  & & J \\
c\ \ e\ \ h\ \ j\ \ h\ \ g\ \ j\ \ i\ \ g\ \ h\ \  & &  \\
c\ \ h\ \ e\ \ j\ \ h\ \ g\ \ j\ \ i\ \ g\ \ h\ \  & &  \\
c\ \ h\ \ e\ \ j\ \ g\ \ h\ \ j\ \ i\ \ g\ \ h\ \  & &  \\
c\ \ h\ \ e\ \ j\ \ g\ \ h\ \ i\ \ j\ \ g\ \ h\ \  & & 
\end{eqnarray*}
\baselineskip=14pt

We list the defining inequalities for the $144$ regions in Table III. We use
the following notation. The line
$$ABC,BCDEH || A,E,BCDE,BEH,AB^2CDEH$$
will denote the region:
$$\begin{array}{cccc}
\alpha_A+\alpha_B+\alpha_C\leq 0 & & & \alpha_A\geq 0 \\
\alpha_B+\alpha_C+\alpha_D+\alpha_E+\alpha_H\leq 0 & & & \alpha_E\geq 0 \\
& & & \alpha_B+\alpha_C+\alpha_D+\alpha_E\geq 0 \\
& & & \alpha_B+\alpha_E+\alpha_H\geq 0 \\
& & & \alpha_A+2\alpha_B+\alpha_C+\alpha_D+\alpha_E+\alpha_H\geq 0
\end{array}$$

The values of $R(v)$ for $v$ in each of these $144$ regions of linearity can be
obtained from Table IV. We have $R(v)=g(v)+\varepsilon(v)$ where
$g(v)=(c,h,e,j,g,h,i,j,g,h)$ and $\varepsilon(v)$ is given in Table IV.

\begin{center} Table III \\ {\em Regions of linearity of $R$}. \\ \ \\
\begin{longtable}{ccl}
Region of linearity & & Defining inequalities
\endfirsthead
Region of linearity & & Defining inequalities
\endhead
\endfoot
\endlastfoot
1   & &   A, BC, CD, E, BH   \ $\|$\    BCD \\
2   & &   A, BC, E, BH   \ $\|$\    B, CD \\
3   & &   ABC, CD, E, BH   \ $\|$\    A, BCD \\
4   & &   A, BC, E, BCDH   \ $\|$\    BCD, BH \\
5   & &   A, CD, E, BH   \ $\|$\    BC, D \\
6   & &   A, BC, CD, BEH   \ $\|$\    BCD, E \\
7   & &   A, BC, BEH   \ $\|$\    B, CD, E \\
8   & &   ABC, E, BH   \ $\|$\    A, B, CD \\
9   & &   ABC, E, BCDH   \ $\|$\    A, BCD, BH \\
10   & &   A, E, BCDH   \ $\|$\    BC, D, BH \\
11   & &   A, CD, BEH   \ $\|$\    BC, D, E \\
12   & &   A, B, C, BEH   \ $\|$\    CD, BE \\
13   & &   B, ABC, E, H   \ $\|$\    AB, CD \\
14   & &   ABC, BCD, E, H   \ $\|$\    ABCD, BH \\
15   & &   A, D, E, BCH   \ $\|$\    BC, BH \\
16   & &   A, C, D, BEH   \ $\|$\    BC, DE \\
17   & &   A, B, BEH   \ $\|$\    C, D, BE \\
18   & &   AB, C, BEH   \ $\|$\    A, CD, BE \\
19   & &   ABC, BE, H   \ $\|$\    AB, CD, E \\
20   & &   B, ABC, E   \ $\|$\    AB, CD, H \\
21   & &   ABC, BCD, E   \ $\|$\    B, ABCD, H \\
22   & &   BCD, E, H   \ $\|$\    ABC, D, BH \\
23   & &   D, E, BCH   \ $\|$\    A, BC, BH \\
24   & &   A, DE, BCH   \ $\|$\    BC, E, BH \\
25   & &   A, C, BDEH   \ $\|$\    BC, DE, BEH \\
26   & &   A, D, BEH   \ $\|$\    B, C, DE \\
27   & &   A, B, D, BEH   \ $\|$\    C, BDE \\
28   & &   AB, C, BE, H   \ $\|$\    CD, ABE \\
29   & &   B, ABC, CD, E   \ $\|$\    ABCD, H \\
30   & &   BC, D, E, H   \ $\|$\    ABC, BH \\
31   & &   A, C, DE, BH   \ $\|$\    BC, BEH \\
32   & &   A, BDEH   \ $\|$\    B, C, DE, BEH \\
33   & &   AB, BEH   \ $\|$\    A, C, D, BE \\
34   & &   ABC, BE   \ $\|$\    AB, CD, E, H \\
35   & &   BCD, E   \ $\|$\    B, ABC, D, H \\
36   & &   DE, BCH   \ $\|$\    A, BC, E, BH \\
37   & &   C, DE, BH   \ $\|$\    A, BC, BEH \\
38   & &   A, DE, BH   \ $\|$\    B, C, BEH \\
39   & &   A, B, BDEH   \ $\|$\    C, BDE, BEH \\
40   & &   AB, D, BEH   \ $\|$\    A, C, BDE \\
41   & &   AB, BE, H   \ $\|$\    C, D, ABE \\
42   & &   AB, C, BE   \ $\|$\    CD, ABE, H \\
43   & &   ABC, CD, BE   \ $\|$\    ABCD, E, H \\
44   & &   B, CD, E   \ $\|$\    ABC, D, H \\
45   & &   BC, D, E   \ $\|$\    B, ABC, H \\
46   & &   BC, DE, H   \ $\|$\    ABC, E, BH \\
47   & &   DE, BH   \ $\|$\    A, B, C, BEH \\
48   & &   AB, BDEH   \ $\|$\    A, C, BDE, BEH \\
49   & &   AB, BE   \ $\|$\    C, D, ABE, H \\
50   & &   CD, BE   \ $\|$\    ABC, D, E, H \\
51   & &   BC, DE   \ $\|$\    B, ABC, E, H \\
52   & &   B, DE, H   \ $\|$\    AB, C, BEH \\
53   & &   AB, BDE, H   \ $\|$\    C, ABDE, BEH \\
54   & &   AB, D, BE   \ $\|$\    C, ABDE, H \\
55   & &   C, D, BE   \ $\|$\    ABC, DE, H \\
56   & &   B, C, DE   \ $\|$\    ABC, BE, H \\
57   & &   B, DE   \ $\|$\    AB, C, BE, H \\
58   & &   BDE, H   \ $\|$\    AB, C, DE, BEH \\
59   & &   AB, BDE   \ $\|$\    C, BE, ABDE, H \\
60   & &   D, BE   \ $\|$\    AB, C, DE, H \\
61   & &   C, BDE   \ $\|$\    ABC, BE, DE, H \\
62   & &   BDE   \ $\|$\    AB, C, BE, DE, H \\
63   & &   B, AB, C, BE, ABE, H   \ $\|$\    CD \\
64   & &   ABC, CD, BCD, E, H, BH   \ $\|$\    ABCD \\
65   & &   BC, ABC, BCD, ABCD, E, H   \ $\|$\    BH \\
66   & &   A, C, D, DE, BH, BEH   \ $\|$\    BC \\
67   & &   A, C, BC, CD, BCD, BEH   \ $\|$\    BCDE \\
68   & &   BC, ABC, BCDE, ABCDE, H   \ $\|$\    E, BH \\
69   & &   B, C, ABC, BCDE, ABCDE   \ $\|$\    BE, H \\
70   & &   BC, ABC, BCD, ABCD, E   \ $\|$\    B, H \\
71   & &   B, AB, C, BE, ABE   \ $\|$\    CD, H \\
72   & &   B, AB, BDE, ABDE, H   \ $\|$\    C, BEH \\
73   & &   B, AB, BE, ABE, H   \ $\|$\    C, D \\
74   & &   A, D, DE, BH, BEH   \ $\|$\    B, C \\
75   & &   C, D, DE, BH, BEH   \ $\|$\    A, BC \\
76   & &   B, D, DE, H, BEH   \ $\|$\    AB, C \\
77   & &   C, BC, DE, H, BH   \ $\|$\    ABC, BEH \\
78   & &   CD, BCD, E, H, BH   \ $\|$\    ABC, D \\
79   & &   B, AB, D, BE, ABDE   \ $\|$\    C, H \\
80   & &   B, C, D, BE, DE   \ $\|$\    ABC, H \\
81   & &   C, ABC, CD, ABCD, BE   \ $\|$\    ABCDE, H \\
82   & &   C, D, BCDE, H, BEH   \ $\|$\    ABC, DE \\
83   & &   AB, D, BDE, H, BEH   \ $\|$\    C, ABDE \\
84   & &   A, C, BC, BCDEH, $\mbox{B}^2$CDEH   \ $\|$\    BCDE, BEH \\
85   & &   C, ABC, CD, ABCD, BEH   \ $\|$\    A, BCDE \\
86   & &   C, ABC, BCDE, H, A$\mbox{B}^2$CDEH   \ $\|$\    ABCDE, BEH \\
87   & &   ABC, CD, BCDE, H, BEH   \ $\|$\    ABCD, E \\
88   & &   BC, ABC, BCDE, ABCDE   \ $\|$\    B, E, H \\
89   & &   D, DE, BH, BEH   \ $\|$\    A, B, C \\
90   & &   B, AB, BDE, ABDE   \ $\|$\    C, BE, H \\
91   & &   B, AB, BE, ABE   \ $\|$\    C, D, H \\
92   & &   B, D, DE, BE   \ $\|$\    AB, C, H \\
93   & &   C, ABC, BCDE, A$\mbox{B}^2$CDE   \ $\|$\    BE, ABCDE, H \\
94   & &   C, BCDE, H, BDEH   \ $\|$\    ABC, DE, BEH \\
95   & &   D, BDE, H, BEH   \ $\|$\    AB, C, DE \\
96   & &   CD, BCDE, H, BEH   \ $\|$\    ABC, D, E \\
97   & &   C, ABC, BCDEH, A$\mbox{B}^2$CDEH   \ $\|$\    A, BCDE, BEH \\
98   & &   A, D, E   \ $\|$\    B, BC, BH, BCH \\
99   & &   B, E, H   \ $\|$\    AB, ABC, D, CD \\
100   & &   CD, E, BH   \ $\|$\    A, ABC, D, BCD \\
101   & &   A, E, BH   \ $\|$\    B, BC, D, CD \\
102   & &   A, BC, E   \ $\|$\    B, BCD, BH, BCDH \\
103   & &   ABC, CD, BEH   \ $\|$\    A, ABCD, E, BCDE \\
104   & &   ABC, BCDE, H   \ $\|$\    E, ABCDE, BEH, A$\mbox{B}^2$CDEH \\
105   & &   A, BC, BCDEH   \ $\|$\    E, BCDE, BEH, $\mbox{B}^2$CDEH \\
106   & &   C, D, BEH   \ $\|$\    A, ABC, DE, BCDE \\
107   & &   A, B, C   \ $\|$\    BE, BCDE, BEH, BCDEH \\
108   & &   B, E   \ $\|$\    AB, ABC, D, CD, H \\
109   & &   A, E   \ $\|$\    B, BC, D, BH, BCDH \\
110   & &   D, E   \ $\|$\    A, B, BC, BH, BCH \\
111   & &   A, DE   \ $\|$\    B, BC, E, BH, BCH \\
112   & &   E, BCDH   \ $\|$\    A, ABC, D, BCD, BH \\
113   & &   E, BH   \ $\|$\    A, B, ABC, D, CD \\
114   & &   ABC, E   \ $\|$\    A, B, BCD, BH, BCDH \\
115   & &   ABC, BEH   \ $\|$\    A, AB, CD, E, BE \\
116   & &   ABC, BCDEH   \ $\|$\    A, E, BCDE, BEH, A$\mbox{B}^2$CDEH \\
117   & &   BE, H   \ $\|$\    AB, ABC, D, CD, E \\
118   & &   A, BEH   \ $\|$\    B, BC, D, CD, E \\
119   & &   CD, BEH   \ $\|$\    A, ABC, D, E, BCDE \\
120   & &   D, BEH   \ $\|$\    A, AB, C, DE, BDE \\
121   & &   C, BDEH   \ $\|$\    A, ABC, DE, BCDE, BEH \\
122   & &   BCDE, H   \ $\|$\    ABC, E, DE, BEH, BDEH \\
123   & &   A, BCDEH   \ $\|$\    BC, E, DE, BEH, BDEH \\
124   & &   ABC, BCDE   \ $\|$\    ABCDE, A$\mbox{B}^2$CDE, E, BE, H \\
125   & &   AB, C   \ $\|$\    A, BE, BCDE, BEH, BCDEH \\
126   & &   A, BC   \ $\|$\    B, E, BCDE, BEH, BCDEH \\
127   & &   A, B   \ $\|$\    C, BE, BDE, BEH, BDEH \\
128   & &   DE   \ $\|$\    A, B, BC, E, BH, BCH \\
129   & &   BDEH   \ $\|$\    A, AB, C, DE, BDE, BEH \\
130   & &   AB   \ $\|$\    A, C, BE, BDE, BEH, BDEH \\
131   & &   BE   \ $\|$\    AB, ABC, D, CD, E, H \\
132   & &   BCDE   \ $\|$\    ABC, E, BE, DE, BDE, H \\
133   & &   C, BC, ABC, BCDE, ABCDE, H, BH   \ $\|$\    BEH \\
134   & &   B, C, ABC, CD, ABCD, BE, ABCDE   \ $\|$\    H \\
135   & &   B, AB, D, BDE, ABDE, H, BEH   \ $\|$\    C \\
136   & &   C, BC, D, DE, H, BH, BEH   \ $\|$\    ABC \\
137   & &   C, ABC, CD, ABCD, BCDE, H, BEH   \ $\|$\    ABCDE \\
138   & &   E   \ $\|$\    A, B, ABC, D, BCD, BH, BCDH \\
139   & &   BEH   \ $\|$\    A, AB, ABC, D, CD, E, BE \\
140   & &   BCDEH   \ $\|$\    A, ABC, E, DE, BCDE, BEH, BDEH \\
141   & &   ABC   \ $\|$\    A, AB, E, BE, BCDE, BEH, BCDEH \\
142   & &   A   \ $\|$\    B, BC, E, DE, BEH, BDEH, BCDEH \\
143   & &   C, BC, ABC, CD, BCD, ABCD, BCDE, ABCDE, H, BH, BEH   \ $\|$\ - \\
144   & &   - \ $\|$\    A, AB, ABC, E, BE, DE, BDE, BCDE, BEH, BDEH, BCDEH.
\end{longtable}
\end{center}

\begin{center} Table IV \\
{\em The linear function $\varepsilon(v)$ on each region.}

\begin{longtable}{cl}
Region of linearity & \hspace{5cm} $\varepsilon(v)$
\endfirsthead
Region of linearity & \hspace{5cm} $\varepsilon(v)$
\endhead
\endfoot
\endlastfoot
                          $1$ &

$(-\alpha_A - \alpha_E, -\alpha_B + \alpha_D - \alpha_E - \alpha_H, \alpha_A, \alpha_C + \alpha_E, -\alpha_C - \alpha_D, -\alpha_B - \alpha_C,$ \\
& \ \ \ \ \ \ \ \ $    \alpha_B + \alpha_C + \alpha_D + \alpha_H, \alpha_B + 2 \alpha_C + \alpha_D, 0, 0), $ \\

                          $2$ &

$(-\alpha_A - \alpha_E, -\alpha_B + \alpha_D - \alpha_E - \alpha_H, \alpha_A, \alpha_C + \alpha_E, 0, -\alpha_B - \alpha_C, \alpha_B + \alpha_H,$ \\
& \ \ \ \ \ \ \ \ $    \alpha_B + \alpha_C, \alpha_C + \alpha_D, 0), $ \\

                          $3$ &

$(-\alpha_E, -\alpha_B + \alpha_D - \alpha_E - \alpha_H, 0, \alpha_C + \alpha_E, -\alpha_C - \alpha_D, -\alpha_A - \alpha_B - \alpha_C,$ \\
& \ \ \ \ \ \ \ \ $    \alpha_B + \alpha_C + \alpha_D + \alpha_H, \alpha_A + \alpha_B + 2 \alpha_C + \alpha_D, 0, 0), $ \\

                          $4$ &

$(-\alpha_A - \alpha_E, \alpha_D - \alpha_E, \alpha_A, \alpha_C + \alpha_E, -\alpha_B - \alpha_C - \alpha_D - \alpha_H, -\alpha_B - \alpha_C,$ \\
& \ \ \ \ \ \ \ \ $    \alpha_B + \alpha_C + \alpha_D + \alpha_H, 2 \alpha_B + 2 \alpha_C + \alpha_D + \alpha_H, 0, 0), $ \\

                          $5$ &

$(-\alpha_A - \alpha_E, -\alpha_B + \alpha_D - \alpha_E - \alpha_H, \alpha_A, \alpha_C + \alpha_E, -\alpha_C - \alpha_D, 0,$ \\
& \ \ \ \ \ \ \ \ $    \alpha_B + \alpha_C + \alpha_D + \alpha_H, \alpha_C + \alpha_D, 0, \alpha_B + \alpha_C), $ \\

                          $6$ &

$(-\alpha_A, -\alpha_B + \alpha_D - \alpha_E - \alpha_H, \alpha_A, \alpha_C, -\alpha_C - \alpha_D, -\alpha_B - \alpha_C,$ \\
& \ \ \ \ \ \ \ \ $    \alpha_B + \alpha_C + \alpha_D + \alpha_E + \alpha_H, \alpha_B + 2 \alpha_C + \alpha_D, 0, 0), $ \\

                          $7$ &

$(-\alpha_A, -\alpha_B + \alpha_D - \alpha_E - \alpha_H, \alpha_A, \alpha_C, 0, -\alpha_B - \alpha_C, \alpha_B + \alpha_E + \alpha_H, \alpha_B + \alpha_C,$ \\
& \ \ \ \ \ \ \ \ $    \alpha_C + \alpha_D, 0), $ \\

                          $8$ &

$(-\alpha_E, -\alpha_B + \alpha_D - \alpha_E - \alpha_H, 0, \alpha_C + \alpha_E, 0, -\alpha_A - \alpha_B - \alpha_C, \alpha_B + \alpha_H,$ \\
& \ \ \ \ \ \ \ \ $    \alpha_A + \alpha_B + \alpha_C, \alpha_C + \alpha_D, 0), $ \\

                          $9$ &

$(-\alpha_E, \alpha_D - \alpha_E, 0, \alpha_C + \alpha_E, -\alpha_B - \alpha_C - \alpha_D - \alpha_H, -\alpha_A - \alpha_B - \alpha_C,$ \\
& \ \ \ \ \ \ \ \ $    \alpha_B + \alpha_C + \alpha_D + \alpha_H, \alpha_A + 2 \alpha_B + 2 \alpha_C + \alpha_D + \alpha_H, 0, 0), $ \\

                          $10$ &

$(-\alpha_A - \alpha_E, \alpha_D - \alpha_E, \alpha_A, \alpha_C + \alpha_E, -\alpha_B - \alpha_C - \alpha_D - \alpha_H, 0,$ \\
& \ \ \ \ \ \ \ \ $    \alpha_B + \alpha_C + \alpha_D + \alpha_H, \alpha_B + \alpha_C + \alpha_D + \alpha_H, 0, \alpha_B + \alpha_C), $ \\

                          $11$ &

$(-\alpha_A, -\alpha_B + \alpha_D - \alpha_E - \alpha_H, \alpha_A, \alpha_C, -\alpha_C - \alpha_D, 0,$ \\
& \ \ \ \ \ \ \ \ $    \alpha_B + \alpha_C + \alpha_D + \alpha_E + \alpha_H, \alpha_C + \alpha_D, 0, \alpha_B + \alpha_C), $ \\

                          $12$ &

$(-\alpha_A, -\alpha_B + \alpha_D - \alpha_E - \alpha_H, \alpha_A, -\alpha_B + \alpha_C, \alpha_B, -\alpha_C, \alpha_B + \alpha_E + \alpha_H, \alpha_C,$ \\
& \ \ \ \ \ \ \ \ $    \alpha_C + \alpha_D, 0), $ \\

                          $13$ &

$(-\alpha_E, -\alpha_B + \alpha_D - \alpha_E - \alpha_H, -\alpha_B, \alpha_C + \alpha_E, \alpha_B, -\alpha_A - \alpha_B - \alpha_C, \alpha_H,$ \\
& \ \ \ \ \ \ \ \ $    \alpha_A + \alpha_B + \alpha_C, \alpha_C + \alpha_D, 0), $ \\

                          $14$ &

$(-\alpha_E, \alpha_D - \alpha_E, -\alpha_B - \alpha_C - \alpha_D, \alpha_C + \alpha_E, -\alpha_H, -\alpha_A - \alpha_B - \alpha_C, \alpha_H,$ \\
& \ \ \ \ \ \ \ \ $    \alpha_A + 2 \alpha_B + 2 \alpha_C + \alpha_D + \alpha_H, 0, 0), $ \\

                          $15$ &

$(-\alpha_A - \alpha_D - \alpha_E, \alpha_D - \alpha_E, \alpha_A, \alpha_C + \alpha_E, -\alpha_B - \alpha_C - \alpha_H, \alpha_D, \alpha_B + \alpha_C + \alpha_H,$ \\
& \ \ \ \ \ \ \ \ $    \alpha_B + \alpha_C + \alpha_H, 0, \alpha_B + \alpha_C), $ \\

                          $16$ &

$(-\alpha_A, -\alpha_B + \alpha_D - \alpha_E - \alpha_H, \alpha_A, \alpha_C - \alpha_D, -\alpha_C, \alpha_D,$ \\
& \ \ \ \ \ \ \ \ $    \alpha_B + \alpha_C + \alpha_D + \alpha_E + \alpha_H, \alpha_C, 0, \alpha_B + \alpha_C), $ \\

                          $17$ &

$(-\alpha_A, -\alpha_B + \alpha_D - \alpha_E - \alpha_H, \alpha_A, -\alpha_B + \alpha_C, \alpha_B, 0, \alpha_B + \alpha_E + \alpha_H, 0,$ \\
& \ \ \ \ \ \ \ \ $    \alpha_C + \alpha_D, \alpha_C), $ \\

                          $18$ &

$(0, -\alpha_B + \alpha_D - \alpha_E - \alpha_H, 0, -\alpha_A - \alpha_B + \alpha_C, \alpha_A + \alpha_B, -\alpha_C, \alpha_B + \alpha_E + \alpha_H,$ \\
& \ \ \ \ \ \ \ \ $    \alpha_C, \alpha_C + \alpha_D, 0), $ \\

                          $19$ &

$(0, -\alpha_B + \alpha_D - \alpha_E - \alpha_H, -\alpha_B - \alpha_E, \alpha_C, \alpha_B + \alpha_E, -\alpha_A - \alpha_B - \alpha_C, \alpha_H,$ \\
& \ \ \ \ \ \ \ \ $    \alpha_A + \alpha_B + \alpha_C, \alpha_C + \alpha_D, 0), $ \\

                          $20$ &

$(-\alpha_E, -\alpha_B + \alpha_D - \alpha_E, -\alpha_B, \alpha_C + \alpha_E, \alpha_B, -\alpha_A - \alpha_B - \alpha_C, 0, \alpha_A + \alpha_B + \alpha_C,$ \\
& \ \ \ \ \ \ \ \ $    \alpha_C + \alpha_D + \alpha_H, 0), $ \\

                          $21$ &

$(-\alpha_E, \alpha_D - \alpha_E, -\alpha_B - \alpha_C - \alpha_D, \alpha_C + \alpha_E, 0, -\alpha_A - \alpha_B - \alpha_C, 0,$ \\
& \ \ \ \ \ \ \ \ $    \alpha_A + 2 \alpha_B + 2 \alpha_C + \alpha_D, \alpha_H, 0), $ \\

                          $22$ &

$(-\alpha_E, \alpha_D - \alpha_E, -\alpha_B - \alpha_C - \alpha_D, \alpha_C + \alpha_E, -\alpha_H, 0, \alpha_H, \alpha_B + \alpha_C + \alpha_D + \alpha_H,$ \\
& \ \ \ \ \ \ \ \ $    0, \alpha_A + \alpha_B + \alpha_C), $ \\

                          $23$ &

$(-\alpha_D - \alpha_E, \alpha_D - \alpha_E, 0, \alpha_C + \alpha_E, -\alpha_B - \alpha_C - \alpha_H, \alpha_D, \alpha_B + \alpha_C + \alpha_H,$ \\
& \ \ \ \ \ \ \ \ $    \alpha_B + \alpha_C + \alpha_H, 0, \alpha_A + \alpha_B + \alpha_C), $ \\

                          $24$ &

$(-\alpha_A - \alpha_D - \alpha_E, \alpha_D, \alpha_A, \alpha_C, -\alpha_B - \alpha_C - \alpha_H, \alpha_D + \alpha_E, \alpha_B + \alpha_C + \alpha_H,$ \\
& \ \ \ \ \ \ \ \ $    \alpha_B + \alpha_C + \alpha_H, 0, \alpha_B + \alpha_C), $ \\

                          $25$ &

$(-\alpha_A, \alpha_D, \alpha_A, -\alpha_B + \alpha_C - \alpha_D - \alpha_E - \alpha_H, -\alpha_C, \alpha_B + \alpha_D + \alpha_E + \alpha_H,$ \\
& \ \ \ \ \ \ \ \ $    \alpha_B + \alpha_C + \alpha_D + \alpha_E + \alpha_H, \alpha_C, 0, \alpha_B + \alpha_C), $ \\

                          $26$ &

$(-\alpha_A, -\alpha_B + \alpha_D - \alpha_E - \alpha_H, \alpha_A, \alpha_C - \alpha_D, 0, \alpha_D, \alpha_B + \alpha_D + \alpha_E + \alpha_H, 0,$ \\
& \ \ \ \ \ \ \ \ $    \alpha_C, \alpha_B + \alpha_C), $ \\

                          $27$ &

$(-\alpha_A, -\alpha_B + \alpha_D - \alpha_E - \alpha_H, \alpha_A, -\alpha_B + \alpha_C - \alpha_D, \alpha_B, \alpha_D,$ \\
& \ \ \ \ \ \ \ \ $    \alpha_B + \alpha_D + \alpha_E + \alpha_H, 0, \alpha_C, \alpha_C), $ \\

                          $28$ &

$(0, -\alpha_B + \alpha_D - \alpha_E - \alpha_H, -\alpha_B - \alpha_E, -\alpha_A - \alpha_B + \alpha_C, \alpha_A + 2 \alpha_B + \alpha_E, -\alpha_C,$ \\
& \ \ \ \ \ \ \ \ $    \alpha_H, \alpha_C, \alpha_C + \alpha_D, 0), $ \\

                          $29$ &

$(-\alpha_E, -\alpha_B + \alpha_D - \alpha_E, -\alpha_B - \alpha_C - \alpha_D, \alpha_C + \alpha_E, \alpha_B, -\alpha_A - \alpha_B - \alpha_C, 0,$ \\
& \ \ \ \ \ \ \ \ $    \alpha_A + \alpha_B + 2 \alpha_C + \alpha_D, \alpha_H, 0), $ \\

                          $30$ &

$(-\alpha_D - \alpha_E, \alpha_D - \alpha_E, -\alpha_B - \alpha_C, \alpha_C + \alpha_E, -\alpha_H, \alpha_D, \alpha_H, \alpha_B + \alpha_C + \alpha_H, 0,$ \\
& \ \ \ \ \ \ \ \ $    \alpha_A + \alpha_B + \alpha_C), $ \\

                          $31$ &

$(-\alpha_A - \alpha_D - \alpha_E, \alpha_D, \alpha_A, -\alpha_B + \alpha_C - \alpha_H, -\alpha_C, \alpha_B + \alpha_D + \alpha_E + \alpha_H,$ \\
& \ \ \ \ \ \ \ \ $    \alpha_B + \alpha_C + \alpha_H, \alpha_C, 0, \alpha_B + \alpha_C), $ \\

                          $32$ &

$(-\alpha_A, \alpha_D, \alpha_A, -\alpha_B + \alpha_C - \alpha_D - \alpha_E - \alpha_H, 0, \alpha_B + \alpha_D + \alpha_E + \alpha_H,$ \\
& \ \ \ \ \ \ \ \ $    \alpha_B + \alpha_D + \alpha_E + \alpha_H, 0, \alpha_C, \alpha_B + \alpha_C), $ \\

                          $33$ &

$(0, -\alpha_B + \alpha_D - \alpha_E - \alpha_H, 0, -\alpha_A - \alpha_B + \alpha_C, \alpha_A + \alpha_B, 0, \alpha_B + \alpha_E + \alpha_H,$ \\
& \ \ \ \ \ \ \ \ $    0, \alpha_C + \alpha_D, \alpha_C), $ \\

                          $34$ &

$(0, -\alpha_B + \alpha_D - \alpha_E, -\alpha_B - \alpha_E, \alpha_C, \alpha_B + \alpha_E, -\alpha_A - \alpha_B - \alpha_C, 0,$ \\
& \ \ \ \ \ \ \ \ $    \alpha_A + \alpha_B + \alpha_C, \alpha_C + \alpha_D + \alpha_H, 0), $ \\

                          $35$ &

$(-\alpha_E, \alpha_D - \alpha_E, -\alpha_B - \alpha_C - \alpha_D, \alpha_C + \alpha_E, 0, 0, 0, \alpha_B + \alpha_C + \alpha_D, \alpha_H,$ \\
& \ \ \ \ \ \ \ \ $    \alpha_A + \alpha_B + \alpha_C), $ \\

                          $36$ &

$(-\alpha_D - \alpha_E, \alpha_D, 0, \alpha_C, -\alpha_B - \alpha_C - \alpha_H, \alpha_D + \alpha_E, \alpha_B + \alpha_C + \alpha_H,$ \\
& \ \ \ \ \ \ \ \ $    \alpha_B + \alpha_C + \alpha_H, 0, \alpha_A + \alpha_B + \alpha_C), $ \\

                          $37$ &

$(-\alpha_D - \alpha_E, \alpha_D, 0, -\alpha_B + \alpha_C - \alpha_H, -\alpha_C, \alpha_B + \alpha_D + \alpha_E + \alpha_H,$ \\
& \ \ \ \ \ \ \ \ $    \alpha_B + \alpha_C + \alpha_H, \alpha_C, 0, \alpha_A + \alpha_B + \alpha_C), $ \\

                          $38$ &

$(-\alpha_A - \alpha_D - \alpha_E, \alpha_D, \alpha_A, -\alpha_B + \alpha_C - \alpha_H, 0, \alpha_B + \alpha_D + \alpha_E + \alpha_H, \alpha_B + \alpha_H,$ \\
& \ \ \ \ \ \ \ \ $    0, \alpha_C, \alpha_B + \alpha_C), $ \\

                          $39$ &

$(-\alpha_A, \alpha_D, \alpha_A, -2 \alpha_B + \alpha_C - \alpha_D - \alpha_E - \alpha_H, \alpha_B, \alpha_B + \alpha_D + \alpha_E + \alpha_H,$ \\
& \ \ \ \ \ \ \ \ $    \alpha_B + \alpha_D + \alpha_E + \alpha_H, 0, \alpha_C, \alpha_C), $ \\

                          $40$ &

$(0, -\alpha_B + \alpha_D - \alpha_E - \alpha_H, 0, -\alpha_A - \alpha_B + \alpha_C - \alpha_D, \alpha_A + \alpha_B, \alpha_D,$ \\
& \ \ \ \ \ \ \ \ $    \alpha_B + \alpha_D + \alpha_E + \alpha_H, 0, \alpha_C, \alpha_C), $ \\

                          $41$ &

$(0, -\alpha_B + \alpha_D - \alpha_E - \alpha_H, -\alpha_B - \alpha_E, -\alpha_A - \alpha_B + \alpha_C, \alpha_A + 2 \alpha_B + \alpha_E, 0,$ \\
& \ \ \ \ \ \ \ \ $    \alpha_H, 0, \alpha_C + \alpha_D, \alpha_C), $ \\

                          $42$ &

$(0, -\alpha_B + \alpha_D - \alpha_E, -\alpha_B - \alpha_E, -\alpha_A - \alpha_B + \alpha_C, \alpha_A + 2 \alpha_B + \alpha_E, -\alpha_C, 0,$ \\
& \ \ \ \ \ \ \ \ $    \alpha_C, \alpha_C + \alpha_D + \alpha_H, 0), $ \\

                          $43$ &

$(0, -\alpha_B + \alpha_D - \alpha_E, -\alpha_B - \alpha_C - \alpha_D - \alpha_E, \alpha_C, \alpha_B + \alpha_E, -\alpha_A - \alpha_B - \alpha_C, 0,$ \\
& \ \ \ \ \ \ \ \ $    \alpha_A + \alpha_B + 2 \alpha_C + \alpha_D, \alpha_H, 0), $ \\

                          $44$ &

$(-\alpha_E, -\alpha_B + \alpha_D - \alpha_E, -\alpha_B - \alpha_C - \alpha_D, \alpha_C + \alpha_E, \alpha_B, 0, 0, \alpha_C + \alpha_D, \alpha_H,$ \\
& \ \ \ \ \ \ \ \ $    \alpha_A + \alpha_B + \alpha_C), $ \\

                          $45$ &

$(-\alpha_D - \alpha_E, \alpha_D - \alpha_E, -\alpha_B - \alpha_C, \alpha_C + \alpha_E, 0, \alpha_D, 0, \alpha_B + \alpha_C, \alpha_H,$ \\
& \ \ \ \ \ \ \ \ $    \alpha_A + \alpha_B + \alpha_C), $ \\

                          $46$ &

$(-\alpha_D - \alpha_E, \alpha_D, -\alpha_B - \alpha_C, \alpha_C, -\alpha_H, \alpha_D + \alpha_E, \alpha_H, \alpha_B + \alpha_C + \alpha_H, 0,$ \\
& \ \ \ \ \ \ \ \ $    \alpha_A + \alpha_B + \alpha_C), $ \\

                          $47$ &

$(-\alpha_D - \alpha_E, \alpha_D, 0, -\alpha_B + \alpha_C - \alpha_H, 0, \alpha_B + \alpha_D + \alpha_E + \alpha_H, \alpha_B + \alpha_H, 0,$ \\
& \ \ \ \ \ \ \ \ $    \alpha_C, \alpha_A + \alpha_B + \alpha_C), $ \\

                          $48$ &

$(0, \alpha_D, 0, -\alpha_A - 2 \alpha_B + \alpha_C - \alpha_D - \alpha_E - \alpha_H, \alpha_A + \alpha_B,$ \\
& \ \ \ \ \ \ \ \ $    \alpha_B + \alpha_D + \alpha_E + \alpha_H, \alpha_B + \alpha_D + \alpha_E + \alpha_H, 0, \alpha_C, \alpha_C), $ \\

                          $49$ &

$(0, -\alpha_B + \alpha_D - \alpha_E, -\alpha_B - \alpha_E, -\alpha_A - \alpha_B + \alpha_C, \alpha_A + 2 \alpha_B + \alpha_E, 0, 0, 0,$ \\
& \ \ \ \ \ \ \ \ $    \alpha_C + \alpha_D + \alpha_H, \alpha_C), $ \\

                          $50$ &

$(0, -\alpha_B + \alpha_D - \alpha_E, -\alpha_B - \alpha_C - \alpha_D - \alpha_E, \alpha_C, \alpha_B + \alpha_E, 0, 0, \alpha_C + \alpha_D,$ \\
& \ \ \ \ \ \ \ \ $    \alpha_H, \alpha_A + \alpha_B + \alpha_C), $ \\

                          $51$ &

$(-\alpha_D - \alpha_E, \alpha_D, -\alpha_B - \alpha_C, \alpha_C, 0, \alpha_D + \alpha_E, 0, \alpha_B + \alpha_C, \alpha_H,$ \\
& \ \ \ \ \ \ \ \ $    \alpha_A + \alpha_B + \alpha_C), $ \\

                          $52$ &

$(-\alpha_D - \alpha_E, \alpha_D, -\alpha_B, -\alpha_B + \alpha_C - \alpha_H, \alpha_B, \alpha_B + \alpha_D + \alpha_E + \alpha_H, \alpha_H, 0, \alpha_C,$ \\
& \ \ \ \ \ \ \ \ $    \alpha_A + \alpha_B + \alpha_C), $ \\

                          $53$ &

$(0, \alpha_D, -\alpha_B - \alpha_D - \alpha_E, -\alpha_A - 2 \alpha_B + \alpha_C - \alpha_D - \alpha_E - \alpha_H,$ \\
& \ \ \ \ \ \ \ \ $    \alpha_A + 2 \alpha_B + \alpha_D + \alpha_E, \alpha_B + \alpha_D + \alpha_E + \alpha_H, \alpha_H, 0, \alpha_C, \alpha_C), $ \\

                          $54$ &

$(0, -\alpha_B + \alpha_D - \alpha_E, -\alpha_B - \alpha_D - \alpha_E, -\alpha_A - \alpha_B + \alpha_C - \alpha_D,$ \\
& \ \ \ \ \ \ \ \ $    \alpha_A + 2 \alpha_B + \alpha_D + \alpha_E, \alpha_D, 0, 0, \alpha_C + \alpha_H, \alpha_C), $ \\

                          $55$ &

$(0, -\alpha_B + \alpha_D - \alpha_E, -\alpha_B - \alpha_C - \alpha_D - \alpha_E, \alpha_C - \alpha_D, \alpha_B + \alpha_D + \alpha_E, \alpha_D,$ \\
& \ \ \ \ \ \ \ \ $    0, \alpha_C, \alpha_H, \alpha_A + \alpha_B + \alpha_C), $ \\

                          $56$ &

$(-\alpha_D - \alpha_E, \alpha_D, -\alpha_B - \alpha_C, -\alpha_B + \alpha_C, \alpha_B, \alpha_B + \alpha_D + \alpha_E, 0, \alpha_C, \alpha_H,$ \\
& \ \ \ \ \ \ \ \ $    \alpha_A + \alpha_B + \alpha_C), $ \\

                          $57$ &

$(-\alpha_D - \alpha_E, \alpha_D, -\alpha_B, -\alpha_B + \alpha_C, \alpha_B, \alpha_B + \alpha_D + \alpha_E, 0, 0, \alpha_C + \alpha_H,$ \\
& \ \ \ \ \ \ \ \ $    \alpha_A + \alpha_B + \alpha_C), $ \\

                          $58$ &

$(0, \alpha_D, -\alpha_B - \alpha_D - \alpha_E, -\alpha_B + \alpha_C - \alpha_D - \alpha_E - \alpha_H, \alpha_B + \alpha_D + \alpha_E,$ \\
& \ \ \ \ \ \ \ \ $    \alpha_B + \alpha_D + \alpha_E + \alpha_H, \alpha_H, 0, \alpha_C, \alpha_A + \alpha_B + \alpha_C), $ \\

                          $59$ &

$(0, \alpha_D, -\alpha_B - \alpha_D - \alpha_E, -\alpha_A - 2 \alpha_B + \alpha_C - \alpha_D - \alpha_E,$ \\
& \ \ \ \ \ \ \ \ $    \alpha_A + 2 \alpha_B + \alpha_D + \alpha_E, \alpha_B + \alpha_D + \alpha_E, 0, 0, \alpha_C + \alpha_H, \alpha_C), $ \\

                          $60$ &

$(0, -\alpha_B + \alpha_D - \alpha_E, -\alpha_B - \alpha_D - \alpha_E, \alpha_C - \alpha_D, \alpha_B + \alpha_D + \alpha_E, \alpha_D, 0,$ \\
& \ \ \ \ \ \ \ \ $    0, \alpha_C + \alpha_H, \alpha_A + \alpha_B + \alpha_C), $ \\

                          $61$ &

$(0, \alpha_D, -\alpha_B - \alpha_C - \alpha_D - \alpha_E, -\alpha_B + \alpha_C - \alpha_D - \alpha_E, \alpha_B + \alpha_D + \alpha_E,$ \\
& \ \ \ \ \ \ \ \ $    \alpha_B + \alpha_D + \alpha_E, 0, \alpha_C, \alpha_H, \alpha_A + \alpha_B + \alpha_C), $ \\

                          $62$ &

$(0, \alpha_D, -\alpha_B - \alpha_D - \alpha_E, -\alpha_B + \alpha_C - \alpha_D - \alpha_E, \alpha_B + \alpha_D + \alpha_E,$ \\
& \ \ \ \ \ \ \ \ $    \alpha_B + \alpha_D + \alpha_E, 0, 0, \alpha_C + \alpha_H, \alpha_A + \alpha_B + \alpha_C), $ \\

                          $63$ &

$(-\alpha_A - \alpha_B - \alpha_E, -\alpha_B + \alpha_D - \alpha_E - \alpha_H, \alpha_A, \alpha_C + \alpha_E, \alpha_B, -\alpha_C, \alpha_H, \alpha_C,$ \\
& \ \ \ \ \ \ \ \ $    \alpha_C + \alpha_D, 0), $ \\

                          $64$ &

$(-\alpha_E, -\alpha_B + \alpha_D - \alpha_E - \alpha_H, -\alpha_B - \alpha_C - \alpha_D, \alpha_C + \alpha_E, \alpha_B, -\alpha_A - \alpha_B - \alpha_C,$ \\
& \ \ \ \ \ \ \ \ $    \alpha_H, \alpha_A + \alpha_B + 2 \alpha_C + \alpha_D, 0, 0), $ \\

                          $65$ &

$(-\alpha_A - \alpha_B - \alpha_C - \alpha_D - \alpha_E, \alpha_D - \alpha_E, \alpha_A, \alpha_C + \alpha_E, -\alpha_H, \alpha_D, \alpha_H,$ \\
& \ \ \ \ \ \ \ \ $    \alpha_B + \alpha_C + \alpha_H, 0, 0), $ \\

                          $66$ &

$(-\alpha_A - \alpha_D - \alpha_E, -\alpha_B + \alpha_D - \alpha_E - \alpha_H, \alpha_A, \alpha_C + \alpha_E, -\alpha_C, \alpha_D,$ \\
& \ \ \ \ \ \ \ \ $    \alpha_B + \alpha_C + \alpha_H, \alpha_C, 0, \alpha_B + \alpha_C), $ \\

                          $67$ &

$(-\alpha_A, -\alpha_B + \alpha_D - \alpha_E - \alpha_H, \alpha_A, -\alpha_B - \alpha_D, \alpha_B, \alpha_D,$ \\
& \ \ \ \ \ \ \ \ $    \alpha_B + \alpha_C + \alpha_D + \alpha_E + \alpha_H, \alpha_C, 0, 0), $ \\

                          $68$ &

$(-\alpha_A - \alpha_B - \alpha_C - \alpha_D - \alpha_E, \alpha_D, \alpha_A, \alpha_C, -\alpha_H, \alpha_D + \alpha_E, \alpha_H, \alpha_B + \alpha_C + \alpha_H,$ \\
& \ \ \ \ \ \ \ \ $    0, 0), $ \\

                          $69$ &

$(-\alpha_A - \alpha_B - \alpha_C - \alpha_D - \alpha_E, \alpha_D, \alpha_A, -\alpha_B + \alpha_C, \alpha_B, \alpha_B + \alpha_D + \alpha_E, 0, \alpha_C,$ \\
& \ \ \ \ \ \ \ \ $    \alpha_H, 0), $ \\

                          $70$ &

$(-\alpha_A - \alpha_B - \alpha_C - \alpha_D - \alpha_E, \alpha_D - \alpha_E, \alpha_A, \alpha_C + \alpha_E, 0, \alpha_D, 0, \alpha_B + \alpha_C,$ \\
& \ \ \ \ \ \ \ \ $    \alpha_H, 0), $ \\

                          $71$ &

$(-\alpha_A - \alpha_B - \alpha_E, -\alpha_B + \alpha_D - \alpha_E, \alpha_A, \alpha_C + \alpha_E, \alpha_B, -\alpha_C, 0, \alpha_C,$ \\
& \ \ \ \ \ \ \ \ $    \alpha_C + \alpha_D + \alpha_H, 0), $ \\

                          $72$ &

$(-\alpha_A - \alpha_B - \alpha_D - \alpha_E, \alpha_D, \alpha_A, -\alpha_B + \alpha_C - \alpha_H, \alpha_B, \alpha_B + \alpha_D + \alpha_E + \alpha_H, \alpha_H,$ \\
& \ \ \ \ \ \ \ \ $    0, \alpha_C, \alpha_C), $ \\

                          $73$ &

$(-\alpha_A - \alpha_B - \alpha_E, -\alpha_B + \alpha_D - \alpha_E - \alpha_H, \alpha_A, \alpha_C + \alpha_E, \alpha_B, 0, \alpha_H, 0,$ \\
& \ \ \ \ \ \ \ \ $    \alpha_C + \alpha_D, \alpha_C), $ \\

                          $74$ &

$(-\alpha_A - \alpha_D - \alpha_E, -\alpha_B + \alpha_D - \alpha_E - \alpha_H, \alpha_A, \alpha_C + \alpha_E, 0, \alpha_D, \alpha_B + \alpha_H, 0,$ \\
& \ \ \ \ \ \ \ \ $    \alpha_C, \alpha_B + \alpha_C), $ \\

                          $75$ &

$(-\alpha_D - \alpha_E, -\alpha_B + \alpha_D - \alpha_E - \alpha_H, 0, \alpha_C + \alpha_E, -\alpha_C, \alpha_D, \alpha_B + \alpha_C + \alpha_H, \alpha_C,$ \\
& \ \ \ \ \ \ \ \ $    0, \alpha_A + \alpha_B + \alpha_C), $ \\

                          $76$ &

$(-\alpha_D - \alpha_E, -\alpha_B + \alpha_D - \alpha_E - \alpha_H, -\alpha_B, \alpha_C + \alpha_E, \alpha_B, \alpha_D, \alpha_H, 0, \alpha_C,$ \\
& \ \ \ \ \ \ \ \ $    \alpha_A + \alpha_B + \alpha_C), $ \\

                          $77$ &

$(-\alpha_D - \alpha_E, \alpha_D, -\alpha_B - \alpha_C, -\alpha_B + \alpha_C - \alpha_H, \alpha_B, \alpha_B + \alpha_D + \alpha_E + \alpha_H, \alpha_H, \alpha_C,$ \\
& \ \ \ \ \ \ \ \ $    0, \alpha_A + \alpha_B + \alpha_C), $ \\

                          $78$ &

$(-\alpha_E, -\alpha_B + \alpha_D - \alpha_E - \alpha_H, -\alpha_B - \alpha_C - \alpha_D, \alpha_C + \alpha_E, \alpha_B, 0, \alpha_H, \alpha_C + \alpha_D,$ \\
& \ \ \ \ \ \ \ \ $    0, \alpha_A + \alpha_B + \alpha_C), $ \\

                          $79$ &

$(-\alpha_A - \alpha_B - \alpha_D - \alpha_E, -\alpha_B + \alpha_D - \alpha_E, \alpha_A, \alpha_C + \alpha_E, \alpha_B, \alpha_D, 0, 0,$ \\
& \ \ \ \ \ \ \ \ $    \alpha_C + \alpha_H, \alpha_C), $ \\

                          $80$ &

$(-\alpha_D - \alpha_E, -\alpha_B + \alpha_D - \alpha_E, -\alpha_B - \alpha_C, \alpha_C + \alpha_E, \alpha_B, \alpha_D, 0, \alpha_C, \alpha_H,$ \\
& \ \ \ \ \ \ \ \ $    \alpha_A + \alpha_B + \alpha_C), $ \\

                          $81$ &

$(0, -\alpha_B + \alpha_D - \alpha_E, -\alpha_B - \alpha_C - \alpha_D - \alpha_E, -\alpha_A - \alpha_B - \alpha_D,$ \\
& \ \ \ \ \ \ \ \ $    \alpha_A + 2 \alpha_B + \alpha_C + \alpha_D + \alpha_E, \alpha_D, 0, \alpha_C, \alpha_H, 0), $ \\

                          $82$ &

$(0, -\alpha_B + \alpha_D - \alpha_E - \alpha_H, -\alpha_B - \alpha_C - \alpha_D - \alpha_E, \alpha_C - \alpha_D, \alpha_B + \alpha_D + \alpha_E,$ \\
& \ \ \ \ \ \ \ \ $    \alpha_D, \alpha_H, \alpha_C, 0, \alpha_A + \alpha_B + \alpha_C), $ \\

                          $83$ &

$(0, -\alpha_B + \alpha_D - \alpha_E - \alpha_H, -\alpha_B - \alpha_D - \alpha_E, -\alpha_A - \alpha_B + \alpha_C - \alpha_D,$ \\
& \ \ \ \ \ \ \ \ $    \alpha_A + 2 \alpha_B + \alpha_D + \alpha_E, \alpha_D, \alpha_H, 0, \alpha_C, \alpha_C), $ \\

                          $84$ &

$(-\alpha_A, \alpha_D, \alpha_A, -2 \alpha_B - \alpha_D - \alpha_E - \alpha_H, \alpha_B, \alpha_B + \alpha_D + \alpha_E + \alpha_H,$ \\
& \ \ \ \ \ \ \ \ $    \alpha_B + \alpha_C + \alpha_D + \alpha_E + \alpha_H, \alpha_C, 0, 0), $ \\

                          $85$ &

$(0, -\alpha_B + \alpha_D - \alpha_E - \alpha_H, 0, -\alpha_A - \alpha_B - \alpha_D, \alpha_A + \alpha_B, \alpha_D,$ \\
& \ \ \ \ \ \ \ \ $    \alpha_B + \alpha_C + \alpha_D + \alpha_E + \alpha_H, \alpha_C, 0, 0), $ \\

                          $86$ &

$(0, \alpha_D, -\alpha_B - \alpha_C - \alpha_D - \alpha_E, -\alpha_A - 2 \alpha_B - \alpha_D - \alpha_E - \alpha_H,$ \\
& \ \ \ \ \ \ \ \ $    \alpha_A + 2 \alpha_B + \alpha_C + \alpha_D + \alpha_E, \alpha_B + \alpha_D + \alpha_E + \alpha_H, \alpha_H, \alpha_C, 0, 0),$ \\

                          $87$ &

$(0, -\alpha_B + \alpha_D - \alpha_E - \alpha_H, -\alpha_B - \alpha_C - \alpha_D - \alpha_E, \alpha_C, \alpha_B + \alpha_E,$ \\
& \ \ \ \ \ \ \ \ $    -\alpha_A - \alpha_B - \alpha_C, \alpha_H, \alpha_A + \alpha_B + 2 \alpha_C + \alpha_D, 0, 0), $ \\

                          $88$ &

$(-\alpha_A - \alpha_B - \alpha_C - \alpha_D - \alpha_E, \alpha_D, \alpha_A, \alpha_C, 0, \alpha_D + \alpha_E, 0, \alpha_B + \alpha_C, \alpha_H, 0
    ), $ \\

                          $89$ &

$(-\alpha_D - \alpha_E, -\alpha_B + \alpha_D - \alpha_E - \alpha_H, 0, \alpha_C + \alpha_E, 0, \alpha_D, \alpha_B + \alpha_H, 0, \alpha_C,$ \\
& \ \ \ \ \ \ \ \ $    \alpha_A + \alpha_B + \alpha_C), $ \\

                          $90$ &

$(-\alpha_A - \alpha_B - \alpha_D - \alpha_E, \alpha_D, \alpha_A, -\alpha_B + \alpha_C, \alpha_B, \alpha_B + \alpha_D + \alpha_E, 0, 0,$ \\
& \ \ \ \ \ \ \ \ $    \alpha_C + \alpha_H, \alpha_C), $ \\

                          $91$ &

$(-\alpha_A - \alpha_B - \alpha_E, -\alpha_B + \alpha_D - \alpha_E, \alpha_A, \alpha_C + \alpha_E, \alpha_B, 0, 0, 0,$ \\
& \ \ \ \ \ \ \ \ $    \alpha_C + \alpha_D + \alpha_H, \alpha_C), $ \\

                          $92$ &

$(-\alpha_D - \alpha_E, -\alpha_B + \alpha_D - \alpha_E, -\alpha_B, \alpha_C + \alpha_E, \alpha_B, \alpha_D, 0, 0, \alpha_C + \alpha_H,$ \\
& \ \ \ \ \ \ \ \ $    \alpha_A + \alpha_B + \alpha_C), $ \\

                          $93$ &

$(0, \alpha_D, -\alpha_B - \alpha_C - \alpha_D - \alpha_E, -\alpha_A - 2 \alpha_B - \alpha_D - \alpha_E,$ \\
& \ \ \ \ \ \ \ \ $    \alpha_A + 2 \alpha_B + \alpha_C + \alpha_D + \alpha_E, \alpha_B + \alpha_D + \alpha_E, 0, \alpha_C, \alpha_H, 0), $ \\

                          $94$ &

$(0, \alpha_D, -\alpha_B - \alpha_C - \alpha_D - \alpha_E, -\alpha_B + \alpha_C - \alpha_D - \alpha_E - \alpha_H, \alpha_B + \alpha_D + \alpha_E,$ \\
& \ \ \ \ \ \ \ \ $    \alpha_B + \alpha_D + \alpha_E + \alpha_H, \alpha_H, \alpha_C, 0, \alpha_A + \alpha_B + \alpha_C), $ \\

                          $95$ &

$(0, -\alpha_B + \alpha_D - \alpha_E - \alpha_H, -\alpha_B - \alpha_D - \alpha_E, \alpha_C - \alpha_D, \alpha_B + \alpha_D + \alpha_E, \alpha_D,$ \\
& \ \ \ \ \ \ \ \ $    \alpha_H, 0, \alpha_C, \alpha_A + \alpha_B + \alpha_C), $ \\

                          $96$ &

$(0, -\alpha_B + \alpha_D - \alpha_E - \alpha_H, -\alpha_B - \alpha_C - \alpha_D - \alpha_E, \alpha_C, \alpha_B + \alpha_E, 0, \alpha_H,$ \\
& \ \ \ \ \ \ \ \ $    \alpha_C + \alpha_D, 0, \alpha_A + \alpha_B + \alpha_C), $ \\

                          $97$ &

$(0, \alpha_D, 0, -\alpha_A - 2 \alpha_B - \alpha_D - \alpha_E - \alpha_H, \alpha_A + \alpha_B, \alpha_B + \alpha_D + \alpha_E + \alpha_H,$ \\
& \ \ \ \ \ \ \ \ $    \alpha_B + \alpha_C + \alpha_D + \alpha_E + \alpha_H, \alpha_C, 0, 0), $ \\

                          $98$ &

$(-\alpha_A - \alpha_D - \alpha_E, \alpha_D - \alpha_E, \alpha_A, \alpha_C + \alpha_E, 0, \alpha_D, 0, 0, \alpha_B + \alpha_C + \alpha_H,$ \\
& \ \ \ \ \ \ \ \ $    \alpha_B + \alpha_C), $ \\

                          $99$ &

$(-\alpha_E, -\alpha_B + \alpha_D - \alpha_E - \alpha_H, -\alpha_B, \alpha_C + \alpha_E, \alpha_B, 0, \alpha_H, 0, \alpha_C + \alpha_D,$ \\
& \ \ \ \ \ \ \ \ $    \alpha_A + \alpha_B + \alpha_C), $ \\

                         $100$ &

$(-\alpha_E, -\alpha_B + \alpha_D - \alpha_E - \alpha_H, 0, \alpha_C + \alpha_E, -\alpha_C - \alpha_D, 0,$ \\
& \ \ \ \ \ \ \ \ $    \alpha_B + \alpha_C + \alpha_D + \alpha_H, \alpha_C + \alpha_D, 0, \alpha_A + \alpha_B + \alpha_C), $ \\

                         $101$ &

$(-\alpha_A - \alpha_E, -\alpha_B + \alpha_D - \alpha_E - \alpha_H, \alpha_A, \alpha_C + \alpha_E, 0, 0, \alpha_B + \alpha_H, 0,$ \\
& \ \ \ \ \ \ \ \ $    \alpha_C + \alpha_D, \alpha_B + \alpha_C), $ \\

                         $102$ &

$(-\alpha_A - \alpha_E, \alpha_D - \alpha_E, \alpha_A, \alpha_C + \alpha_E, 0, -\alpha_B - \alpha_C, 0, \alpha_B + \alpha_C,$ \\
& \ \ \ \ \ \ \ \ $    \alpha_B + \alpha_C + \alpha_D + \alpha_H, 0), $ \\

                         $103$ &

$(0, -\alpha_B + \alpha_D - \alpha_E - \alpha_H, 0, \alpha_C, -\alpha_C - \alpha_D, -\alpha_A - \alpha_B - \alpha_C,$ \\
& \ \ \ \ \ \ \ \ $    \alpha_B + \alpha_C + \alpha_D + \alpha_E + \alpha_H, \alpha_A + \alpha_B + 2 \alpha_C + \alpha_D, 0, 0), $ \\

                         $104$ &

$(0, \alpha_D, -\alpha_B - \alpha_C - \alpha_D - \alpha_E, \alpha_C, -\alpha_H, -\alpha_A - \alpha_B - \alpha_C, \alpha_H,$ \\
& \ \ \ \ \ \ \ \ $    \alpha_A + 2 \alpha_B + 2 \alpha_C + \alpha_D + \alpha_E + \alpha_H, 0, 0), $ \\

                         $105$ &

$(-\alpha_A, \alpha_D, \alpha_A, \alpha_C, -\alpha_B - \alpha_C - \alpha_D - \alpha_E - \alpha_H, -\alpha_B - \alpha_C,$ \\
& \ \ \ \ \ \ \ \ $    \alpha_B + \alpha_C + \alpha_D + \alpha_E + \alpha_H, 2 \alpha_B + 2 \alpha_C + \alpha_D + \alpha_E + \alpha_H, 0, 0),$ \\

                         $106$ &

$(0, -\alpha_B + \alpha_D - \alpha_E - \alpha_H, 0, \alpha_C - \alpha_D, -\alpha_C, \alpha_D,$ \\
& \ \ \ \ \ \ \ \ $    \alpha_B + \alpha_C + \alpha_D + \alpha_E + \alpha_H, \alpha_C, 0, \alpha_A + \alpha_B + \alpha_C), $ \\

                         $107$ &

$(-\alpha_A, \alpha_D, \alpha_A, -\alpha_B + \alpha_C, \alpha_B, -\alpha_C, 0, \alpha_C, \alpha_B + \alpha_C + \alpha_D + \alpha_E + \alpha_H, 0),$ \\

                         $108$ &

$(-\alpha_E, -\alpha_B + \alpha_D - \alpha_E, -\alpha_B, \alpha_C + \alpha_E, \alpha_B, 0, 0, 0, \alpha_C + \alpha_D + \alpha_H,$ \\
& \ \ \ \ \ \ \ \ $    \alpha_A + \alpha_B + \alpha_C), $ \\

                         $109$ &

$(-\alpha_A - \alpha_E, \alpha_D - \alpha_E, \alpha_A, \alpha_C + \alpha_E, 0, 0, 0, 0, \alpha_B + \alpha_C + \alpha_D + \alpha_H,$ \\
& \ \ \ \ \ \ \ \ $    \alpha_B + \alpha_C), $ \\

                         $110$ &

$(-\alpha_D - \alpha_E, \alpha_D - \alpha_E, 0, \alpha_C + \alpha_E, 0, \alpha_D, 0, 0, \alpha_B + \alpha_C + \alpha_H,$ \\
& \ \ \ \ \ \ \ \ $    \alpha_A + \alpha_B + \alpha_C), $ \\

                         $111$ &

$(-\alpha_A - \alpha_D - \alpha_E, \alpha_D, \alpha_A, \alpha_C, 0, \alpha_D + \alpha_E, 0, 0, \alpha_B + \alpha_C + \alpha_H, \alpha_B + \alpha_C
    ), $ \\

                         $112$ &

$(-\alpha_E, \alpha_D - \alpha_E, 0, \alpha_C + \alpha_E, -\alpha_B - \alpha_C - \alpha_D - \alpha_H, 0, \alpha_B + \alpha_C + \alpha_D + \alpha_H,$ \\
& \ \ \ \ \ \ \ \ $    \alpha_B + \alpha_C + \alpha_D + \alpha_H, 0, \alpha_A + \alpha_B + \alpha_C), $ \\

                         $113$ &

$(-\alpha_E, -\alpha_B + \alpha_D - \alpha_E - \alpha_H, 0, \alpha_C + \alpha_E, 0, 0, \alpha_B + \alpha_H, 0, \alpha_C + \alpha_D,$ \\
& \ \ \ \ \ \ \ \ $    \alpha_A + \alpha_B + \alpha_C), $ \\

                         $114$ &

$(-\alpha_E, \alpha_D - \alpha_E, 0, \alpha_C + \alpha_E, 0, -\alpha_A - \alpha_B - \alpha_C, 0, \alpha_A + \alpha_B + \alpha_C,$ \\
& \ \ \ \ \ \ \ \ $    \alpha_B + \alpha_C + \alpha_D + \alpha_H, 0), $ \\

                         $115$ &

$(0, -\alpha_B + \alpha_D - \alpha_E - \alpha_H, 0, \alpha_C, 0, -\alpha_A - \alpha_B - \alpha_C, \alpha_B + \alpha_E + \alpha_H,$ \\
& \ \ \ \ \ \ \ \ $    \alpha_A + \alpha_B + \alpha_C, \alpha_C + \alpha_D, 0), $ \\

                         $116$ &

$(0, \alpha_D, 0, \alpha_C, -\alpha_B - \alpha_C - \alpha_D - \alpha_E - \alpha_H, -\alpha_A - \alpha_B - \alpha_C,$ \\
& \ \ \ \ \ \ \ \ $    \alpha_B + \alpha_C + \alpha_D + \alpha_E + \alpha_H, \alpha_A + 2 \alpha_B + 2 \alpha_C + \alpha_D + \alpha_E + \alpha_H, 0, 0)
    , $ \\

                         $117$ &

$(0, -\alpha_B + \alpha_D - \alpha_E - \alpha_H, -\alpha_B - \alpha_E, \alpha_C, \alpha_B + \alpha_E, 0, \alpha_H, 0, \alpha_C + \alpha_D,$ \\
& \ \ \ \ \ \ \ \ $    \alpha_A + \alpha_B + \alpha_C), $ \\

                         $118$ &

$(-\alpha_A, -\alpha_B + \alpha_D - \alpha_E - \alpha_H, \alpha_A, \alpha_C, 0, 0, \alpha_B + \alpha_E + \alpha_H, 0, \alpha_C + \alpha_D,$ \\
& \ \ \ \ \ \ \ \ $    \alpha_B + \alpha_C), $ \\

                         $119$ &

$(0, -\alpha_B + \alpha_D - \alpha_E - \alpha_H, 0, \alpha_C, -\alpha_C - \alpha_D, 0, \alpha_B + \alpha_C + \alpha_D + \alpha_E + \alpha_H,$ \\
& \ \ \ \ \ \ \ \ $    \alpha_C + \alpha_D, 0, \alpha_A + \alpha_B + \alpha_C), $ \\

                         $120$ &

$(0, -\alpha_B + \alpha_D - \alpha_E - \alpha_H, 0, \alpha_C - \alpha_D, 0, \alpha_D, \alpha_B + \alpha_D + \alpha_E + \alpha_H, 0,$ \\
& \ \ \ \ \ \ \ \ $    \alpha_C, \alpha_A + \alpha_B + \alpha_C), $ \\

                         $121$ &

$(0, \alpha_D, 0, -\alpha_B + \alpha_C - \alpha_D - \alpha_E - \alpha_H, -\alpha_C, \alpha_B + \alpha_D + \alpha_E + \alpha_H,$ \\
& \ \ \ \ \ \ \ \ $    \alpha_B + \alpha_C + \alpha_D + \alpha_E + \alpha_H, \alpha_C, 0, \alpha_A + \alpha_B + \alpha_C), $ \\

                         $122$ &

$(0, \alpha_D, -\alpha_B - \alpha_C - \alpha_D - \alpha_E, \alpha_C, -\alpha_H, 0, \alpha_H, \alpha_B + \alpha_C + \alpha_D + \alpha_E + \alpha_H,$ \\
& \ \ \ \ \ \ \ \ $    0, \alpha_A + \alpha_B + \alpha_C), $ \\

                         $123$ &

$(-\alpha_A, \alpha_D, \alpha_A, \alpha_C, -\alpha_B - \alpha_C - \alpha_D - \alpha_E - \alpha_H, 0, \alpha_B + \alpha_C + \alpha_D + \alpha_E + \alpha_H,$ \\
& \ \ \ \ \ \ \ \ $    \alpha_B + \alpha_C + \alpha_D + \alpha_E + \alpha_H, 0, \alpha_B + \alpha_C), $ \\

                         $124$ &

$(0, \alpha_D, -\alpha_B - \alpha_C - \alpha_D - \alpha_E, \alpha_C, 0, -\alpha_A - \alpha_B - \alpha_C, 0,$ \\
& \ \ \ \ \ \ \ \ $    \alpha_A + 2 \alpha_B + 2 \alpha_C + \alpha_D + \alpha_E, \alpha_H, 0), $ \\

                         $125$ &

$(0, \alpha_D, 0, -\alpha_A - \alpha_B + \alpha_C, \alpha_A + \alpha_B, -\alpha_C, 0, \alpha_C,$ \\
& \ \ \ \ \ \ \ \ $    \alpha_B + \alpha_C + \alpha_D + \alpha_E + \alpha_H, 0), $ \\

                         $126$ &

$(-\alpha_A, \alpha_D, \alpha_A, \alpha_C, 0, -\alpha_B - \alpha_C, 0, \alpha_B + \alpha_C, \alpha_B + \alpha_C + \alpha_D + \alpha_E + \alpha_H, 0
    ), $ \\

                         $127$ &

$(-\alpha_A, \alpha_D, \alpha_A, -\alpha_B + \alpha_C, \alpha_B, 0, 0, 0, \alpha_B + \alpha_C + \alpha_D + \alpha_E + \alpha_H, \alpha_C),$ \\

                         $128$ &

$(-\alpha_D - \alpha_E, \alpha_D, 0, \alpha_C, 0, \alpha_D + \alpha_E, 0, 0, \alpha_B + \alpha_C + \alpha_H, \alpha_A + \alpha_B + \alpha_C
    ), $ \\

                         $129$ &

$(0, \alpha_D, 0, -\alpha_B + \alpha_C - \alpha_D - \alpha_E - \alpha_H, 0, \alpha_B + \alpha_D + \alpha_E + \alpha_H,$ \\
& \ \ \ \ \ \ \ \ $    \alpha_B + \alpha_D + \alpha_E + \alpha_H, 0, \alpha_C, \alpha_A + \alpha_B + \alpha_C), $ \\

                         $130$ &

$(0, \alpha_D, 0, -\alpha_A - \alpha_B + \alpha_C, \alpha_A + \alpha_B, 0, 0, 0, \alpha_B + \alpha_C + \alpha_D + \alpha_E + \alpha_H,$ \\
& \ \ \ \ \ \ \ \ $    \alpha_C), $ \\

                         $131$ &

$(0, -\alpha_B + \alpha_D - \alpha_E, -\alpha_B - \alpha_E, \alpha_C, \alpha_B + \alpha_E, 0, 0, 0, \alpha_C + \alpha_D + \alpha_H,$ \\
& \ \ \ \ \ \ \ \ $    \alpha_A + \alpha_B + \alpha_C), $ \\

                         $132$ &

$(0, \alpha_D, -\alpha_B - \alpha_C - \alpha_D - \alpha_E, \alpha_C, 0, 0, 0, \alpha_B + \alpha_C + \alpha_D + \alpha_E, \alpha_H,$ \\
& \ \ \ \ \ \ \ \ $    \alpha_A + \alpha_B + \alpha_C), $ \\

                         $133$ &

$(-\alpha_A - \alpha_B - \alpha_C - \alpha_D - \alpha_E, \alpha_D, \alpha_A, -\alpha_B + \alpha_C - \alpha_H, \alpha_B,$ \\
& \ \ \ \ \ \ \ \ $    \alpha_B + \alpha_D + \alpha_E + \alpha_H, \alpha_H, \alpha_C, 0, 0), $ \\

                         $134$ &

$(-\alpha_A - \alpha_B - \alpha_C - \alpha_D - \alpha_E, -\alpha_B + \alpha_D - \alpha_E, \alpha_A, \alpha_C + \alpha_E, \alpha_B, \alpha_D, 0, \alpha_C,$ \\
& \ \ \ \ \ \ \ \ $    \alpha_H, 0), $ \\

                         $135$ &

$(-\alpha_A - \alpha_B - \alpha_D - \alpha_E, -\alpha_B + \alpha_D - \alpha_E - \alpha_H, \alpha_A, \alpha_C + \alpha_E, \alpha_B, \alpha_D, \alpha_H, 0,$ \\
& \ \ \ \ \ \ \ \ $    \alpha_C, \alpha_C), $ \\

                         $136$ &

$(-\alpha_D - \alpha_E, -\alpha_B + \alpha_D - \alpha_E - \alpha_H, -\alpha_B - \alpha_C, \alpha_C + \alpha_E, \alpha_B, \alpha_D, \alpha_H, \alpha_C, 0,$ \\
& \ \ \ \ \ \ \ \ $    \alpha_A + \alpha_B + \alpha_C), $ \\

                         $137$ &

$(0, -\alpha_B + \alpha_D - \alpha_E - \alpha_H, -\alpha_B - \alpha_C - \alpha_D - \alpha_E, -\alpha_A - \alpha_B - \alpha_D,$ \\
& \ \ \ \ \ \ \ \ $    \alpha_A + 2 \alpha_B + \alpha_C + \alpha_D + \alpha_E, \alpha_D, \alpha_H, \alpha_C, 0, 0), $ \\

                         $138$ &

$(-\alpha_E, \alpha_D - \alpha_E, 0, \alpha_C + \alpha_E, 0, 0, 0, 0, \alpha_B + \alpha_C + \alpha_D + \alpha_H,$ \\
& \ \ \ \ \ \ \ \ $    \alpha_A + \alpha_B + \alpha_C), $ \\

                         $139$ &

$(0, -\alpha_B + \alpha_D - \alpha_E - \alpha_H, 0, \alpha_C, 0, 0, \alpha_B + \alpha_E + \alpha_H, 0, \alpha_C + \alpha_D,$ \\
& \ \ \ \ \ \ \ \ $    \alpha_A + \alpha_B + \alpha_C), $ \\

                         $140$ &

$(0, \alpha_D, 0, \alpha_C, -\alpha_B - \alpha_C - \alpha_D - \alpha_E - \alpha_H, 0, \alpha_B + \alpha_C + \alpha_D + \alpha_E + \alpha_H,$ \\
& \ \ \ \ \ \ \ \ $    \alpha_B + \alpha_C + \alpha_D + \alpha_E + \alpha_H, 0, \alpha_A + \alpha_B + \alpha_C), $ \\

                         $141$ &

$(0, \alpha_D, 0, \alpha_C, 0, -\alpha_A - \alpha_B - \alpha_C, 0, \alpha_A + \alpha_B + \alpha_C,$ \\
& \ \ \ \ \ \ \ \ $    \alpha_B + \alpha_C + \alpha_D + \alpha_E + \alpha_H, 0), $ \\

                         $142$ &

$(-\alpha_A, \alpha_D, \alpha_A, \alpha_C, 0, 0, 0, 0, \alpha_B + \alpha_C + \alpha_D + \alpha_E + \alpha_H, \alpha_B + \alpha_C),$ \\

                         $143$ &

$(-\alpha_A - \alpha_B - \alpha_C - \alpha_D - \alpha_E, -\alpha_B + \alpha_D - \alpha_E - \alpha_H, \alpha_A, \alpha_C + \alpha_E, \alpha_B, \alpha_D, \alpha_H,$ \\
& \ \ \ \ \ \ \ \ $    \alpha_C, 0, 0), $ \\

                         $144$ &

$(0, \alpha_D, 0, \alpha_C, 0, 0, 0, 0, \alpha_B + \alpha_C + \alpha_D + \alpha_E + \alpha_H, \alpha_A + \alpha_B + \alpha_C).$ \\
\end{longtable}
\end{center}

We note that not all the regions of linearity are defined by the same number of
inequalities. There are $62$ regions ($1$ to $62$) defined by
$6$ inequalities, $70$ regions ($63$ to $132$) defined by $7$ inequalities,
$10$ regions defined by $8$ inequalities ($133$ to $142$) and $2$ regions
($143$ to $144$) defined by $11$ inequalities.

We note in particular that the number ($62$) of regions defined by the minimum
number of inequalities is equal to the number of equivalence classes of
reduced expressions for $w_0$. We shall show that there is in fact a natural
bijection between these two sets.

Now the vectors which parametrize the PBW-basis elements of $U^-$ have all
coordinates in $\N$. Thus we shall add the inequalities
$$a\geq 0,\ b\geq 0,\ c\geq 0,\ d\geq 0,\ e\geq 0,\ f\geq 0,\ g\geq 0,\ h\geq 0,\ i\geq 0,\ j\geq 0$$
to those of Table III.

We define a {\em simplicial region} to be the subset of $\R^{10}$ defined by one of the
sets of $6$ inequalities in Table III together with the above inequalities asserting
that all coordinates are nonnegative. Of these 10 latter inequalities there will be
$4$ from which the remaining $6$ follow using the $6$ inequalities from Table III.
Thus each simplicial region will have $10$ walls, $6$ of which are given by Table III and
the remaining $4$ will have form that some coordinate is greater than or equal to
zero. The $62$ simplicial regions together with their walls are shown in Table V
(where, for example, $abgh$ indicates the four walls $a\geq 0$, $b\geq 0$,
$g\geq 0$ and $h\geq 0$).

\begin{center} Table V \\
{\em The simplicial regions and their walls.}

\begin{longtable}{clc}
Region & \ \ \ \ \ Main walls & Coordinate walls
\endfirsthead
Region & Main walls & Coordinate walls
\endhead
\endfoot
\endlastfoot
1   &   [A, BC, CD, E, BH]   $\|$   [BCD]   &   abgh   \\
2   &   [A, BC, E, BH]   $\|$   [B, CD]   &   abgh   \\
3   &   [ABC, CD, E, BH]   $\|$   [A, BCD]   &   begh   \\
4   &   [A, BC, E, BCDH]   $\|$   [BCD, BH]   &   abgh   \\
5   &   [A, CD, E, BH]   $\|$   [BC, D]   &   abgh   \\
6   &   [A, BC, CD, BEH]   $\|$   [BCD, E]   &   afgh   \\
7   &   [A, BC, BEH]   $\|$   [B, CD, E]   &   afgh   \\
8   &   [ABC, E, BH]   $\|$   [A, B, CD]   &   begh   \\
9   &   [ABC, E, BCDH]   $\|$   [A, BCD, BH]   &   begh   \\
10   &   [A, E, BCDH]   $\|$   [BC, D, BH]   &   abgh   \\
11   &   [A, CD, BEH]   $\|$   [BC, D, E]   &   afgh   \\
12   &   [A, B, C, BEH]   $\|$   [CD, BE]   &   acfh   \\
13   &   [B, ABC, E, H]   $\|$   [AB, CD]   &   bceh   \\
14   &   [ABC, BCD, E, H]   $\|$   [ABCD, BH]   &   begh   \\
15   &   [A, D, E, BCH]   $\|$   [BC, BH]   &   abdg   \\
16   &   [A, C, D, BEH]   $\|$   [BC, DE]   &   adfg   \\
17   &   [A, B, BEH]   $\|$   [C, D, BE]   &   achj   \\
18   &   [AB, C, BEH]   $\|$   [A, CD, BE]   &   cefh   \\
19   &   [ABC, BE, H]   $\|$   [AB, CD, E]   &   cefh   \\
20   &   [B, ABC, E]   $\|$   [AB, CD, H]   &   bchi   \\
21   &   [ABC, BCD, E]   $\|$   [B, ABCD, H]   &   bghi   \\
22   &   [BCD, E, H]   $\|$   [ABC, D, BH]   &   begh   \\
23   &   [D, E, BCH]   $\|$   [A, BC, BH]   &   bdeg   \\
24   &   [A, DE, BCH]   $\|$   [BC, E, BH]   &   adfg   \\
25   &   [A, C, BDEH]   $\|$   [BC, DE, BEH]   &   adfg   \\
26   &   [A, D, BEH]   $\|$   [B, C, DE]   &   adgj   \\
27   &   [A, B, D, BEH]   $\|$   [C, BDE]   &   acdj   \\
28   &   [AB, C, BE, H]   $\|$   [CD, ABE]   &   cefh   \\
29   &   [B, ABC, CD, E]   $\|$   [ABCD, H]   &   bchi   \\
30   &   [BC, D, E, H]   $\|$   [ABC, BH]   &   bdeg   \\
31   &   [A, C, DE, BH]   $\|$   [BC, BEH]   &   adfg   \\
32   &   [A, BDEH]   $\|$   [B, C, DE, BEH]   &   adgj   \\
33   &   [AB, BEH]   $\|$   [A, C, D, BE]   &   cehj   \\
34   &   [ABC, BE]   $\|$   [AB, CD, E, H]   &   cfhi   \\
35   &   [BCD, E]   $\|$   [B, ABC, D, H]   &   bghi   \\
36   &   [DE, BCH]   $\|$   [A, BC, E, BH]   &   defg   \\
37   &   [C, DE, BH]   $\|$   [A, BC, BEH]   &   defg   \\
38   &   [A, DE, BH]   $\|$   [B, C, BEH]   &   adgj   \\
39   &   [A, B, BDEH]   $\|$   [C, BDE, BEH]   &   acdj   \\
40   &   [AB, D, BEH]   $\|$   [A, C, BDE]   &   cdej   \\
41   &   [AB, BE, H]   $\|$   [C, D, ABE]   &   cehj   \\
42   &   [AB, C, BE]   $\|$   [CD, ABE, H]   &   cfhi   \\
43   &   [ABC, CD, BE]   $\|$   [ABCD, E, H]   &   cfhi   \\
44   &   [B, CD, E]   $\|$   [ABC, D, H]   &   bchi   \\
45   &   [BC, D, E]   $\|$   [B, ABC, H]   &   bdgi   \\
46   &   [BC, DE, H]   $\|$   [ABC, E, BH]   &   defg   \\
47   &   [DE, BH]   $\|$   [A, B, C, BEH]   &   degj   \\
48   &   [AB, BDEH]   $\|$   [A, C, BDE, BEH]   &   cdej   \\
49   &   [AB, BE]   $\|$   [C, D, ABE, H]   &   chij   \\
50   &   [CD, BE]   $\|$   [ABC, D, E, H]   &   cfhi   \\
51   &   [BC, DE]   $\|$   [B, ABC, E, H]   &   dfgi   \\
52   &   [B, DE, H]   $\|$   [AB, C, BEH]   &   cdej   \\
53   &   [AB, BDE, H]   $\|$   [C, ABDE, BEH]   &   cdej   \\
54   &   [AB, D, BE]   $\|$   [C, ABDE, H]   &   cdij   \\
55   &   [C, D, BE]   $\|$   [ABC, DE, H]   &   cdfi   \\
56   &   [B, C, DE]   $\|$   [ABC, BE, H]   &   cdfi   \\
57   &   [B, DE]   $\|$   [AB, C, BE, H]   &   cdij   \\
58   &   [BDE, H]   $\|$   [AB, C, DE, BEH]   &   cdej   \\
59   &   [AB, BDE]   $\|$   [C, BE, ABDE, H]   &   cdij   \\
60   &   [D, BE]   $\|$   [AB, C, DE, H]   &   cdij   \\
61   &   [C, BDE]   $\|$   [ABC, BE, DE, H]   &   cdfi   \\
62   &   [BDE]   $\|$   [AB, C, BE, DE, H]   &   cdij   \\
\end{longtable}
\end{center}

In order to define a natural bijection between the set of equivalence relations of
reduced words for $w_0$ and the set of simplicial regions we consider the transforms
under $S_{\mathbf{i}}^{\mathbf{j}}$ of the spanning vectors of $C_{\mathbf{i}}$
considered in $\S4$. We have one such vector $v_P$ for each partial quiver $P$ and
one vector $v_j$ for each $j=1,2,3,4$. These vectors $v_P,v_j$ can, for example,
be calculated by the rectangle algorithm described in $\S4$. These vectors are
listed in Table VI.

We consider the incidence properties relating these vectors $v_P$ and $v_j$ to the
walls of the simplicial regions. There are altogether $32$ hyperplanes in $\R^{10}$ which
arise as walls of simplicial regions. $22$ of these are hyperplanes appearing in Table III
and $10$ are the hyperplanes obtained by putting one coordinate equal to $0$. There
are $22$ vectors of form $v_P$ and $4$ vectors of form $v_j$. These are obtained from
the rectangle algorithm which associates with each partial quiver $P$ or each
$j\in \{1,2,3,4\}$ a set of positive roots. Such a set of positive roots is
translated into a vector with components $0$ or $1$ by means of the ordering
$$(\alpha_1,\alpha_3,\alpha_1+\alpha_2+\alpha_3,\alpha_3+\alpha_4,\alpha_2+\alpha_3,
\alpha_1+\alpha_2+\alpha_3+\alpha_4,\alpha_2+\alpha_3+\alpha_4,\alpha_1+\alpha_2,
\alpha_4,\alpha_2)$$ determined by the reduced word $\mathbf{j}=
(1,3,2,4,1,3,2,4,1,3)$. In Table VII
we indicate which vectors $v_P$ and $v_j$
lie on which walls of simplicial regions. A cross in a particular position indicates that
the given vector does not lie on the given wall.

\begin{center} Table VI \\ {\em Vectors obtained by the rectangle algorithm.}

\renewcommand{\arraystretch}{1.8}
\begin{longtable}{cc}
Partial quiver $P$ & \ \ \ \ \ \ \ \ Vector $v_P$ \\
\endfirsthead
Partial quiver $P$ & \ \ \ \ \ \ \ \ Vector $v_P$ \\
\endhead
\endfoot
\endlastfoot
LLL & (1,0,0,0,0,0,0,0,0,0) \\
LLR & (0,0,1,0,0,0,0,0,1,0) \\
LRL & (1,1,0,0,0,0,1,1,0,0) \\
RLL & (1,0,0,1,0,0,0,0,0,1) \\
LRR & (0,1,0,0,0,0,0,1,1,0) \\
RLR & (0,0,1,1,0,0,0,0,1,1) \\
RRL & (1,0,0,0,0,0,1,0,0,0) \\
RRR & (0,0,0,0,0,0,0,0,1,0) \\
LL- & (0,0,0,0,0,0,0,1,0,0) \\
LR- & (0,1,0,0,0,1,0,0,0,0) \\
RL- & (0,0,0,0,0,1,0,0,0,1) \\
RR- & (0,0,0,1,0,0,0,0,0,0) \\
\0LL & (1,0,0,0,0,0,0,0,0,1) \\
\0LR & (0,0,1,0,0,0,0,0,1,1) \\
\0RL & (1,1,0,0,0,0,1,0,0,0) \\
\0RR & (0,1,0,0,0,0,0,0,1,0) \\
L\0- & (0,1,0,0,0,0,0,1,0,0) \\
R\0- & (0,0,0,1,0,0,0,0,0,1) \\
\0L- & (0,0,1,0,0,0,0,0,0,1) \\
\0R- & (0,1,0,0,0,0,1,0,0,0) \\
\0\0L & (1,0,0,0,1,0,0,0,0,0) \\
\0\0R & (0,0,0,0,1,0,0,0,1,0) \\
 & \\
Integer $j$ & \ \ \ \ \ \ \ \ Vector $v_j$ \\
 & \\
$1$ & (1,0,0,0,1,0,0,0,1,0) \\
$2$ & (0,0,1,0,0,0,1,0,0,0) \\
$3$ & (0,1,0,0,0,1,0,0,0,1) \\
$4$ & (0,0,0,1,0,0,0,1,0,0)
\end{longtable}
\end{center}

\renewcommand{\arraystretch}{1.2}
\begin{center} Table VII \\ {\em Incidence table: spanning vectors and walls
of simplicial regions.} \\ x indicates that the vector lies on the given
wall. \\ \ \\
\begin{tabular}{c|@{\hspace{6pt}}c@{\hspace{6pt}}c@{\hspace{6pt}}c@{\hspace{6pt}}c@{\hspace{6pt}}c@{\hspace{6pt}}c@{\hspace{6pt}}c@{\hspace{6pt}}c@{\hspace{6pt}}c@{\hspace{6pt}}c@{\hspace{6pt}}c@{\hspace{6pt}}c@{\hspace{6pt}}c@{\hspace{6pt}}c@{\hspace{6pt}}c@{\hspace{6pt}}c@{\hspace{6pt}}c@{\hspace{6pt}}c@{\hspace{6pt}}c@{\hspace{6pt}}c@{\hspace{6pt}}c@{\hspace{6pt}}c|@{\hspace{6pt}}c@{\hspace{6pt}}c@{\hspace{6pt}}c@{\hspace{6pt}}c}
 & L & L & L & R & L & R & R & R & L & L & R & R & - & - & - & - & L & R & - & - & - & - &   &   &   &   \\
 & L & L & R & L & R & L & R & R & L & R & L & R & L & L & R & R & - & - & L & R & - & - & 1 & 2 & 3 & 4 \\
 & L & R & L & L & R & R & L & R & - & - & - & - & L & R & L & R & - & - & - & - & L & R
 &   &   &   &   \\
\hline
A & . & x & . & . & x & x & . & x & x & x & x & x & . & x & . & x & x & x & x & x & x & . & x & x & x & x \\

B & x & . & . & x & x & . & . & x & x & x & x & x & x & . & . & x & x & x & . & . & x & x & x & x & x & x \\

C & x & x & x & . & x & . & x & x & x & . & x & x & . & . & x & x & x & . & . & x & x & x & x & x & x & x \\

D & x & x & . & . & . & . & x & x & . & x & x & . & x & x & x & x & . & . & x & x & x & x & x & x & x & x \\

E & x & x & . & x & . & x & x & x & x & x & . & x & x & x & . & . & . & x & x & . & x & x & x & x & x & x \\

H & x & . & x & x & . & . & x & . & x & x & x & x & x & . & x & . & x & x & x & x & . & x & x & x & x & x \\

AB & . & . & x & . & x & . & x & x & x & x & x & x & . & . & x & x & x & x & . & . & x & . & x & x & x & x \\

BC & x & . & . & . & x & x & . & x & x & . & x & x & . & x & . & x & x & . & x & . & x & x & x & x & x & x \\

CD & x & x & . & x & . & x & x & x & . & . & x & . & . & . & x & x & . & x & . & x & x & x & x & x & x & x \\

DE & x & x & x & . & x & . & x & x & . & x & . & . & x & x & . & . & x & . & x & . & x & x & x & x & x & x \\

BE & x & . & x & x & . & . & . & x & x & x & . & x & x & . & x & . & . & x & . & x & x & x & x & x & x & x \\

BH & x & x & . & x & . & x & . & . & x & x & x & x & x & x & . & . & x & x & . & . & . & x & x & x & x & x \\

ABC & . & . & x & x & x & x & x & x & x & . & x & x & x & x & x & x & x & . & x & . & x & . & x & x & x & x \\

BCD & x & . & . & x & . & . & . & x & . & . & x & . & . & x & . & x & . & x & x & . & x & x & x & x & x & x \\

ABE & . & . & . & . & . & . & x & x & x & x & . & x & . & . & . & . & . & x & . & x & x & . & x & x & x & x \\

BDE & x & . & . & . & x & . & . & x & . & x & . & . & x & . & x & . & x & . & . & x & x & x & x & x & x & x \\

BCH & x & x & . & . & . & . & . & . & x & . & x & x & . & . & . & . & x & . & x & . & . & x & x & x & x & x \\

BEH & x & x & x & x & x & x & . & . & x & x & . & x & x & x & x & x & . & x & . & x & . & x & x & x & x & x \\

ABCD & . & . & . & . & . & . & x & x & . & . & x & . & x & x & x & x & . & x & x & . & x & . & x & x & x & x \\

ABDE & . & . & x & . & x & . & x & x & . & x & . & . & . & . & . & . & x & . & . & x & x & . & x & x & x & x \\

BCDH & x & x & . & x & . & x & . & . & . & . & x & . & . & . & . & . & . & x & x & . & . & x & x & x & x & x \\

BDEH & x & x & . & . & . & . & . & . & . & x & . & . & x & x & x & x & x & . & . & x & . & x & x & x & x & x \\
\hline
a & . & x & . & . & x & x & . & x & x & x & x & x & . & x & . & x & x & x & x & x & . & x & . & x & x & x \\

b & x & x & . & x & . & x & x & x & x & . & x & x & x & x & . & . & . & x & x & . & x & x & x & x & . & x \\

c & x & . & x & x & x & . & x & x & x & x & x & x & x & . & x & x & x & x & . & x & x & x & x & . & x & x \\

d & x & x & x & . & x & . & x & x & x & x & x & . & x & x & x & x & x & . & x & x & x & x & x & x & x & . \\

e & x & x & x & x & x & x & x & x & x & x & x & x & x & x & x & x & x & x & x & x & . & . & . & x & x & x \\

f & x & x & x & x & x & x & x & x & x & . & . & x & x & x & x & x & x & x & x & x & x & x & x & x & . & x \\

g & x & x & . & x & x & x & . & x & x & x & x & x & x & x & . & x & x & x & x & . & x & x & x & . & x & x \\

h & x & x & . & x & . & x & x & x & . & x & x & x & x & x & x & x & . & x & x & x & x & x & x & x & x & . \\

i & x & . & x & x & . & . & x & . & x & x & x & x & x & . & x & . & x & x & x & x & x & . & . & x & x & x \\

j & x & x & x & . & x & . & x & x & x & x & . & x & . & . & x & x & x & . & . & x & x & x & x & x & . & x
\end{tabular}
\end{center}
\pagebreak
\renewcommand{\arraystretch}{1.4}

The incidence table enables us to describe a bijection between equivalence classes
of reduced words and simplicial regions. Let $\mathcal{W}$ be the set of hyperplanes in
$\R^{10}$ which arise as walls of simplicial regions. We have seen that $|\mathcal{W}|=22+10=
32$.

\begin{prop} \label{iso}
Let $\mathbf{i}$ be a reduced expression for $w_0$. Let $\mathcal{P}(\mathbf{i})$ be
the set of $6$ partial quivers obtained from $\mathbf{i}$ using the chamber sets as
in $\S3$. Consider the $10$ vectors $\{v_P,\ P\in \mathcal{P}(\mathbf{i});\ v_j,\ 
j=1,2,3,4\}$. For each vector in the set consider the set of $9$ remaining vectors.
Then there is a unique hyperplane $W\in\mathcal{W}$ such that these $9$ vectors lie
on $W$. Moreover, the set of $10$ elements of $\mathcal{W}$ obtained in this way
are the $10$ boundary walls of a simplicial region $X_{\mathbf{i}}$. The original $10$
vectors $v_P,\ P\in \mathcal{P}(\mathbf{i});\ v_j,\ j=1,2,3,4$ all lie in
$X_{\mathbf{i}}$. Also, the map $\mathbf{i}\rightarrow X_{\mathbf{i}}$ gives a
bijection between equivalence classes of reduced words and simplicial regions.
\end{prop}

The bijection $\mathbf{i}\rightarrow X_{\mathbf{i}}$ is described in Table VIII.

\begin{cor} \label{spvects}
The vectors $v_P,\ P\in \mathcal{P}(\mathbf{i});\ v_j,\ j=1,2,3,4$ are spanning
vectors of the region $X_{\mathbf{i}}$.
\end{cor}

\noindent {\bf Proof:} This follows from the fact that each vector lies in $X_{\mathbf{i}}$ and
that each wall of $X_{\mathbf{i}}$ contains $9$ of the $10$ vectors.~$\Box$

\renewcommand{\arraystretch}{2}
\begin{center} Table VIII \\ {\em The correspondence between reduced words and simplicial regions.}
\\ \ \\
\begin{longtable}{cl}
Simplicial region & \ \ \ \ \ \ \ \ Reduced word \\ &
\endfirsthead
Simplicial region & \ \ \ \ \ \ \ \ Reduced word \\ &
\endhead
\endfoot
\endlastfoot
                  1  &  [1, 3, 2, 1, 3, 4, 3, 2, 1, 3]  \\

                  2  &  [1, 3, 2, 1, 3, 4, 2, 3, 2, 1]  \\

                  3  &  [3, 2, 1, 2, 3, 4, 3, 2, 1, 3]  \\

                  4  &  [1, 3, 2, 4, 1, 3, 2, 4, 1, 3]  \\

                  5  &  [1, 3, 2, 3, 4, 3, 2, 1, 2, 3]  \\

                  6  &  [1, 2, 3, 2, 1, 4, 3, 2, 1, 3]  \\

                  7  &  [1, 2, 3, 2, 1, 4, 2, 3, 2, 1]  \\

                  8  &  [3, 2, 1, 2, 3, 4, 2, 3, 2, 1]  \\

                  9  &  [3, 2, 1, 2, 4, 3, 2, 4, 1, 3]  \\

                  10  &  [1, 3, 2, 4, 3, 2, 1, 2, 4, 3]  \\

                  11  &  [1, 2, 3, 2, 4, 3, 2, 1, 2, 3]  \\

                  12  &  [1, 2, 3, 1, 2, 1, 4, 3, 2, 1]  \\

                  13  &  [3, 2, 1, 3, 2, 3, 4, 3, 2, 1]  \\

                  14  &  [3, 2, 1, 4, 3, 2, 3, 4, 1, 3]  \\

                  15  &  [1, 3, 4, 3, 2, 3, 1, 2, 4, 3]  \\

                  16  &  [1, 2, 3, 4, 3, 2, 3, 1, 2, 3]  \\

                  17  &  [1, 2, 3, 1, 2, 4, 3, 2, 1, 2]  \\

                  18  &  [2, 1, 2, 3, 2, 1, 4, 3, 2, 1]  \\

                  19  &  [2, 3, 2, 1, 2, 3, 4, 3, 2, 1]  \\

                  20  &  [3, 2, 1, 3, 4, 2, 3, 2, 4, 1]  \\

                  21  &  [3, 2, 1, 4, 3, 2, 4, 3, 4, 1]  \\

                  22  &  [3, 2, 4, 3, 2, 1, 2, 3, 4, 3]  \\

                  23  &  [3, 4, 3, 2, 1, 2, 3, 2, 4, 3]  \\

                  24  &  [1, 4, 3, 4, 2, 3, 1, 2, 4, 3]  \\

                  25  &  [1, 2, 4, 3, 4, 2, 3, 1, 2, 3]  \\

                  26  &  [1, 2, 3, 4, 3, 2, 1, 2, 3, 2]  \\

                  27  &  [1, 2, 3, 4, 3, 1, 2, 1, 3, 2]  \\

                  28  &  [2, 3, 1, 2, 1, 3, 4, 3, 2, 1]  \\

                  29  &  [3, 2, 1, 3, 4, 3, 2, 3, 4, 1]  \\

                  30  &  [3, 4, 3, 2, 3, 1, 2, 3, 4, 3]  \\

                  31  &  [1, 4, 3, 2, 3, 4, 3, 1, 2, 3]  \\

                  32  &  [1, 2, 4, 3, 4, 2, 1, 2, 3, 2]  \\

                  33  &  [2, 1, 2, 3, 2, 4, 3, 2, 1, 2]  \\

                  34  &  [2, 3, 2, 1, 4, 2, 3, 2, 4, 1]  \\

                  35  &  [3, 2, 4, 3, 2, 1, 2, 4, 3, 4]  \\

                  36  &  [4, 3, 4, 2, 1, 2, 3, 2, 4, 3]  \\

                  37  &  [4, 3, 2, 1, 2, 3, 4, 3, 2, 3]  \\

                  38  &  [1, 4, 3, 2, 3, 4, 1, 2, 3, 2]  \\

                  39  &  [1, 2, 4, 3, 4, 1, 2, 1, 3, 2]  \\

                  40  &  [2, 1, 2, 3, 4, 3, 2, 1, 3, 2]  \\

                  41  &  [2, 3, 1, 2, 3, 4, 3, 2, 1, 2]  \\

                  42  &  [2, 3, 1, 2, 1, 4, 3, 2, 4, 1]  \\

                  43  &  [2, 3, 2, 1, 4, 3, 2, 3, 4, 1]  \\

                  44  &  [3, 2, 3, 4, 3, 2, 1, 2, 3, 4]  \\

                  45  &  [3, 4, 3, 2, 3, 1, 2, 4, 3, 4]  \\

                  46  &  [4, 3, 4, 2, 3, 1, 2, 3, 4, 3]  \\

                  47  &  [4, 3, 2, 1, 2, 3, 4, 2, 3, 2]  \\

                  48  &  [2, 1, 2, 4, 3, 4, 2, 1, 3, 2]  \\

                  49  &  [2, 3, 1, 2, 4, 3, 2, 1, 2, 4]  \\

                  50  &  [2, 3, 2, 4, 3, 2, 1, 2, 3, 4]  \\

                  51  &  [4, 3, 4, 2, 3, 1, 2, 4, 3, 4]  \\

                  52  &  [4, 3, 2, 1, 3, 2, 3, 4, 3, 2]  \\

                  53  &  [2, 1, 4, 3, 2, 3, 4, 1, 3, 2]  \\

                  54  &  [2, 3, 4, 3, 1, 2, 1, 3, 2, 4]  \\

                  55  &  [2, 3, 4, 3, 2, 3, 1, 2, 3, 4]  \\

                  56  &  [4, 3, 2, 3, 4, 3, 1, 2, 3, 4]  \\

                  57  &  [4, 3, 2, 3, 4, 1, 2, 3, 2, 4]  \\

                  58  &  [4, 2, 3, 2, 1, 2, 3, 4, 3, 2]  \\

                  59  &  [4, 2, 3, 4, 1, 2, 1, 3, 2, 4]  \\

                  60  &  [2, 3, 4, 3, 2, 1, 2, 3, 2, 4]  \\

                  61  &  [4, 2, 3, 2, 4, 3, 1, 2, 3, 4]  \\

                  62  &  [4, 2, 3, 2, 4, 1, 2, 3, 2, 4]
\end{longtable}
\end{center}

\section{The transition function $S_{\mathbf{i}}^{\mathbf{j}}$ on the
Lusztig cone $C_{\mathbf{i}}$}

We shall now describe an algorithm for calculating
$S_{\mathbf{i}}^{\mathbf{j}}(\mathbf{a})$ for $\mathbf{a}\in C_{\mathbf{i}}$. As
before $\mathbf{j}$ is the reduced word $(1,3,2,4,1,3,2,4,1,3)$.
The algorithm is best indicated by an example. We illustrate in the case when
$\mathbf{i}=(3,2,1,4,3,2,3,4,1,3)$. The Lusztig cone for this reduced word
$\mathbf{i}$ is given by
$$C_{\mathbf{i}}=\{\mathbf{a}\in \Z^{10}\,:\,P_{\mathbf{i}}\mathbf{a}\geq 0\},$$
where
\renewcommand{\arraystretch}{1}
$$P_{\mathbf{i}}=\left( \begin{array}{cccccccccc}
-1 &   1 &   0 &   1 &  -1 &   0 &   0 &   0 &   0 &   0 \\
 0 &  -1 &   1 &   0 &   1 &  -1 &   0 &   0 &   0 &   0 \\
 0 &   0 &  -1 &   0 &   0 &   1 &   0 &   0 &  -1 &   0 \\
 0 &   0 &   0 &  -1 &   1 &   0 &   1 &  -1 &   0 &   0 \\
 0 &   0 &   0 &   0 &  -1 &   1 &  -1 &   0 &   0 &   0 \\
 0 &   0 &   0 &   0 &   0 &   0 &  -1 &   1 &   0 &  -1 \\
 0 &   0 &   0 &   0 &   0 &   0 &   1 &   0 &   0 &   0 \\
 0 &   0 &   0 &   0 &   0 &   0 &   0 &   0 &   0 &   1 \\
 1 &   0 &   0 &   0 &   0 &   0 &   0 &   0 &   0 &   0 \\
 0 &   0 &   0 &   0 &   0 &   0 &   0 &   0 &   1 &   0
\end{array}\right).$$
Thus $\mathbf{a}\in C_{\mathbf{i}}$ if and only if
$\mathbf{a}=a\mathbf{c}_1+b\mathbf{c}_2+c\mathbf{c}_3+d\mathbf{c}_4+e\mathbf{c}_5+f\mathbf{c}_6+g\mathbf{c}_7
+h\mathbf{c}_8+i\mathbf{c}_9+j\mathbf{c}_{10}$, where $\mathbf{c}_1,\mathbf{c}_2,\ldots \mathbf{c}_{10}$ are the columns of
$Q_{\mathbf{i}}=P_{\mathbf{i}}^{-1}$ and $a\geq 0,\ b\geq 0,\ldots ,j\geq 0$. We
have
$$Q_{\mathbf{i}}=\left( \begin{array}{cccccccccc}
0 &  0 &  0 &  0 &  0 &  0 &  0 &  0 &  1 &  0 \\
1 &  0 &  0 &  1 &  0 &  1 &  0 &  1 &  1 &  0 \\
1 &  1 &  0 &  1 &  1 &  1 &  1 &  1 &  1 &  0 \\
1 &  1 &  1 &  0 &  0 &  0 &  0 &  0 &  1 &  1 \\ 
1 &  1 &  1 &  1 &  0 &  1 &  0 &  1 &  1 &  1 \\
1 &  1 &  1 &  1 &  1 &  1 &  1 &  1 &  1 &  1 \\
0 &  0 &  0 &  0 &  0 &  0 &  1 &  0 &  0 &  0 \\
0 &  0 &  0 &  0 &  0 &  1 &  1 &  1 &  0 &  0 \\
0 &  0 &  0 &  0 &  0 &  0 &  0 &  0 &  0 &  1 \\
0 &  0 &  0 &  0 &  0 &  0 &  0 &  1 &  0 &  0
\end{array}\right).$$

Thus a typical point $\mathbf{a}\in C_{\mathbf{i}}$ is given by
\begin{eqnarray*}
\mathbf{a} & = & (i, a + d + f + h + i, a + b + d + e + f + g + h + i,
a + b + c + i + j, \\
& & a + b + c + d + f + h + i + j,
a + b + c + d + e + f + g + h + i + j, g, f + g + h, j, h),
\end{eqnarray*}
where $a\geq 0,\ b\geq 0,\ldots ,j\geq 0$.

We now consider $S_{\mathbf{i}}^{\mathbf{j}}(\mathbf{a})$. We recall that
$S_{\mathbf{i}}^{\mathbf{j}}=\phi_{\mathbf{j}}\psi_{\mathbf{i}}^{-1}$ where
$\psi_{\mathbf{i}}:\B\rightarrow Y_{\mathbf{i}}$ and $\phi_{\mathbf{j}}:B\rightarrow
\N^k$ are bijections. Let $\psi_{\mathbf{i}}^{-1}(\mathbf{a})=b$ and
$\phi_{\mathbf{i}}^{-1}(b)=\mathbf{c}$.
Then $S_{\mathbf{i}}^{\mathbf{j}}(\mathbf{a})=\mathbf{c}$. We also recall that
$$b\equiv \TF_{i_1}^{a_1}\TF_{i_2}^{a_2}\cdots \TF_{i_k}^{a_k}\cdot 1 \mod
v{\mathcal{L}}',\mbox{\ \ and}$$
$$b\equiv F_{\mathbf{j}}^{\mathbf{c}}\mod v\mathcal{L},$$
where $\mathcal{L}$ is the $\Z[v]$-lattice spanned by $\B$ and $\mathcal{L}'$ is
the $\CA$-lattice spanned by $\B$. Since $\Z[v]\subseteq \CA$ we have
$\mathcal{L}\subseteq \mathcal{L}'$. Thus
$$\TF_{i_1}^{a_1}\TF_{i_2}^{a_2}\cdots \TF_{i_k}^{a_k}\cdot 1\equiv
F_{\mathbf{j}}^{\mathbf{c}}\mod v\mathcal{L}'.$$

We shall describe an algorithm for obtaining the vector $\mathbf{c}$
determined by
this condition. This is based on the fact that, since the reduced expression
$\mathbf{j}$ is adapted to the quiver $Q_1$ in Figure $11$,
(in the sense of~\cite{lusztig2}),
the actions of $\TF_1$ and $\TF_3$ on elements
$F_{\mathbf{j}}^{\mathbf{c}}$, taken$\mod v\mathcal{L}'$, are given by a simple
formula. By~\cite[Cor. 2.5]{lusztig3} we have
\begin{eqnarray*}
\TF_1^xF_{\mathbf{j}}^{\mathbf{c}} & \equiv & F_{\mathbf{j}}^{\mathbf{c}'}\mod v
\mathcal{L}', \\
\TF_3^xF_{\mathbf{j}}^{\mathbf{c}} & \equiv & F_{\mathbf{j}}^{\mathbf{c}''}\mod v
\mathcal{L}',
\end{eqnarray*}
where
\begin{eqnarray*}
\mathbf{c'} & = & \mathbf{c}+(x,0,0,0,0,0,0,0,0,0), \\
\mathbf{c''} & = & \mathbf{c}+(0,x,0,0,0,0,0,0,0,0).
\end{eqnarray*}
Similarly, the reduced expression $\mathbf{j}'$ is adapted to the quiver
$Q_2$ in Figure $11$,
so the actions of $\TF_2$ and $\TF_4$ on elements $F_{\mathbf{j}'}^{\mathbf{c}},$
taken mod $v\mathcal{L}'$, are given by
\begin{eqnarray*}
\TF_2^xF_{\mathbf{j}'}^{\mathbf{c}} & \equiv & F_{\mathbf{j}'}^{\mathbf{c}'}\mod v
\mathcal{L}', \\
\TF_4^xF_{\mathbf{j}'}^{\mathbf{c}} & \equiv & F_{\mathbf{j}'}^{\mathbf{c}''}\mod v
\mathcal{L}',
\end{eqnarray*}
where
\begin{eqnarray*}
\mathbf{c'} & = & \mathbf{c}+(x,0,0,0,0,0,0,0,0,0), \\
\mathbf{c''} & = & \mathbf{c}+(0,x,0,0,0,0,0,0,0,0).
\end{eqnarray*}

\beginpicture

\setcoordinatesystem units <1cm,1cm>             
\setplotarea x from -1 to 6, y from 0 to 4.5      

\put{Figure $11$: Two quivers of type $A_4$.}[c] at 7.4 0.1
\put{$Q_1$}[c] at 4.7 1.0
\put{$Q_2$}[c] at 10.6 1.0

\scriptsize{

\multiput {$\circ$} at 9   2 *3 1 0 /      %

\linethickness=1pt           

\putrule from 9.05 2 to 9.95 2  %
\putrule from 10.05 2 to 10.95 2  
\putrule from 11.05 2 to 11.95 2  %

\put{$1$}   [l] at 8.9 1.75
\put{$2$}   [l] at 9.9 1.75
\put{$3$}   [l] at 10.9 1.75
\put{$4$}   [l] at 11.9 1.75

\setlinear \plot 9.4 2.1 9.5 2 / %
\setlinear \plot 9.4 1.9 9.5 2 / %
\setlinear \plot 10.6 2.1 10.5 2 / %
\setlinear \plot 10.6 1.9 10.5 2 / %
\setlinear \plot 11.4 2.1 11.5 2 / %
\setlinear \plot 11.4 1.9 11.5 2 / 

}

\scriptsize{

\multiput {$\circ$} at 3   2 *3 1 0 /      %

\linethickness=1pt           

\putrule from 3.05 2 to 3.95 2  %
\putrule from 4.05 2 to 4.95 2  
\putrule from 5.05 2 to 5.95 2  %

\put{$1$}   [l] at 2.9 1.75
\put{$2$}   [l] at 3.9 1.75
\put{$3$}   [l] at 4.9 1.75
\put{$4$}   [l] at 5.9 1.75

\setlinear \plot 3.6 2.1 3.5 2 / %
\setlinear \plot 3.6 1.9 3.5 2 / %
\setlinear \plot 4.4 2.1 4.5 2 / %
\setlinear \plot 4.4 1.9 4.5 2 / %
\setlinear \plot 5.6 2.1 5.5 2 / %
\setlinear \plot 5.6 1.9 5.5 2 / 

}

\endpicture
In order to calculate
$\TF_3^{a_1}\TF_2^{a_2}\TF_1^{a_3}\TF_4^{a_4}\TF_3^{a_5}\TF_2^{a_6}\TF_3^{a_7}
\TF_4^{a_8}\TF_1^{a_9}\TF_3^{a_{10}}\cdot 1\mod v\mathcal{L}'$ we shall write
the element being acted on at each stage in the form $F_{\mathbf{j}}^{\mathbf{c}}$
or $F_{\mathbf{j}'}^{\mathbf{c}'}$ as appropriate, in order to be able to use the
above formulae. The elements can be transformed from form
$F_{\mathbf{j}}^{\mathbf{c}}$ to $F_{\mathbf{j}'}^{\mathbf{c}'}$ by using the
transition function $R=R_{\mathbf{j}}^{\mathbf{j'}}$ and from form
$F_{\mathbf{j}'}^{\mathbf{c}'}$ to $F_{\mathbf{j}}^{\mathbf{c}}$ by
$R^{-1}=R_{\mathbf{j}'}^{\mathbf{j}}$. It follows from the description of the
transition function in~\cite{lusztig2} and the form of $\mathbf{j}$ and
$\mathbf{j}'$ that
$R_{\mathbf{j}'}^{\mathbf{j}}=\tau R_{\mathbf{j}}^{\mathbf{j}'}\tau$ where
$\tau=(1\,2)(3\,4)(5\,6)(7\,8)(9\,10)$. Thus $R^{-1}=\tau R \tau$ where
$\tau(a,b,c,d,e,f,g,h,i,j)=(b,a,d,c,f,e,h,g,j,i)$.

We now apply this algorithm to our given example where
$$(a_1,a_2,a_3,a_4,a_5,a_6,a_7,a_8,a_9,a_{10})=
(i,a+d+f+h+i,a+b+d+e+f+g+h+i, a+b+c+i+j,$$
$$\ \ \ \ \ \ a+b+c+d+f+h+i+j,a+b+c+d+e+f+g+h+i+j,g,
f+g+h,j,h).$$
We have
$$\TF_3^{a_{10}}\cdot 1=\TF_3^{a_{10}}\cdot F_{\mathbf{j}}^{\mathbf{0}}\equiv
F_{\mathbf{j}}^{\mathbf{c_1}}\mod v\mathcal{L}',$$
where $\mathbf{c_1}=(0,h,0,0,0,0,0,0,0,0)$. Also,
$$\TF_1^{a_9}F_{\mathbf{j}}^{\mathbf{c_1}}\equiv F_{\mathbf{j}}^{\mathbf{c_2}}\mod
v\mathcal{L}',$$
where $\mathbf{c_2}=(j,h,0,0,0,0,0,0,0,0)$. One can then check that
$\mathbf{c_2}$ lies in region $144$ of $\S5$ and it follows from Table IV that
$F_{\mathbf{j}}^{\mathbf{c_2}}\equiv F_{\mathbf{j}'}^{\mathbf{c_3}}\mod v
\mathcal{L}',$
where $\mathbf{c_3}=(0,0,0,0,0,0,0,0,h,j)$.
Hence $\TF_4^{a_8}F_{\mathbf{j}'}^{\mathbf{c_3}}\equiv F_{\mathbf{j}'}^{\mathbf{c_4}}
\mod v\mathcal{L}',$
where $\mathbf{c_4}=(0,f+g+h,0,0,0,0,0,0,h,j)$.

We now apply $R_{\mathbf{j'}}^{\mathbf{j}}$ to $\mathbf{c_4}$.
Applying $\tau$ to this vector we
obtain $\mathbf{c_5}=(f+g+h,0,0,0,0,0,0,0,j,h)$. We check that this vector
lies in region $136$ of $\S5$. Hence
$R_{\mathbf{j}}^{\mathbf{j'}}(\mathbf{c_5})=\mathbf{c_6}$,
where $\mathbf{c_6}=(0,j,h,0,0,0,0,0,0,f+g)$.
Applying $\tau$ again we obtain
$\mathbf{c_7}=(j,0,0,h,0,0,0,0,f+g,0)$.
Hence $$F_{\mathbf{j}'}^{\mathbf{c_4}}\equiv F_{\mathbf{j}}^{\mathbf{c_7}}\mod
v\mathcal{L}'.$$
Then $$\TF_3^{a_7}F_{\mathbf{j}}^{\mathbf{c_7}}\equiv
F_{\mathbf{j}}^{\mathbf{c_8}}\mod v\mathcal{L}',$$
where $\mathbf{c_8}=(j,g,0,h,0,0,0,0,f+g,0)$.
This vector lies in region $139$ of $\S5$. Hence $$F_{\mathbf{j}}^{\mathbf{c_8}}
\equiv F_{\mathbf{j}'}^{\mathbf{c_9}}\mod v\mathcal{L}',$$
where $\mathbf{c_9}=(0,f+h,0,0,0,0,g,0,h,j)$. Then
$$\TF_2^{a_6}F_{\mathbf{j}'}^{\mathbf{c_9}}\equiv F_{\mathbf{j}'}^{\mathbf{c_{10}}}
\mod v\mathcal{L}',$$
where
$\mathbf{c_{10}}=(a+b+c+d+e+f+g+h+i+j,f+h,0,0,0,0,g,0,h,j)$.

We now apply $R_{\mathbf{j'}}^{\mathbf{j}}$ to $\mathbf{c_{10}}$.
Applying $\tau$ we obtain
$\mathbf{c_{11}}=(f+h,a+b+c+d+e+f+g+h+i+j,0,0,0,0,0,g,j,h)$. This lies in
region $144$ of $\S5$.  Hence
$R_{\mathbf{j}}^{\mathbf{j'}}(\mathbf{c_{11}})=\mathbf{c_{12}}$,
where $\mathbf{c_{12}}=(0,0,0,0,0,g,j,h,a+b+c+d+e+f+i,f+g)$. Applying $\tau$
again we obtain $\mathbf{c_{13}}=(0,0,0,0,g,0,h,j,f+g,a+b+c+d+e+f+i)$. Hence
$$F_{\mathbf{j}'}^{\mathbf{c_{10}}}\equiv F_{\mathbf{j}}^{\mathbf{c_{13}}}\mod
v\mathcal{L}'.$$
Then $\TF_3^{a_5}F_{\mathbf{j}}^{\mathbf{c_{13}}}\equiv F_{\mathbf{j}}^{\mathbf{c_{14}}}\mod v\mathcal{L}',$
where $\mathbf{c_{14}}=(0,a+b+c+d+f+h+i+j,0,0,g,0,h,j,f+g,a+b+c+d+e+f+i)$.
This lies in region $133$ of $\S5$. Hence
$$F_{\mathbf{j}}^{\mathbf{c_{14}}}
\equiv F_{\mathbf{j}'}^{\mathbf{c_{15}}}\mod v\mathcal{L}',$$
where $\mathbf{c_{15}}=(e+g,0,0,f+h,0,a+b+c+d+i+j,g,0,h,j)$. Thus
$$\TF_4^{a_4}F_{\mathbf{j}'}^{\mathbf{c_{15}}}\equiv F_{\mathbf{j}'}^{\mathbf{c_{16}}}
\mod v\mathcal{L}',$$
where
$\mathbf{c_{16}}=(e+g,a+b+c+i+j,0,f+h,0,a+b+c+d+i+j,g,0,h,j)$.

We now apply $R_{\mathbf{j'}}^{\mathbf{j}}$ to $\mathbf{c_{16}}$.
Applying $\tau$ gives
$\mathbf{c_{17}}=(a+b+c+i+j,e+g,f+h,0,a+b+c+d+i+j,0,0,g,j,h)$.
This lies in region $142$ of $\S5$.  Hence
$R_{\mathbf{j}}^{\mathbf{j'}}(\mathbf{c_{17}})=\mathbf{c_{18}}$,
where $\mathbf{c_{18}}=(d+f+h,0,a+b+c+i+j,0,0,g,j,h,a+b+c+d+e+f+i,f+g)$.
Applying $\tau$ again
gives $\mathbf{c_{19}}=(0,d+f+h,0,a+b+c+i+j,g,0,h,j,f+g,a+b+c+d+e+f+i)$. Hence
$$F_{\mathbf{j}'}^{\mathbf{c_{16}}}\equiv F_{\mathbf{j}}^{\mathbf{c_{19}}}\mod
v\mathcal{L}'.$$
We have $\TF_1^{a_3}F_{\mathbf{j}}^{\mathbf{c_{19}}}\equiv F_{\mathbf{j}}^{\mathbf{c_{20}}}\mod v\mathcal{L}',$
where $\mathbf{c_{20}}=(a+b+d+e+f+g+h+i,d+f+h,0,a+b+c+i+j,g,0,h,j,f+g,a+b+c+d+e+f+i)$.
This lies in region $104$. Hence
$$F_{\mathbf{j}}^{\mathbf{c_{20}}}
\equiv F_{\mathbf{j}'}^{\mathbf{c_{21}}}\mod v\mathcal{L}',$$
where $\mathbf{c_{21}}=(0,a+b+c+i+j,e+g,0,f+h,c+j,g,a+b+d+i,h,j)$. Thus
$$\TF_2^{a_2}F_{\mathbf{j}'}^{\mathbf{c_{21}}}\equiv F_{\mathbf{j}'}^{\mathbf{c_{22}}}
\mod v\mathcal{L}',$$
where
$\mathbf{c_{22}}=(a+d+f+h+i,a+b+c+i+j,e+g,0,f+h,c+j,g,a+b+d+i,h,j)$.

We now apply $R_{\mathbf{j'}}^{\mathbf{j}}$ to $\mathbf{c_{22}}$.
Applying $\tau$ gives
$\mathbf{c_{23}}=(a+b+c+i+j,a+d+f+h+i,0,e+g,c+j,f+h,a+b+d+i,g,j,h)$.
This lies in region $49$.  Hence
$R_{\mathbf{j}}^{\mathbf{j'}}(\mathbf{c_{23}})=\mathbf{c_{24}}$,
where $\mathbf{c_{24}}=(0,b+e+g,b+c+j,d+f+h,a+i,g,j,h,a+b+c+d+e+f+i,f+g)$.
Applying $\tau$ again
gives $\mathbf{c_{25}}=(b+e+g,0,d+f+h,b+c+j,g,a+i,h,j,f+g,a+b+c+d+e+f+i)$.
Thus
$$F_{\mathbf{j}'}^{\mathbf{c_{22}}}\equiv F_{\mathbf{j}}^{\mathbf{c_{25}}}\mod
v\mathcal{L}'.$$
Thus $\TF_3^{a_1}F_{\mathbf{j}}^{\mathbf{c_{25}}}\equiv F_{\mathbf{j}}^{\mathbf{c_{26}}}\mod v\mathcal{L}',$
where $\mathbf{c_{26}}=(b+e+g,i,d+f+h,b+c+j,g,a+i,h,j,f+g,a+b+c+d+e+f+i)$.

Thus we finally arrive at the conclusion that
$$\TF_3^{a_1}\TF_2^{a_2}\TF_1^{a_3}\TF_4^{a_4}\TF_3^{a_5}\TF_2^{a_6}
\TF_3^{a_7}\TF_4^{a_8}\TF_1^{a_9}\TF_3^{a_{10}}\cdot 1\equiv
F_{\mathbf{j}}^{\mathbf{c}}\mod v\mathcal{L}',$$ 
where $\mathbf{c}=(b+e+g,i,d+f+h,b+c+j,g,a+i,h,j,f+g,a+b+c+d+e+f+i)$.

This procedure can be carried out for each reduced expression $\mathbf{i}$
for $w_0$, i.e.\,for representatives of the $62$ commutation classes of
such reduced expressions. This was carried out by computer and in each case
the coordinates of $\mathbf{c}$ are seen to be positive linear combinations
of $a,b,c,d,e,f,g,h,i,j$. Thus we have obtained:

\begin{prop} \label{cixi}
Let $\mathbf{i}$ be a reduced expression for $w_0$ in type $A_4$. Then
$S_{\mathbf{i}}^{\mathbf{j}}(C_{\mathbf{i}})=X_{\mathbf{i}}$. Moreover the
map $S_{\mathbf{i}}^{\mathbf{j}}:C_{\mathbf{i}}\rightarrow X_{\mathbf{i}}$ is
linear.
\end{prop}

\noindent {\bf Proof:} We know from Corollary~\ref{spvects} that
$S_{\mathbf{i}}^{\mathbf{j}}$ maps the spanning vectors of $C_{\mathbf{i}}$ to
spanning vectors of $X_{\mathbf{i}}$. Since $S_{\mathbf{i}}^{\mathbf{j}}$
is linear on $C_{\mathbf{i}}$ it follows that
$S_{\mathbf{i}}^{\mathbf{j}}(C_{\mathbf{i}})=X_{\mathbf{i}}$.~$\Box$

\section{An Isomorphism of Graphs}
We shall now show that the bijection given by
$\mathbf{i}\mapsto X_{\mathbf{i}}$ in Proposition~\ref{iso} between
equivalence classes of reduced words for $w_0$ and simplicial regions is an
isomorphism of graphs.

We first introduce a graph structure on the set of equivalence classes of
reduced words for $w_0$. We say that two such classes are adjacent if there
exist reduced words $\mathbf{i},\mathbf{j}$ in these classes such that
$\mathbf{j}$ can be obtained from $\mathbf{i}$ by applying a long braid
relation.

For example, the classes containing the reduced words $(2,3,1,2,1,3,4,3,2,1)$,
\\ $(2,3,1,2,3,4,3,2,1,2)$ are adjacent. For the former is equivalent to
$(2,3,1,2,3,4,3,1,2,1)$, which can be obtained from the latter by applying
the long braid relation $212\rightarrow 121$.

Secondly, we introduce a graph structure on the set of simplicial regions. Two
simplicial regions are said to be adjacent if they have a common wall with
equation $Q=0$ and there is a bijection between their remaining walls
such that corresponding walls have equations $Q_i=0$, $Q_i+\lambda_iQ=0$
for some $\lambda_i\in \R$.

For example, the simplicial regions $28$ and $41$ are adjacent. We see from
Table V that their walls are as follows:
\begin{eqnarray*}
28   &   [AB, C, BE, H]   \|   [CD, ABE]   &   cefh   \\
41   &   [AB, BE, H]   \|   [C, D, ABE]    &   cehj.
\end{eqnarray*}
The common wall is taken as $C$, i.e. $\alpha_C=0$. The bijection between
the remaining walls satisfies $CD\leftrightarrow D$, $f\leftrightarrow j$
and $X\leftrightarrow X$ for all other walls $X$. This bijection is of the
required form since $\alpha_C(v)=f-j$.

We shall obtain an isomorphism of graphs by showing that $\mathbf{i}$,
$\mathbf{j}$ are adjacent if and only if $X_{\mathbf{i}},X_{\mathbf{j}}$
are adjacent. We first need the following property of a pair of reduced
words:

\begin{prop} \label{differ}
Two classes of reduced words $[\mathbf{i}],[\mathbf{i}']$ for $w_0$ differ
by a simple long braid relation if and only if their two families of chamber
sets differ by just one chamber set.
\end{prop}

We illustrate this proposition with the reduced words $(2,3,1,2,3,4,3,1,2,1)$
\\ and $(2,3,1,2,3,4,3,2,1,2)$, which differ by a long braid relation
$121\rightarrow 212$. Their respective chamber diagrams are shown in
Figure $12$,
\begin{figure}[htbp]

\beginpicture

\setcoordinatesystem units <0.5cm,0.6cm>             
\setplotarea x from -2 to 38, y from -2 to 6.5       

\linethickness=0.5pt           

\put{$1$}[c] at -0.3 4 %
\put{$2$}[c] at -0.3 3 %
\put{$3$}[c] at -0.3 2 %
\put{$4$}[c] at -0.3 1 %
\put{$5$}[c] at -0.3 0 %

\put{$2$}[c] at 1.5 -1 %
\put{$3$}[c] at 2.5 -1 %
\put{$1$}[c] at 3.5 -1 %
\put{$2$}[c] at 4.5 -1 %
\put{$3$}[c] at 5.5 -1 %
\put{$4$}[c] at 6.5 -1 %
\put{$3$}[c] at 7.5 -1 %
\put{$1$}[c] at 8.5 -1 %
\put{$2$}[c] at 9.5 -1 %
\put{$1$}[c] at 10.5 -1 %

\setlinear \plot 0 0  6 0 / %
\setlinear \plot 6 0  8 2 / %
\setlinear \plot 8 2  9 2 / %
\setlinear \plot 9 2  11 4 / %
\setlinear \plot 11 4  12 4 / %

\setlinear \plot 0 1  2 1 / %
\setlinear \plot 2 1  3 2 / %
\setlinear \plot 3 2  4 2 / %
\setlinear \plot 4 2  5 3 / %
\setlinear \plot 5 3  8 3 / %
\setlinear \plot 8 3  9 4 / %
\setlinear \plot 9 4  10 4 / %
\setlinear \plot 10 4  11 3 / %
\setlinear \plot 11 3  12 3 / %

\setlinear \plot 0 2  1 2 / %
\setlinear \plot 1 2  2 3 / %
\setlinear \plot 2 3  3 3 / %
\setlinear \plot 3 3  4 4 / %
\setlinear \plot 4 4  8 4 / %
\setlinear \plot 8 4  10 2 / %
\setlinear \plot 10 2 12 2 / %

\setlinear \plot 0 3  1 3 / %
\setlinear \plot 1 3  3 1 / %
\setlinear \plot 3 1  5 1 / %
\setlinear \plot 5 1  6 2 / %
\setlinear \plot 6 2  7 2 / %
\setlinear \plot 7 2  8 1 / %
\setlinear \plot 8 1  12 1 / %

\setlinear \plot 0 4  3 4 / %
\setlinear \plot 3 4  7 0 / %
\setlinear \plot 7 0  12 0 / %

\put{$245$}[c] at 3 2.5 %
\put{$25$}[c] at 4 1.5 %
\put{$1245$}[c] at 6 3.5 %
\put{$125$}[c] at 7 2.5 %
\put{$15$}[c] at 6.5 1.5 %
\put{$1235$}[c] at 9.5 3.5 %


\put{$1$}[c] at 16.7 4 %
\put{$2$}[c] at 16.7 3 %
\put{$3$}[c] at 16.7 2 %
\put{$4$}[c] at 16.7 1 %
\put{$5$}[c] at 16.7 0 %

\put{$2$}[c] at 18.5 -1 %
\put{$3$}[c] at 19.5 -1 %
\put{$1$}[c] at 20.5 -1 %
\put{$2$}[c] at 21.5 -1 %
\put{$3$}[c] at 22.5 -1 %
\put{$4$}[c] at 23.5 -1 %
\put{$3$}[c] at 24.5 -1 %
\put{$2$}[c] at 25.5 -1 %
\put{$1$}[c] at 26.5 -1 %
\put{$2$}[c] at 27.5 -1 %

\setlinear \plot 17 0  23 0 / %
\setlinear \plot 23 0  27 4 / %
\setlinear \plot 27 4  29 4 / %

\setlinear \plot 17 1  19 1 / %
\setlinear \plot 19 1  20 2 / %
\setlinear \plot 20 2  21 2 / %
\setlinear \plot 21 2  22 3 / %
\setlinear \plot 22 3  25 3 / %
\setlinear \plot 25 3  26 2 / %
\setlinear \plot 26 2  27 2 / %
\setlinear \plot 27 2  28 3 / %
\setlinear \plot 28 3  29 3 / %

\setlinear \plot 17 2  18 2 / %
\setlinear \plot 18 2  19 3 / %
\setlinear \plot 19 3  20 3 / %
\setlinear \plot 20 3  21 4 / %
\setlinear \plot 21 4  26 4 / %
\setlinear \plot 26 4  28 2 / %
\setlinear \plot 28 2 29 2 / %

\setlinear \plot 17 3  18 3 / %
\setlinear \plot 18 3  20 1 / %
\setlinear \plot 20 1  22 1 / %
\setlinear \plot 22 1  23 2 / %
\setlinear \plot 23 2  24 2 / %
\setlinear \plot 24 2  25 1 / %
\setlinear \plot 25 1  29 1 / %

\setlinear \plot 17 4  20 4 / %
\setlinear \plot 20 4  24 0 / %
\setlinear \plot 24 0  29 0 / %

\put{$245$}[c] at 20 2.5 %
\put{$25$}[c] at 21 1.5 %
\put{$1245$}[c] at 23 3.5 %
\put{$125$}[c] at 24 2.5 %
\put{$15$}[c] at 23.5 1.5 %
\put{$124$}[c] at 26.5 2.5 %

\put{Figure $12$: Two words differing by a long braid relation.}[c] at 14 -2.5 %
\endpicture

\end{figure}
and their respective chamber sets are:
$$25,\ 15,\ 245,\ 125,\ 1245,\ 1235,$$
and
$$25,\ 15,\ 245,\ 125,\ 1245,\ 124.$$
\noindent {\bf Proof:} First suppose that the reduced words $\mathbf{i},\mathbf{i}'$
differ by a single long braid relation. Consider the crossings in their
respective chamber diagrams. The only difference in the order of the
crossings is in a set of three consecutive crossings. We may assume
\begin{eqnarray*}
\mathbf{i} & = & \ldots k,k+1,k,\ldots \\
\mathbf{i}' & = & \ldots k+1,k,k+1, \ldots
\end{eqnarray*}
where the dotted letters are the same in $\mathbf{i},\mathbf{i}'$.
The corresponding crossings have form:
\begin{eqnarray*}
\mathbf{i}: &  & \ldots (pq),(pr),(qr) \ldots \\
\mathbf{i}': &  & \ldots (qr),(pr),(pq) \ldots
\end{eqnarray*}
where the dotted crossings are the same for $\mathbf{i},\mathbf{i}'$,
and $p<q<r$. Now the three crossings $(pq),(pr),(qr)$ bound a chamber
in each chamber diagram. This chamber has form as shown in Figure $13$.

\begin{figure}[htbp]

\beginpicture

\setcoordinatesystem units <0.8cm,0.6cm>             
\setplotarea x from -5 to 12, y from -2 to 7       

\linethickness=0.5pt           

\put{$p$}[c] at -1.3 6 %
\put{$q$}[c] at -1.3 5 %
\put{$r$}[c] at -1.3 4 %

\setlinear \plot -1 4  0 4 / %
\setlinear \plot 0 4  1 4 / %
\setlinear \plot 1 4  3 6 / %
\setlinear \plot 3 6  4 6 / %

\setlinear \plot -1 5  0 5 / %
\setlinear \plot 0 5  1 6 / %
\setlinear \plot 1 6  2 6 / %
\setlinear \plot 2 6  3 5 / %
\setlinear \plot 3 5  4 5 / %

\setlinear \plot -1 6  0 6 / %
\setlinear \plot 0 6  1 5 / %
\setlinear \plot 1 5  2 4 / %
\setlinear \plot 2 4  4 4 / %

\put{in CD($\mathbf{i}$),}[l] at 7 5 %


\put{$p$}[c] at -1.3 2 %
\put{$q$}[c] at -1.3 1 %
\put{$r$}[c] at -1.3 0 %

\setlinear \plot -1 0  0 0 / %
\setlinear \plot 0 0  2 2 / %
\setlinear \plot 2 2  4 2 / %

\setlinear \plot -1 1  0 1 / %
\setlinear \plot 0 1  1 0 / %
\setlinear \plot 1 0  2 0 / %
\setlinear \plot 2 0  3 1 / %
\setlinear \plot 3 1  4 1 / %

\setlinear \plot -1 2  0 2 / %
\setlinear \plot 0 2  1 2 / %
\setlinear \plot 1 2  3 0 / %
\setlinear \plot 3 0  4 0 / %

\put{in CD($\mathbf{i}'$).}[l] at 7 1 %

\put{Figure $13$: Chamber bounded by $(pq),(pr)$ and $(qr)$.} at 3.5 -1.5

\endpicture

\end{figure}
The chamber sets corresponding to this chamber in \cdi\ and \cdii\ are
different, but all other chamber sets in \cdi, \cdii\ are the same. Thus
the families of chamber sets in \cdi, \cdii\ differ by just one chamber
set.

Now suppose conversely that the chamber sets in \cdi, \cdii\ differ by just
one chamber set. Let $C$ be the chamber in \cdi\ whose chamber set does not
appear in \cdi. We note that each crossing $(ij)$ appears just once in each
chamber diagram and, if $i<j$, the chamber set for the chamber on the left of
$(ij)$ gives the chamber set on the right of $(ij)$ replacing $j$ by $i$.

\noindent (a) Suppose first that chambers $C,C'$ occupy the same positions in their
respective chamber diagrams. Let the left and right hand crossings of $C$
be $(ij),(kl)$, respectively, with $i<j$, $k<l$. Since all other crossings
are the same in the two chamber diagrams, the left and right hand end
crossings of $C'$ must be $(kl),(ij)$, respectively. In \cdi\ string $k$
is above string $j$ at crossing $(ij)$, whereas string $j$ is above string
$k$ at crossing $(kl)$. Thus the crossing $(jk)$ must be above chamber $C$
in \cdi. A similar argument show that crossing $(jk)$ is below chamber
$C'$ in \cdii. Since crossing $(jk)$ appears in the same row in both
chamber diagrams we obtain a contradiction. Thus $C,C'$ cannot occupy the
same positions in \cdi, \cdii\ respectively.

\noindent (b) Consider the sequence of chambers in the row of \cdi\ containing C. This
sequence has form
$$\ldots C_1,C,C_2,\ldots$$
and their chamber sets have form
$$\ldots S(C_1),S(C),S(C_2),\ldots.$$
When $C$ is removed, we know that $S(C_1),S(C_2)$ will be consecutive
chamber sets in the same row of \cdii. Thus $S(C_1),S(C_2)$ differ by a
single number, i.e.\,
$$|S(C_1)\cap S(C_2)|=|S(C_1)|-1=|S(C_2)|-1.$$
This implies that the left and right hand end crossings of $C$ must have
form
\begin{eqnarray*}
(ij),\ (jk) & & \mbox{for}\ i<j<k, \\
\mbox{or\ \ \ \ }(jk),\ (ij) & & \mbox{for}\ i<j<k.
\end{eqnarray*}
So, in passing from $\ldots, C_1,C,C_2,\ldots$ to $\ldots ,C_1,C_2,\ldots,$
crossings $(ij)$ and $(jk)$ are removed and crossing $(ik)$ is added. Since
the set of crossings in \cdi, \cdii\ is the same, the crossing $(ik)$ must be
removed from \cdi\ in adjoining $C'$ and crossings $(ij)$ and $(jk)$ added.
Thus the row of \cdii\ containing $C'$ is the row of \cdi\ containing the
crossing $(ik)$.

\noindent (c) Suppose the left and right end crossings of $C$ are $(ij),(jk)$
respectively with $i<j<k$. Then the crossing $(ik)$ must be below the level
of $C$ in \cdi. Alternatively, if the left and right end crossings of $C$
are $(jk),(ij)$, respectively with $i<j<k$, the crossing $(ik)$ will be
above the level of $C$ in \cdi. These two possibilities are illustrated in
Figure $14$.

\beginpicture
\setcoordinatesystem units <0.8cm,0.6cm>             
\setplotarea x from -2 to 16, y from -3 to 16       
\linethickness=0.5pt           

\put{$i$}[c] at -0.3 12 %
\put{$j$}[c] at -0.3 11 %

\put{$k$}[c] at 6.2 12 %
\put{$j$}[c] at 6.2 11 %

\setlinear \plot 0 12  3.5 8.5 / %

\setlinear \plot 0 11  1 12 / %
\setlinear \plot 5 12  6 11 / %

\setlinear \plot 2.5 8.5  6 12 / %

\put{$(ik)$}[l] at 3.4 9
\put{Case $1$}[c] at 8 14
\put{$C$} at 3 11.5 


\put{$j$}[c] at 9.7 9.5 %
\put{$k$}[c] at 9.7 8.5 %

\put{$j$}[c] at 16.2 9.5 %
\put{$i$}[c] at 16.2 8.5 %

\setlinear \plot 12.5 12  16 8.5 / %

\setlinear \plot 10 9.5  11 8.5 / %
\setlinear \plot 15 8.5  16 9.5 / %

\setlinear \plot 10 8.5  13.5 12 / %

\put{$(ik)$}[l] at 13.4 11.5 %
\put{$C'$} at 13 9


\put{$i$}[c] at 9.7 4 %
\put{$j$}[c] at 9.7 3 %

\put{$k$}[c] at 16.2 4 %
\put{$j$}[c] at 16.2 3 %

\setlinear \plot 10 4  13.5 0.5 / %

\setlinear \plot 10 3  11 4 / %
\setlinear \plot 15 4  16 3 / %

\setlinear \plot 12.5 0.5  16 4 / %

\put{$(ik)$}[l] at 13.4 1
\put{Case $2$}[c] at 8 6
\put{$C'$} at 13 3.5 


\put{$j$}[c] at -0.3 1.5 %
\put{$k$}[c] at -0.3 0.5 %

\put{$j$}[c] at 6.2 1.5 %
\put{$i$}[c] at 6.2 0.5 %

\setlinear \plot 2.5 4  6 0.5 / %

\setlinear \plot 0 1.5  1 0.5 / %
\setlinear \plot 5 0.5  6 1.5 / %

\setlinear \plot 0 0.5  3.5 4 / %

\put{$(ik)$}[l] at 3.4 3.5 %
\put{$C$} at 3 1 
\put{Figure $14$: The possibilities for crossing $(ik)$.} at 8 -1.5
\endpicture

\noindent (d) We next show that crossing $(ik)$ in \cdi\ lies either in the row
immediately below that of $(ij),(jk)$, or the row immediately above.

It will be sufficient to consider Case $1$ and show that $(ik)$ is in the row
below $(ij),(jk)$. A similar argument will work in Case $2$. So suppose
if possible that $(ik)$ is not in the line immediately
below $(ij)$ and $(jk)$. Then there
will be a crossing $(im)$ with $i<m$ between $(ij)$ and $(ik)$. Crossing
$(jm)$ cannot occur between $(ij)$ and $(jk)$ since then $(ij)$ and $(jk)$
could not bound the same chamber. So there must be a crossing $(mk)$
between $(ik)$ and $(jk)$. Also all crossings involving $m$ between $(im)$
and $(mk)$ are on a level below the level of $C$. See Figure $15$.

\beginpicture
\setcoordinatesystem units <0.8cm,0.6cm>             
\setplotarea x from -6 to 16, y from 8.8 to 9       

\linethickness=0.5pt           

\put{$i$}[c] at -0.3 6 %
\put{$j$}[c] at -0.3 5 %

\put{$k$}[c] at 6.2 6 %
\put{$j$}[c] at 6.2 5 %

\setlinear \plot 0 6  3.5 2.5 / %

\setlinear \plot 0 5  2 6 / %
\setlinear \plot 2 6  2.25 5.75 / %
\setlinear \plot 2.25 5.75  2.5 6 / %
\setlinear \plot 2.5 6  2.75 5.75 / %
\setlinear \plot 2.75 5.75  3 6 / %
\setlinear \plot 3 6  	3.25 5.75 / %
\setlinear \plot 3.25 5.75  3.5 6 / %
\setlinear \plot 3.5 6  3.75 5.75 / %
\setlinear \plot 3.75 5.75  4 6 / %
\setlinear \plot 4 6  6 5 / %

\setlinear \plot 2.5 2.5  6 6 / %

\setlinear \plot 1.5 3.5  2.5 4.5 / %
\setlinear \plot 2.5 4.5  2.67 4.25 / %
\setlinear \plot 2.67 4.25  2.83 4.5 / %
\setlinear \plot 2.83 4.5  3 4.25 / %
\setlinear \plot 3 4.25  3.16 4.5 / %
\setlinear \plot 3.16 4.5  3.33 4.25 / %
\setlinear \plot 3.33 4.25  3.5 4.5 / %
\setlinear \plot 3.5 4.5  4.5 3.5 / %

\put{$(im)$}[r] at 1.6 4
\put{$(mk)$}[l] at 4.4 4
\put{$(ik)$}[c] at 3 2.3
\put{Figure $15$: All crossings involving $m$ between $(im)$ and $(mk)$ are
on a level below the level of $C$.} at 3 0.5
\endpicture

Now consider the chamber diagram \cdii. The crossings involving $m$ in \cdii\ 
are in the same positions as in \cdi\ since $m\not\in\{i,j,k\}$. Crossing
$(ik)$ lies to the right of $(im)$ and to the left of $(mk)$ in \cdii.
Now string $i$ remains below string $m$ to the right of $(im)$ and
string $k$ remains below string $m$ to the left of $(mk)$ in \cdii.
Thus crossing $(ik)$ lies below string $m$ in \cdii. Since all crossings
involving $m$ between $(im)$ and $(mk)$ are on a level below that of $C$,
the level of crossing $(ik)$ in \cdii\ is below the level of $C$ in \cdi.
But the level of $(ik)$
in \cdii\ is the same as the level of $C$ in \cdi,
so we have a contradiction.
Thus crossing $(ik)$ in \cdi\ lies on the line below $(ij)$ and $(jk)$.
In Case $2$ $(ik)$ lies on the line above $(ij)$ and $(jk)$.

\noindent (e) We show next there is no crossing involving $j$ between $(ij)$ and $(jk)$
in \cdi. Suppose we are in Case $1$. If there is a crossing $(mj)$ with
$m<j$ between $(ij)$ and $(jk)$ then crossing $(mk)$ would be to the right
of $(mj)$ and the left of $(jk)$. But then $(ij)$ and $(jk)$ could not bound
the same chamber. Similarly if there is a crossing $(jm)$ with $j<m$ between
$(ij)$ and $(jk)$ then crossing $(im)$ would be to the right of $(ij)$ and
to the left of $(jm)$. Then $(ij)$ and $(jk)$ could not bound the same
chamber. A similar argument applies in Case $2$.

\noindent (f) Now we show there is no crossing involving $i$ between $(ij)$ and
$(ik)$ and no crossing involving $k$ between $(ik)$ and $(jk)$. Suppose we
are in Case $1$. Suppose there were a crossing $(im)$, $i<m$, between
$(ij)$ and $(ik)$. Then string $m$ remains below string $j$ and above string
$i$ between $(im)$ and $(jk)$. Since string $i$ is in the row below string
$j$ at the crossing $(ik)$, by (e), we have a contradiction. Similarly there
is no crossing $(mk)$, $m<k$, between $(ik)$ and $(jk)$. A similar argument
holds in Case $2$.

\noindent (g) We now know that the chambers $C$ and $C'$ are as in Figure
$16$.
\begin{center}
\begin{figure}[htbp]
\beginpicture

\setcoordinatesystem units <0.8cm,0.6cm>             
\setplotarea x from -6 to 16, y from -2 to 6.5       

\linethickness=0.5pt           

\put{Case 1:}[c] at -3 3

\put{$i$}[c] at -0.3 6 %
\put{$j$}[c] at -0.3 5 %
\put{$k$}[c] at 2.2 4 %

\put{$k$}[c] at 6.2 6 %
\put{$j$}[c] at 6.2 5 %
\put{$i$}[c] at 3.8 4 %

\setlinear \plot 0 6  1 5 / %
\setlinear \plot 1 5  2.5 5 / %
\setlinear \plot 2.5 5  3.5 4 / %

\setlinear \plot 0 5  1 6 / %
\setlinear \plot 1 6  5 6 / %
\setlinear \plot 5 6  6 5 / %

\setlinear \plot 2.5 4  3.5 5 / %
\setlinear \plot 3.5 5  5 5 / %
\setlinear \plot 5 5  6 6 / %

\put{$C$}[c] at 3 5.5
\put{in \cdi,} at 9 5

\put{$i$}[c] at -0.3 0 %
\put{$j$}[c] at -0.3 1 %
\put{$k$}[c] at 2.2 2 %

\put{$k$}[c] at 6.2 0 %
\put{$j$}[c] at 6.2 1 %
\put{$i$}[c] at 3.8 2 %

\setlinear \plot 0 0  1 1 / %
\setlinear \plot 1 1  2.5 1 / %
\setlinear \plot 2.5 1  3.5 2 / %

\setlinear \plot 0 1  1 0 / %
\setlinear \plot 1 0  5 0 / %
\setlinear \plot 5 0  6 1 / %

\setlinear \plot 2.5 2  3.5 1 / %
\setlinear \plot 3.5 1  5 1 / %
\setlinear \plot 5 1  6 0 / %

\put{$C'$}[c] at 3 0.5
\put{in \cdii.} at 9 1

\endpicture
\beginpicture

\setcoordinatesystem units <0.8cm,0.6cm>             
\setplotarea x from -6 to 16, y from -2 to 6.5       

\linethickness=0.5pt           

\put{Case 2:}[c] at -3 3

\put{$i$}[c] at -0.3 2 %
\put{$j$}[c] at -0.3 1 %
\put{$k$}[c] at 2.2 0 %

\put{$k$}[c] at 6.2 2 %
\put{$j$}[c] at 6.2 1 %
\put{$i$}[c] at 3.8 0 %

\setlinear \plot 0 2  1 1 / %
\setlinear \plot 1 1  2.5 1 / %
\setlinear \plot 2.5 1  3.5 0 / %

\setlinear \plot 0 1  1 2 / %
\setlinear \plot 1 2  5 2 / %
\setlinear \plot 5 2  6 1 / %

\setlinear \plot 2.5 0  3.5 1 / %
\setlinear \plot 3.5 1  5 1 / %
\setlinear \plot 5 1  6 2 / %

\put{$C'$}[c] at 3 1.5
\put{in \cdii.} at 9 1

\put{$i$}[c] at -0.3 4 %
\put{$j$}[c] at -0.3 5 %
\put{$k$}[c] at 2.2 6 %

\put{$k$}[c] at 6.2 4 %
\put{$j$}[c] at 6.2 5 %
\put{$i$}[c] at 3.8 6 %

\setlinear \plot 0 4  1 5 / %
\setlinear \plot 1 5  2.5 5 / %
\setlinear \plot 2.5 5  3.5 6 / %

\setlinear \plot 0 5  1 4 / %
\setlinear \plot 1 4  5 4 / %
\setlinear \plot 5 4  6 5 / %

\setlinear \plot 2.5 6  3.5 5 / %
\setlinear \plot 3.5 5  5 5 / %
\setlinear \plot 5 5  6 4 / %

\put{$C$}[c] at 3 4.5
\put{in \cdi,} at 9 5

\put{Figure $16$: The chambers $C$ and $C'$.}[c] at 3 -1.5
\endpicture
\end{figure}
\end{center}
In Case $1$ we may choose a total order on the crossings in \cdi\ and \cdii\ 
compatible with the diagrams of the form
\begin{eqnarray*}
\ldots (ij),(ik),(jk),\ldots & & \ \ \mbox{in}\ \cdi\\
\ldots (jk),(ik),(ij),\ldots & & \ \ \mbox{in}\ \cdii,
\end{eqnarray*}
where the crossings agree apart from the given triples. The crossings in this
order give rise to reduced words
\begin{eqnarray*}
\cdots & s_ps_{p+1}s_p & \cdots \\
\cdots & s_{p+1}s_ps_{p+1} & \cdots
\end{eqnarray*}
commutation equivalent to $\mathbf{i},\mathbf{i}'$ respectively and
differing by a long braid relation. A similar agreement holds in Case $2$,
so the proposition is proved.

\begin{prop}
Let $\mathbf{i},\mathbf{i}'$ be reduced words for $w_0$. Their commutation
classes $[\mathbf{i}],[\mathbf{i}']$ are adjacent if and only if the
simplicial regions $X_{\mathbf{i}},X_{\mathbf{i}'}$ are adjacent.
\end{prop}

\noindent {\bf Proof:} Suppose $[\mathbf{i}],[\mathbf{i}']$ are adjacent. Then
$\mathbf{i},\mathbf{i}'$ both determine $6$ partial quivers and these two
sets of $6$ have $5$ in common, by Proposition~\ref{differ}.
Let these partial quivers be
\begin{eqnarray*}
\mathcal{P}(\mathbf{i}) & = & \{P_1,P_2,P_3,P_4,P_5,P_6\},\mbox{\ \ \ \ and} \\
\mathcal{P}(\mathbf{i}') & = & \{P_1,P_2,P_3,P_4,P_5,P'_6\}.
\end{eqnarray*}
Then, by Proposition~\ref{iso}, $X_{\mathbf{i}}$ has spanning vectors
$\{v_P,\ P\in \mathcal{P}(\mathbf{i}),\ v_1,v_2,v_3,v_4\}$ and
$X_{\mathbf{i}'}$ has spanning vectors
$\{v_P,\ P\in \mathcal{P}(\mathbf{i}'),\ v_1,v_2,v_3,v_4\}$.
Let $W$ be the common wall of $X_{\mathbf{i}},X_{\mathbf{i}'}$, spanned
by vectors $v_{P_1},v_{P_2},v_{P_3},v_{P_4},v_{P_5}$ and $v_1,v_2,v_3,v_4$.
Let $W_i$, $i=1,2,\ldots ,9$, be the walls of $X_{\mathbf{i}}$ containing all
spanning vectors of $X_{\mathbf{i}}$ except one, where the omitted spanning
vector is not $v_{P_6}$. Let $W'_i$, $i=1,2,\ldots ,9$, be the corresponding
walls of $X_{\mathbf{i}'}$. Let $W$ be given by equation $Q=0$ and
$W_i$ by equation $Q_i=0$. Then we have
$$W_i\cap W_{i'}=W_i\cap W\mbox{\ \ \ \ \ \ \ \ and}$$
$$W_i\cap W=\{\mathbf{x}\,:\,Q(\mathbf{x})=0,\ Q_i(\mathbf{x})=0\}.$$
Since $W_i\cap W\subseteq W'_i$, the corresponding ideals $I(W_i\cap W)$,
$I(W'_i)$ satisfy $I(W'_i)\subseteq I(W_i\cap W)$. Let $W'_i$ have equation
$Q'_i=0$. Then $Q'_i$ lies in the ideal generated by $Q$ and $Q_i$. The only
linear polynomials in this ideal are those of the form $\lambda Q+
\lambda_iQ_i$ with $\lambda,\lambda_i\in \R$. Thus $Q'_i$ has form
$\lambda Q+\lambda_iQ_i$. Now $W'_i$ is not equal to $W$ so $\lambda_i\not=0$.
Without loss of generality we may choose $\lambda_i=1$, so $Q'_i=Q_i+
\lambda Q$ for some $\lambda\in \R$. This shows that the simplicial regions
$X_{\mathbf{i}},X_{\mathbf{i}'}$ are adjacent.

Conversely, suppose $\mathbf{i},\mathbf{i}'$ are reduced words such that
the regions $X_{\mathbf{i}},X_{\mathbf{i}'}$ are adjacent. Thus they
have a common wall $W$ given by equation $Q=0$ and a $1-1$ correspondence
between their remaining walls
$$W_i\leftrightarrow W'_i\ \ \ \ \ \ \ \ i=1,2,\ldots ,9,$$
where $W_i$ has equation $Q_i=0$ and $W'_i$ has equation $Q_i+\lambda_iQ=0$
for $\lambda_i\in \R$. Let $L_i,L'_i$ be the $1$-dimensional subspaces of
$\R^{10}$ given by
$$L_i=\cap_{j,j\not=i}W_j\cap W,\ \ \ \ \ \ L'_i=\cap_{j,j\not=i}W'_j\cap W.$$
Both $L_i$ and $L'_i$ are given by the equations $Q=0$, $Q_j=0$, $j\in
\{1,2,\ldots ,9\}$, $j\not=i$. Thus $L_i=L'_i$. Hence the regions
$X_{\mathbf{i}},X_{\mathbf{i}'}$ have $9$ common spanning vectors, since
a spanning vector has all coordinates $0$ or $1$ so is determined by the
$1$-dimensional subspace containing it. Four of these common spanning vectors
are the vectors $v_1,v_2,v_3,v_4$ of Proposition~\ref{iso}. The remaining
spanning vectors of $X_{\mathbf{i}},X_{\mathbf{i}'}$ have form $v_P$, for
$P\in \mathcal{P}(\mathbf{i})$, $P\in \mathcal{P}(\mathbf{i}')$
respectively. Since distinct partial quivers $P$ give distinct vectors $v_P$
we see that $\mathcal{P}(\mathbf{i})$ and $\mathcal{P}(\mathbf{i}')$ must
have $5$ partial quivers in common. Thus by Proposition~\ref{differ} the
commutation classes $[\mathbf{i}]$ and $[\mathbf{i}']$ of reduced words must
be adjacent.~$\Box$

This proposition shows that the correspondence $[\mathbf{i}]\rightarrow
X_{\mathbf{i}}$ is an isomorphism of graphs. This graph, with $62$ vertices,
is shown in Figure $17$. There is an involution $\iota$ on the set of
reduced expressions for $w_0$, taking a reduced expression $\mathbf{i}=
(i_1,i_2,\ldots ,i_{10})$ to $(5-i_1,5-i_2,\ldots ,5-i_{10})$ (thus applying
the graph automorphism of the Dynkin diagram of type $A_4$).
This induces an involution on the set of commutation classes of reduced
expressions, and thus on the set of simplicial regions. We have numbered the
simplicial regions in such a way that $\iota$ takes region $m$ to region
$63-m$, for $m=1,2,\ldots 62$. In fact, $\iota$ induces an automorphism of the
graph of simplicial regions.

\begin{figure}
\beginpicture

\setcoordinatesystem units <1.2cm,1.3cm>             
\setplotarea x from -2 to 14, y from 0 to 10       

\linethickness=1pt 

\put{$1$} at  3.7      13
\put{$2$} at  -.3      12
\put{$3$} at  1.7      12
\put{$4$} at  3.7      12
\put{$5$} at  5.7      12
\put{$6$} at  7.7      12
\put{$7$} at  -.3      11
\put{$8$} at  1.7      11
\put{$9$} at  3.7      11
\put{$10$} at  5.7      11
\put{$11$} at  7.7      11
\put{$12$} at  -.3      10
\put{$13$} at  1.7      10
\put{$14$} at  3.7      10
\put{$15$} at  5.7      10
\put{$16$} at  7.7      10
\put{$17$} at  -.3      9
\put{$18$} at  .7      8
\put{$19$} at  1.7      9
\put{$20$} at  2.7      8
\put{$21$} at  3.7      9
\put{$22$} at  4.7      8
\put{$23$} at  5.7      9
\put{$24$} at  6.7      8
\put{$25$} at  7.7      9
\put{$26$} at  8.7      8
\put{$27$} at  -.3      7
\put{$28$} at  1.7      7
\put{$29$} at  3.7      7
\put{$30$} at  5.7      7
\put{$31$} at  7.7      7
\put{$32$} at  8.7      6
\put{$33$} at  .7      6
\put{$34$} at  2.7      6
\put{$35$} at  4.7      6
\put{$36$} at  6.7      6
\put{$37$} at  6.7      4
\put{$38$} at  7.7      5
\put{$39$} at  8.7      4
\put{$40$} at  -.3      5
\put{$41$} at  .7      4
\put{$42$} at  1.7      5
\put{$43$} at  2.7      4
\put{$44$} at  3.7      5
\put{$45$} at  4.7      4
\put{$46$} at  5.7      5
\put{$47$} at  7.7      3
\put{$48$} at  -.3      3
\put{$49$} at  1.7      3
\put{$50$} at  3.7      3
\put{$51$} at  5.7      3
\put{$52$} at  7.7      2
\put{$53$} at  -.3      2
\put{$54$} at  1.7      2
\put{$55$} at  3.7      2
\put{$56$} at  5.7      2
\put{$57$} at  7.7      1
\put{$58$} at  -.3      1
\put{$59$} at  1.7      1
\put{$60$} at  3.7      1
\put{$61$} at  5.7      1
\put{$62$} at  3.7      0

\setlinear \plot   4     13      0     12     /  %
\setlinear \plot   4     13      2     12     /  %
\setlinear \plot   4     13      4     12     /  %
\setlinear \plot   4     13      6     12     /  %
\setlinear \plot   4     13      8     12     /  %
\setlinear \plot   0     12      0     11     /  %
\setlinear \plot   0     12      2     11     /  %
\setlinear \plot   2     12      2     11     /  %
\setlinear \plot   2     12      4     11     /  %
\setlinear \plot   4     12      4     11     /  %
\setlinear \plot   4     12      6     11     /  %
\setlinear \plot   6     12      6     11     /  %
\setlinear \plot   6     12      8     11     /  %
\setlinear \plot   8     12      0     11     /  %
\setlinear \plot   8     12      8     11     /  %
\setlinear \plot   0     11      0     10     /  %
\setlinear \plot   2     11      2     10     /  %
\setlinear \plot   4     11      4     10     /  %
\setlinear \plot   6     11      6     10     /  %
\setlinear \plot   8     11      8     10     /  %
\setlinear \plot   0     10      0     9     /  %
\setlinear \plot   0     10      1     8     /  %
\setlinear \plot   2     10      2     9     /  %
\setlinear \plot   2     10      3     8     /  %
\setlinear \plot   4     10      4     9     /  %
\setlinear \plot   4     10      5     8     /  %
\setlinear \plot   6     10      6     9     /  %
\setlinear \plot   6     10      7     8     /  %
\setlinear \plot   8     10      8     9     /  %
\setlinear \plot   8     10      9     8     /  %
\setlinear \plot   0     9      0     7     /  %
\setlinear \plot   0     9      1     6     /  %
\setlinear \plot   1     8      2     7     /  %
\setlinear \plot   1     8      1     6     /  %
\setlinear \plot   2     9      2     7     /  %
\setlinear \plot   2     9      3     6     /  %
\setlinear \plot   3     8      4     7     /  %
\setlinear \plot   3     8      3     6     /  %
\setlinear \plot   4     9      4     7     /  %
\setlinear \plot   4     9      5     6     /  %
\setlinear \plot   5     8      6     7     /  %
\setlinear \plot   5     8      5     6     /  %
\setlinear \plot   6     9      6     7     /  %
\setlinear \plot   6     9      7     6     /  %
\setlinear \plot   7     8      8     7     /  %
\setlinear \plot   7     8      7     6     /  %
\setlinear \plot   8     9      8     7     /  %
\setlinear \plot   8     9      9     6     /  %
\setlinear \plot   9     8      0     7     /  %
\setlinear \plot   9     8      9     6     /  %
\setlinear \plot   0     7      9     4     /  %
\setlinear \plot   0     7      0     5     /  %
\setlinear \plot   2     7      1     4     /  %
\setlinear \plot   2     7      2     5     /  %
\setlinear \plot   4     7      3     4     /  %
\setlinear \plot   4     7      4     5     /  %
\setlinear \plot   6     7      5     4     /  %
\setlinear \plot   6     7      6     5     /  %
\setlinear \plot   8     7      7     4     /  %
\setlinear \plot   8     7      8     5     /  %
\setlinear \plot   9     6      8     5     /  %
\setlinear \plot   9     6      9     4     /  %
\setlinear \plot   1     6      0     5     /  %
\setlinear \plot   1     6      1     4     /  %
\setlinear \plot   3     6      2     5     /  %
\setlinear \plot   3     6      3     4     /  %
\setlinear \plot   5     6      4     5     /  %
\setlinear \plot   5     6      5     4     /  %
\setlinear \plot   7     6      7     4     /  %
\setlinear \plot   7     6      6     5     /  %
\setlinear \plot   7     4      8     3     /  %
\setlinear \plot   8     5      8     3     /  %
\setlinear \plot   9     4      0     3     /  %
\setlinear \plot   0     5      0     3     /  %
\setlinear \plot   1     4      2     3     /  %
\setlinear \plot   2     5      2     3     /  %
\setlinear \plot   3     4      4     3     /  %
\setlinear \plot   4     5      4     3     /  %
\setlinear \plot   5     4      6     3     /  %
\setlinear \plot   6     5      6     3     /  %
\setlinear \plot   8     3      8     2     /  %
\setlinear \plot   0     3      0     2     /  %
\setlinear \plot   2     3      2     2     /  %
\setlinear \plot   4     3      4     2     /  %
\setlinear \plot   6     3      6     2     /  %
\setlinear \plot   8     2      8     1     /  %
\setlinear \plot   8     2      0     1     /  %
\setlinear \plot   0     2      0     1     /  %
\setlinear \plot   0     2      2     1     /  %
\setlinear \plot   2     2      2     1     /  %
\setlinear \plot   2     2      4     1     /  %
\setlinear \plot   4     2      4     1     /  %
\setlinear \plot   4     2      6     1     /  %
\setlinear \plot   6     2      8     1     /  %
\setlinear \plot   6     2      6     1     /  %
\setlinear \plot   8     1      4     0     /  %
\setlinear \plot   0     1      4     0     /  %
\setlinear \plot   2     1      4     0     /  %
\setlinear \plot   4     1      4     0     /  %
\setlinear \plot   6     1      4     0     /  %

\put{The graph of simplicial regions in Type $A_4$}[c] at 4 -1
\put{Figure $17$}[c] at 4 -2

\endpicture

\end{figure}

\section{A Correspondence of Monomials}

We shall now show how monomials in the Kashiwara operators 
given by vectors in a Lusztig cone are related to the corresponding monomials
in divided powers in the generators $F_i$ of $U^-$. To do this
we shall relate both to elements of a basis in $U^-$ of PBW-type.
For this purpose we shall need to consider root vectors of $U^-$. We follow
ideas of Xi in~\cite{xi1}.

Let $\mathbf{l}=(l_1,l_2,\ldots ,l_k)$ be the reduced expression for
$w_0$ given by
$$(l_1,l_2,\ldots ,l_k)=(n,n-1,n,n-2,n-1,n,\ldots ,1,2, \ldots ,n).$$
The elements of the corresponding basis $B_{\mathbf{l}}$ of PBW-type are:
$$F_{\mathbf{l}}^{\mathbf{c}}=F_{l_1}^{(c_1)}
T''_{l_1,-1}(F_{l_2}^{(c_2)})\cdots
T''_{l_1,-1}T''_{l_2,-1}\cdots T''_{l_{k-1},-1}(F_{l_k}^{(c_k)}).$$
The corresponding sequence of positive roots is:
$$\alpha_{l_1},s_{l_1}(\alpha_{l_2}),s_{l_1}s_{l_2}(\alpha_{l_3}),\ldots
,s_{l_1}s_{l_2}\cdots s_{l_{k-1}}(\alpha_{l_k}).$$
We write $\alpha_{ij}=\alpha_i+\cdots +\alpha_{j-1}$ for $i<j$. For each
$p$ with $1\leq p\leq k$ we have
$s_{l_1}s_{l_2}\cdots s_{l_{p-1}}(\alpha_{l_p})=\alpha_{ij}$ for some
$i<j$, and each $\alpha_{ij}$ appears in this way for precisely one $p$.
We then define
$$F_{ij}=T_{l_1,-1}''T_{l_2,-1}''\cdots T_{l_{p-1},-1}''(F_{l_p}).$$
Each $F_{ij}$ is an element of $U^-$ of weight $-\alpha_{ij}$. The elements
$F_{ij}$ will be called the root vectors of $U^-$.

\begin{lemma} \label{pbwexp}
We have:
$$F_{ij}=(-1)^{j-i-1}T_i^{-1}T_{i+1}^{-1}\cdots T_{j-2}^{-1}F_{j-1}$$
for $i<j$.
\end{lemma}

\noindent {\bf Proof:} Let $p$ be the positive integer such that
$s_{l_1}s_{l_2}\cdots s_{l_{p-1}}(\alpha_{l_p})=\alpha_{ij}$. Then we have:
$$s_{l_1}s_{l_2}\cdots s_{l_p}=s_ns_{n-1}s_n\cdots s_{i+1}s_{i+2}\cdots s_n
s_is_{i+1}\cdots s_{n+i+1-j}.$$ It is readily checked that the word
$$s_{j-2}s_{j-3}\cdots s_is_ns_{n-1}s_n\cdots s_{i+1}s_{i+2}\cdots s_ns_i
s_{i+1}\cdots s_{n+i+1-j}$$
is reduced and that this element of $W$ transforms $\alpha_{n+i+1-j}$ into
$\alpha_{j-1}$. It follows from Lusztig~\cite[\S1.3(c)]{lusztig2}, that
$$T''_{j-2,-1}T''_{j-3,-1}\cdots T''_{i,-1}T''_{n,-1}T''_{n-1,-1}T''_{n,-1}
\cdots T''_{n+i+1-j}(F_{n+i+1-j})=F_{j-1}.$$
This asserts that:
$$T''_{j-2,-1}T''_{j-3,-1}\cdots T''_{i,-1}(F_{ij})=F_{j-1},$$
and so
\begin{eqnarray*}
F_{ij} & = & (T''_{i,-1})^{-1}(T''_{i+1,-1})^{-1}\cdots (T''_{j-2,-1})^{-1}(F_{j-1})
\\
 & = & r_i^{-1}T_i^{-1}r_{i+1}^{-1}T_{i+1}^{-1}\cdots r_{j-2}^{-1}T_{j-2}^{-1}
(F_{j-1})
\end{eqnarray*}
since $T''_{i,-1}=T_ir_i$. It follows, using the definition of $r_i$, that
$$F_{ij}=(-1)^{j-i-1}T_i^{-1}T_{i+1}^{-1}\cdots T_{j-2}^{-1}(F_{j-1}),$$
as required.~$\Box$

\begin{lemma}
We have
$$F_{ij}=F_iF_{i+1,j}-vF_{i+1,j}F_i,$$
for $i+1<j$.
\end{lemma}

\noindent {\bf Proof:} First note that
$$F_{i,i+2}=-T_i^{-1}(F_{i+1})=F_iF_{i+1}-vF_{i+1}F_i.$$
Assuming the required result inductively for $F_{i',j'}$ with
$j'-i'<j-i$ we have:
\begin{eqnarray*}
F_{ij} &=& -T_i^{-1}F_{i+1,j} \\
          &=& -T_{i}^{-1}(F_{i+1}F_{i+2,j}-vF_{i+2,j}F_{i+1}) \\
          &=& (F_iF_{i+1}-vF_{i+1}F_i)F_{i+2,j}-vF_{i+2,j}(F_iF_{i+1}-
               vF_{i+1}F_i) \\
          &=& F_i(F_{i+1}F_{i+2,j}-vF_{i+2,j}F_{i+1})-v(F_{i+1}F_{i+2,j}-
               vF_{i+2,j}F_{i+1})F_i \\
          &=& F_iF_{i+1,j}-vF_{i+1,j}F_i.\Box
\end{eqnarray*}

We now obtain commutation relations between root vectors and their divided
powers analagous to those obtained by Xi in~\cite[\S5.6]{xi1}:

\begin{prop} \label{rels}
\begin{eqnarray*}
\mbox{(a)\ \ }
F_{pq}^{(M)}F_{rs}^{(N)} & = &F_{rs}^{(N)}F_{pq}^{(M)},{\rm \ if\ }q<r
{\rm\ or\ }r<p< q<s, \\
\mbox{(b)\ \ }
F_{pq}^{(M)}F_{rs}^{(N)}&=&v^{MN}F_{rs}^{(N)}F_{pq}^{(M)},{\rm \ if\ }
r<p<q=s, \\
\mbox{(c)\ \ }
F_{pq}^{(M)}F_{rs}^{(N)}&=&v^{-MN}F_{rs}^{(N)}F_{pq}^{(M)},{\rm\ if\ }p=r
<q<s{\rm\ or\ }p<r<q=s, \\
\mbox{(d)\ \ }F_{pq}^{(M)}F_{rs}^{(N)}&=&\sum_{0\leq t\leq \min (M,N)}
v^{(M-t)(N-t)}F_{rs}^{(N-t)}F_{ps}^{(t)}F_{pq}^{(M-t)},{\rm\ if\ }q=r,\mbox{\ \ and} \\
\mbox{(e)\ \ }
F_{pq}^{(M)}F_{rs}^{(N)} & = &F_{rs}^{(N)}F_{pq}^{(M)},{\rm \ if\ }s<p
{\rm\ or\ }p<r<s<q. \\
\mbox{(f)\ \ }
F_{pq}^{(M)}F_{rs}^{(N)} & = & \sum_{0\leq t\leq \min (M,N)}
v^{-\frac{1}{2}t(t-1)}(v^{-1}-v)^t[t]! F_{rq}^{(t)}F_{rs}^{(N-t)}
F_{pq}^{(M-t)}F_{ps}^{(t)},{\rm\ if\ }p<r<q<s.
\end{eqnarray*}
\end{prop}

Note that (e) is a restatement of (a).
Using the relations (a)--(e) we may now obtain our main result of this section:

\begin{theorem} \label{kashmon}
Suppose $\mathbf{i}$ is a reduced expression for $w_0$, and suppose ${{\mathbf{a}}}\in
C_{\mathbf{i}}$. Then we have
$$\widetilde{F}_{i_1}^{a_1}\widetilde{F}_{i_2}^{a_2}\cdots \widetilde{F}_{i_k}^{a_k}
\cdot 1\equiv
F_{i_1}^{(a_1)}F_{i_2}^{(a_2)}\cdots F_{i_k}^{(a_k)}
\mod v{{\CL}'}.$$
\end{theorem}

\noindent {\bf Proof}: Let
$b=F_{i_1}^{(a_1)}F_{i_2}^{(a_2)}\cdots F_{i_k}^{(a_k)}\in U^-$.
Then it is shown by Marsh in~\cite[\S4.1]{me7} that in type $A_4$ for
${{\mathbf{a}}}\in C_{\mathbf{i}}$ we have $b\in {\B}$.
So there is
a unique ${\mathbf{c}}\in {\N}^k$ such that $b\equiv
F_{\mathbf{l}}^{\mathbf{c}}\mod v{\CL}$ (with ${\mathbf{l}}$ as
above).

Since $b\in \CL$ and $B_{\mathbf{l}}$ is a ${\Z}[v]$-basis of $\CL$, $b$
is a linear combination of elements $F_{\mathbf{l}}^{\mathbf{d}}$ for
$\mathbf{d}\in\N^k$, with coefficients in $\Z[v]$. We shall obtain information
about this expression for $b$ by means of an algorithm.

This algorithm is as follows. It deals with expressions of the form
$$
\sum_{\mathbf{p}=(p_1,p_2,\ldots ,p_t)\in P}
\lambda_{\mathbf{p},{{\mathbf{a}}}}
F_{\gamma_1}^{(f_1)}F_{\gamma_2}^{(f_2)}
\cdots F_{\gamma_u}^{(f_u)}
$$
where $\gamma_1,\gamma_2,\ldots ,\gamma_u$ are positive roots, $P$ is a
subset of ${\N}^t$ defined by certain linear inequalities,
$\mathbf{a}=(a_1,a_2,\ldots ,a_k)$ is as above,
$f_1,f_2,\ldots ,f_u$ are certain linear functions 
in $a_1,a_2,\ldots ,a_k$ and $p_1,p_2,\ldots ,p_t$, and
$\lambda_{\mathbf{p},{{\mathbf{a}}}}\in \mathbb{Z}[v]$.

The element $F_{i_1}^{(a_1)}F_{i_2}^{(a_2)}\cdots F_{i_k}^{(a_k)}$ is an
expression of this form with just one summand, and the algorithm starts with
this expression.

Let $\alpha^p=s_{l_1}s_{l_2}\cdots s_{l_{p-1}}(\alpha_{l_p})$. Then
$\alpha^1,\alpha^2,\ldots ,\alpha^k$ is the total ordering on the positive
roots determined by $\mathbf{l}$. We write $\alpha^p<\alpha^q$ if
$p<q$. The aim of the algorithm is to transform
$F_{i_1}^{(a_1)}F_{i_2}^{(a_2)}\cdots F_{i_k}^{(a_k)}$ into 
an expression of the above general form in which the positive roots appearing
satisfy $\gamma_1<\gamma_2<\cdots <\gamma_u$.
Consider adjacent pairs
$F_{\gamma_i}^{(f_i)}F_{\gamma_{i+1}}^{(f_{i+1})}$ in each monomial in the
sum, where $\gamma_i
\not=\gamma_{i+1}$. Find the first adjacent pair for which
$\gamma_i>\gamma_{i+1}$. Then, by using one of the relations
(a)--(e) of Proposition~\ref{rels}
it is possible to rewrite
$F_{\gamma_i}^{(f_i)}F_{\gamma_{i+1}}^{(f_{i+1})}$.
In fact just one of these relations can be used. Applying this procedure
to the first adjacent pair in each monomial gives us another expression of
the same general form, so that the algorithm can be repeated.

Whenever an adjacent pair $F_{\gamma_i}^{(f_i)}F_{\gamma_{i+1}}^{(f_{i+1})}$
is obtained with $\gamma_i=\gamma_{i+1}$ it is replaced by the equivalent
expression
\renewcommand{\arraystretch}{0.8}
$\left[ \begin{array}{c} f_i+f_{i+1} \\ f_i \end{array}\right]
F_{\gamma_i}^{(f_i+f_{i+1})}$.
\renewcommand{\arraystretch}{1}

If the algorithm terminates we shall have an expression of $b$ as a
$\Z[v]$-combination of elements in $B_{\mathbf{l}}$. It is not a priori clear
that the algorithm will terminate. However, the algorithm was implemented
in Maple for each reduced expression
$\mathbf{i}$ for $w_0$ in type $A_4$, and did in fact terminate in each case.
The implementation did not calculate the coefficients $\lambda_{\mathbf{p},
\mathbf{a}}$ explicitly, but it did
calculate the smallest power of $v$ appearing in each
$\lambda_{\mathbf{p},{{\mathbf{a}}}}$ (which by
construction must appear with coefficient $1$). This number
is a quadratic expression in $p_1,p_2,\ldots ,p_t,\ a_1,a_2,\ldots ,a_k$.

Next, consider the monomial
$$\widetilde{F}_{i_1}^{a_1}\widetilde{F}_{i_2}^{a_2}\cdots \widetilde{F}_{i_k}^{a_k}
\cdot 1$$ in the Kashiwara root operators $\widetilde{F}_i$. It was shown in
$\S6$ that the map $S_{\mathbf{i}}^{{\mathbf{j}}}$ is
linear on $C_{\mathbf{i}}$, where
$\mathbf{j}=(1,3,2,4,1,3,2,4,1,3)$ and $\mathbf{i}$ is any reduced expression
for $w_0$. In fact, a similar calculation, using Maple, shows that
$S_{\mathbf{i}}^{\mathbf{j}}$ is linear on $C_{\mathbf{i}}$ for all pairs
of reduced expressions $\mathbf{i},\mathbf{j}$. In particular, $S_{\mathbf{i}}^{\mathbf{l}}$ is linear on $C_{\mathbf{i}}$ for all $\mathbf{i}$.
It turns out that $\mathbf{l}$ appears to be a reduced expression for which
the calculation of $S_{\mathbf{i}}^{\mathbf{l}}$ is as simple as possible.
Using Theorem 5.5 in the paper~\cite{me9} by the second author, we have
a description of this map, giving us $S_{\mathbf{i}}^{\mathbf{l}}
(\mathbf{a})=\mathbf{c}$ such that
$$\widetilde{F}_{i_1}^{a_1}\widetilde{F}_{i_2}^{a_2}\cdots \widetilde{F}_{i_k}^{a_k}
\cdot 1\equiv F_{\mathbf{l}}^{\mathbf{c}} \mod v{\CL}'.$$

The coefficient of $F_{\mathbf{l}}^{\mathbf{c}}$ in the expression for
$b$ obtained by the above algorithm was then checked and seen to lie in
$1+v\Z[v]$. It follows from Lusztig~\cite[\S\S2.3, 3.2]{lusztig2}
that the coefficients of all other terms $F_{\mathbf{l}}^{\mathbf{d}}$ in the
sum must lie in $v\Z[v]$, and thus that
$b\equiv F_{{\mathbf{l}}}^{\mathbf{c}}\mod v{{\CL}}$, and so
$b\equiv F_{{\mathbf{l}}}^{\mathbf{c}}\mod v{{\CL}'}$
(by~\cite[\S2.3]{lusztig3}).
It follows that $b=F_{i_1}^{(a_1)}F_{i_2}^{(a_2)}\cdots F_{i_k}^{(a_k)}
\equiv\widetilde{F}_{i_1}^{a_1}\widetilde{F}_{i_2}^{a_2}\cdots \widetilde{F}_{i_k}^{a_k}
\cdot 1 \mod v{{\CL}'}$ as required.~$\square$

We now give an example to demonstrate the above proof.
We take ${\mathbf{i}}=(3, 2, 1, 4, 3, 2, 3, 4, 1, 3)$, and ${\mathbf{a}}\in C_{\mathbf{i}}$.
The corresponding monomial is
$$b=F_3^{(a_1)}F_2^{(a_2)}F_1^{(a_3)}F_4^{(a_4)}F_3^{(a_5)}F_2^{(a_6)}
F_3^{(a_7)}F_4^{(a_8)}F_1^{(a_9)}F_3^{(a_{10})},$$
The ordering on the set of positive
roots induced by $\mathbf{l}=(4,3,4,2,3,4,1,2,3,4)$ is
$$(\alpha_{45},\alpha_{35},\alpha_{34},\alpha_{25},\alpha_{24},\alpha_{23},
\alpha_{15},\alpha_{14},\alpha_{13},\alpha_{12}).$$

We begin the algorithm by writing $b$ in the form
$$b=F_{34}^{(a_1)}F_{23}^{(a_2)}F_{12}^{(a_3)}
F_{45}^{(a_4)}F_{34}^{(a_5)}F_{23}^{(a_6)}
F_{34}^{(a_7)}F_{45}^{(a_8)}F_{12}^{(a_9)}
F_{34}^{(a_{10})},$$
The first adjacent pair of root vectors
appearing in the wrong order is $F_{12}^{(a_3)}F_{45}^{(a_4)}$.
Applying relation (a) of Proposition~\ref{rels} we obtain
$$b=F_{34}^{(a_1)}F_{23}^{(a_2)}F_{45}^{(a_4)}
F_{12}^{(a_3)}F_{34}^{(a_5)}F_{23}^{(a_6)}
F_{34}^{(a_7)}F_{45}^{(a_8)}F_{12}^{(a_9)}
F_{34}^{(a_{10})}.$$
Next apply relation (a) to the adjacent pair
$F_{23}^{(a_2)}F_{45}^{(a_4)}$ and obtain
$$b=F_{34}^{(a_1)}F_{45}^{(a_4)}F_{23}^{(a_2)}
F_{12}^{(a_3)}F_{34}^{(a_5)}F_{23}^{(a_6)}
F_{34}^{(a_7)}F_{45}^{(a_8)}F_{12}^{(a_9)}
F_{34}^{(a_{10})}.$$
Now we apply relation (d) to the adjacent pair
$F_{34}^{(a_1)}F_{45}^{(a_4)}$ and obtain
$$b=\sum
v^{(a_1-p_1)(a_4-p_1)}F_{45}^{(a_4-p_1)}F_{35}^{(p_1)}
F_{34}^{(a_1-p_1)}F_{23}^{(a_2)}
F_{12}^{(a_3)}F_{34}^{(a_5)}F_{23}^{(a_6)}
F_{34}^{(a_7)}F_{45}^{(a_8)}F_{12}^{(a_9)}
F_{34}^{(a_{10})},$$
where the sum is over those $p_1$ satisfying
$0\leq p_1\leq \min (a_1,a_4)$.
Continuing the algorithm in this way we obtain, after 25 steps,

$$b=\sum_{(p_1,p_2,\ldots p_{10})\in P}\lambda_{\mathbf{p},{{\mathbf{a}}}}
F_{45}^{(f_1)}F_{35}^{(f_2)}F_{34}^{(f_3)}
F_{25}^{(f_4)}F_{24}^{(f_5)}F_{23}^{(f_6)}
F_{15}^{(f_7)}F_{14}^{(f_8)}F_{13}^{(f_9)}
F_{12}^{(f_{10})},
$$
where \begin{eqnarray*}
\mathbf{f}& = & (f_1,f_2,f_3,f_4,f_5,f_6,f_7,f_8,f_9,f_{10})=
(a_4 + a_8 - p_1 - p_6 - p_7 - p_8, p_1 + p_8, \\
& & a_1 + a_5 + a_7 + a_{10} - p_1 - p_2 - p_4 - p_5 - p_8 - p_9 - p_{10},
    p_7, p_2 + p_5 - p_7 + p_{10}, \\
& & a_2 + a_6 - p_2 - p_3 - p_5 - p_{10}, p_6, p_4 - p_6 + p_9, p_3 - p_4 - p_9, a_3 + a_9 - p_3).
\end{eqnarray*}

The lowest power of $v$ in $\lambda_{\mathbf{p},{{\mathbf{a}}}}$ is
$\mathbf{x}T\mathbf{x}^t$ where $\mathbf{x}=(a_1,a_2,\ldots ,a_{10},p_1,p_2,\ldots ,p_{10})$ and
$$\scriptsize T=\left( \begin{array}{cccccccccccccccccccc}
0 & 0 & 0 & 1 & -1 & 0 & -1 & 1 & 0 & -1 & -1 & 1 & 0 & 1 & 1 & -1 & -1 & -1 & 1 & 1 \\
0 & 0 & 0 & 0 & 1 & -1 & 1 & 0 & 0 & 1 & 0 & -1 & 1 & -1 & -1 & 0 & 0 & 0 & -1 & -1 \\
0 & 0 & 0 & 0 & 0 & 1 & 0 & 0 & -1 & 0 & 0 & 0 & -1 & 0 & 0 & 0 & 0 & 0 & 0 & 0 \\
0 & 0 & 0 & 0 & 0 & 0 & 0 & -1 & 0 & 0 & 0 & 0 & 0 & 0 & 0 & 1 & 1 & 1 & 0 & 0 \\
0 & 0 & 0 & 0 & 0 & 0 & -1 & 1 & 0 & -1 & 0 & 0 & 0 & 1 & 1 & -1 & -1 & -1 & 1 & 1 \\
0 & 0 & 0 & 0 & 0 & 0 & 1 & 0 & 0 & 1 & 0 & 0 & 0 & -1 & -1 & 0 & 0 & 0 & -1 & -1 \\
0 & 0 & 0 & 0 & 0 & 0 & 0 & 1 & 0 & -1 & 0 & 0 & 0 & 0 & 0 & -1 & -1 & -1 & 1 & 1 \\
0 & 0 & 0 & 0 & 0 & 0 & 0 & 0 & 0 & 0 & 0 & 0 & 0 & 0 & 0 & 0 & 0 & 0 & 0 & 0 \\
0 & 0 & 0 & 0 & 0 & 0 & 0 & 0 & 0 & 0 & 0 & 0 & 0 & 0 & 0 & 0 & 0 & 0 & 0 & 0 \\
0 & 0 & 0 & 0 & 0 & 0 & 0 & 0 & 0 & 0 & 0 & 0 & 0 & 0 & 0 & 0 & 0 & 0 & 0 & 0 \\
0 & 0 & 0 & -1 & 1 & 0 & 1 & -1 & 0 & 1 & 1 & -1 & 0 & -1 & -1 & 1 & 1 & 0 & -1 & -1 \\
0 & 0 & 0 & 0 & -1 & 1 & -1 & 0 & 0 & -1 & 0 & 1 & -1 & 1 & 0 & 0 & 0 & 1 & 1 & 0 \\
0 & 0 & 0 & 0 & 0 & -1 & 0 & 0 & 1 & 0 & 0 & 0 & 1 & 0 & 1 & 0 & 0 & 0 & 0 & 1 \\
0 & 0 & 0 & 0 & 0 & 0 & -1 & 0 & 0 & -1 & 0 & 0 & 0 & 1 & 0 & 0 & 1 & 1 & 0 & -1 \\
0 & 0 & 0 & 0 & 0 & 0 & -1 & 0 & 0 & -1 & 0 & 0 & 0 & 1 & 1 & 0 & 0 & 1 & 1 & 0 \\
0 & 0 & 0 & 0 & 0 & 0 & 0 & -1 & 0 & 1 & 0 & 0 & 0 & 0 & 0 & 1 & 0 & 0 & 0 & 0 \\
0 & 0 & 0 & 0 & 0 & 0 & 0 & -1 & 0 & 1 & 0 & 0 & 0 & 0 & 0 & 1 & 1 & 0 & -1 & 0 \\
0 & 0 & 0 & 0 & 0 & 0 & 0 & -1 & 0 & 1 & 0 & 0 & 0 & 0 & 0 & 1 & 1 & 1 & -1 & -1 \\
0 & 0 & 0 & 0 & 0 & 0 & 0 & 0 & 0 & -1 & 0 & 0 & 0 & 0 & 0 & 0 & 0 & 0 & 1 & 0 \\
0 & 0 & 0 & 0 & 0 & 0 & 0 & 0 & 0 & -1 & 0 & 0 & 0 & 0 & 0 & 0 & 0 & 0 & 1 & 1
\end{array}\right).
$$
Here $P$ is the set of $(p_1,p_2,\ldots ,p_{10})\in {\N}^{10}$ satisfying
$p_1 \leq a_1, p_1 \leq a_4, p_2 \leq a_2, p_2 \leq a_5, p_3 \leq a_3,
p_3 \leq a_6, p_4 \leq p_3, p_4 \leq a_7, p_5\leq p_2 - p_3 + a_2 + a_6,
p_5\leq -p_4 + a_7, p_6\leq p_4, p_6 \leq a_8, p_7\leq p_2 + p_5,
p_7\leq -p_6 + a_8, p_8\leq -p_1 - p_2 - p_4 - p_5 + a_1 + a_5 + a_7,
p_8\leq -p_6 - p_7+ a_8, p_9\leq p_3 - p_4, p_9 \leq a_{10},
p_{10}\leq -p_2 - p_3 - p_5 + a_2 + a_6, p_{10}\leq -p_9 + a_{10}.$

The function $S_{\mathbf{i}}^{{\mathbf{l}}}({{\mathbf{a}}})={\mathbf{c}}$
is given in this example by
$$\mathbf{c}=
(-a_1+a_4, a_1,-a_2+a_5,-a_7+a_8,a_2+a_7-a_8,-a_3+a_6,a_7,a_{10},
a_3-a_7-a_{10},a_9).$$
It can be checked that if $\mathbf{p}$ is taken as
$(a_1, a_2, a_3, a_7, 0, a_7, a_8-a_7, 0, a_{10}, 0)$
then $\mathbf{p}$ satisfies the inequalities as given above and $\mathbf{f}={\mathbf{c}}$.
Furthermore, $\mathbf{x}T\mathbf{x}^t=0$, so that $\lambda_{\mathbf{p},{{\mathbf{a}}}}\in 1+v\Z [v]$.
Since for any $\mathbf{p}$ satisfying the above inequalities we must have
$\lambda_{{{\mathbf{a}}},\mathbf{p}}\in \mathbb{Z}[v]$ and the coefficient of the lowest power
of $v$ in $\lambda_{{{\mathbf{a}}},\mathbf{p}}$ is always $1$, it follows that the coefficient
of $F_{\mathbf{l}}^{\mathbf{c}}$ in the above sum, when terms in the same
PBW basis element
are collected together, is of the form $e+vh$, where $h\in \mathbb{Z}[v]$ and
$e\in {\N}, e\not=0$. Since $b\in {\B}$, we know
by~\cite[\S\S2.3, 3.2]{lusztig2} that $b$ has a unique
expression as a $\mathbb{Z}[v]$-combination of elements $F_{\mathbf{l}}^{\mathbf{d}}$ as $\mathbf{d}$ varies, in which
precisely one of the coefficients lies in $1+v\mathbb{Z}[v]$. It follows that
$e=1$
and $b\equiv F_{{\mathbf{l}}}^{\mathbf{c}}\mod v{{\CL}}$, whence
$b\equiv F_{{\mathbf{l}}}^{\mathbf{c}}\mod v{{\CL}'}$ by~\cite[\S2.3]{lusztig3}.
But we know that
$$\widetilde{F}_{i_1}^{a_1}\widetilde{F}_{i_2}^{a_2}\cdots \widetilde{F}_{i_{10}}^{a_{10}}
\equiv F_{{\mathbf{l}}}^{\mathbf{c}}\mod v{{\CL}'},$$
since $S_{\mathbf{i}}^{{\mathbf{l}}}({{\mathbf{a}}})={\mathbf{c}}$.
It follows that
$$\widetilde{F}_{i_1}^{a_1}\widetilde{F}_{i_2}^{a_2}\cdots \widetilde{F}_{i_{10}}^{a_{10}}
\cdot 1\equiv
F_{i_1}^{(a_1)}F_{i_2}^{(a_2)}\cdots F_{i_{10}}^{(a_{10})}
\mod v{{\CL}'}.\Box$$

As a consequence of Theorem~\ref{kashmon} we may deduce the relationships
between the
set $M_{\mathbf{i}}$ of tight monomials, the Lusztig cone $C_{\mathbf{i}}$,
and the simplicial region of linearity $X_{\mathbf{i}}$ of the piecewise-linear
function $R_{\mathbf{j}}^{\mathbf{j}'}$.

\begin{theorem}
For each reduced expression $\mathbf{i}$ for $w_0$ let $M_{\mathbf{i}}=
\{ F_{i_1}^{(a_1)}F_{i_2}^{(a_2)}\cdots F_{i_k}^{(a_k)}\,:\,\mathbf{a}\in
C_{\mathbf{i}}\}$. Suppose we have type $A_4$. Then $M_{\mathbf{i}}\subseteq
\B$. Moreover, under Kashiwara's parametrization of $\B$ we have
$\psi_{\mathbf{i}}(M_{\mathbf{i}})=C_{\mathbf{i}}$ and under Lusztig's parametrization
of $\B$ we have $\phi_{\mathbf{j}}(M_{\mathbf{i}})=X_{\mathbf{i}}$.
\end{theorem}

\noindent {\bf Proof:} We recall that $M_{\mathbf{i}}\subseteq \B$ was shown in~\cite{me7}.
Theorem~\ref{kashmon} shows that
$\psi_{\mathbf{i}}(M_{\mathbf{i}})=C_{\mathbf{i}}$. Proposition~\ref{cixi}
shows that $S_{\mathbf{i}}^{\mathbf{j}}(C_{\mathbf{i}})=X_{\mathbf{i}}$.
Now $S_{\mathbf{i}}^{\mathbf{j}}=\phi_{\mathbf{j}}\psi_{\mathbf{i}}^{-1}$.
Since $\psi_{\mathbf{i}}^{-1}(C_{\mathbf{i}})=M_{\mathbf{i}}$ it follows that
$\phi_{\mathbf{j}}(M_{\mathbf{i}})=X_{\mathbf{i}}$.~$\Box$

{\bf Acknowledgements}

The research for this paper was supported in part by EPSRC research grant
ref. GR/K55233 (the second author was an EPSRC research
assistant of Professor K. A. Brown at the University of Glasgow, Scotland).
It was typeset in \LaTeX.

\newcommand{\noopsort}[1]{}\newcommand{\singleletter}[1]{#1}

\end{document}